\documentclass[12pt]{article}
\usepackage{amssymb,amsmath,latexsym}
\usepackage{color}%cyan,red;magenta,green;yellow,blue
\allowdisplaybreaks

\begin{document}

\begin{doublespace}

\def\1{{\bf 1}}
\def\ind{{\bf 1}}
\def\nn{\nonumber}
\def\bee{\begin{equation}}
\def\eee{\end{equation}}
\def\sA {{\cal A}} \def\sB {{\cal B}} \def\sC {{\cal C}}
\def\sD {{\cal D}} \def\sE {{\cal E}} \def\sF {{\cal F}}
\def\sG {{\cal G}} \def\sH {{\cal H}} \def\sI {{\cal I}}
\def\sJ {{\cal J}} \def\sK {{\cal K}} \def\sL {{\cal L}}
\def\sM {{\cal M}} \def\sN {{\cal N}} \def\sO {{\cal O}}
\def\sP {{\cal P}} \def\sQ {{\cal Q}} \def\sR {{\cal R}}
\def\sS {{\cal S}} \def\sT {{\cal T}} \def\sU {{\cal U}}
\def\sV {{\cal V}} \def\sW {{\cal W}} \def\sX {{\cal X}}
\def\sY {{\cal Y}} \def\sZ {{\cal Z}}

\def\bA {{\mathbb A}} \def\bB {{\mathbb B}} \def\bC {{\mathbb C}}
\def\bD {{\mathbb D}} \def\bE {{\mathbb E}} \def\bF {{\mathbb F}}
\def\bG {{\mathbb G}} \def\bH {{\mathbb H}} \def\bI {{\mathbb I}}
\def\bJ {{\mathbb J}} \def\bK {{\mathbb K}} \def\bL {{\mathbb L}}
\def\bM {{\mathbb M}} \def\bN {{\mathbb N}} \def\bO {{\mathbb O}}
\def\bP {{\mathbb P}} \def\bQ {{\mathbb Q}} \def\bR {{\mathbb R}}
\def\bS {{\mathbb S}} \def\bT {{\mathbb T}} \def\bU {{\mathbb U}}
\def\bV {{\mathbb V}} \def\bW {{\mathbb W}} \def\bX {{\mathbb X}}
\def\bY {{\mathbb Y}} \def\bZ {{\mathbb Z}}
\def\R {{\mathbb R}} \def\RR {{\mathbb R}} \def\H {{\mathbb H}}
\def\n{{\bf n}} \def\Z {{\mathbb Z}}

\newcommand{\red}{\color{red}}
\newcommand{\normal}{\color{black}}

\newcommand{\expr}[1]{\left( #1 \right)}
\newcommand{\cl}[1]{\overline{#1}}
\newtheorem{thm}{Theorem}[section]
\newtheorem{lemma}[thm]{Lemma}
\newtheorem{defn}[thm]{Definition}
\newtheorem{prop}[thm]{Proposition}
\newtheorem{corollary}[thm]{Corollary}
\newtheorem{remark}[thm]{Remark}
\newtheorem{example}[thm]{Example}
\numberwithin{equation}{section}
\def\ee{\varepsilon}
\def\qed{{\hfill $\Box$ \bigskip}}
\def\NN{{\mathcal N}}
\def\AA{{\mathcal A}}
\def\MM{{\mathcal M}}
\def\BB{{\mathcal B}}
\def\CC{{\mathcal C}}
\def\LL{{\mathcal L}}
\def\DD{{\mathcal D}}
\def\FF{{\mathcal F}}
\def\EE{{\mathcal E}}
\def\QQ{{\mathcal Q}}
\def\SS{{\mathcal S}}
\def\RR{{\mathbb R}}
\def\R{{\mathbb R}}
\def\L{{\bf L}}
\def\K{{\bf K}}
\def\S{{\bf S}}
\def\A{{\bf A}}
\def\E{{\mathbb E}}
\def\F{{\bf F}}
\def\P{{\mathbb P}}
\def\N{{\mathbb N}}
\def\eps{\varepsilon}
\def\wh{\widehat}
\def\wt{\widetilde}
\def\pf{\noindent{\bf Proof.} }
\def\pff{\noindent{\bf Proof} }
\def\cp{\mathrm{Cap}}
\def\UU{\mathcal U}

\title{\Large \bf  On the boundary theory of subordinate killed L\'evy processes }

\author{{\bf Panki Kim}\thanks{This work was supported by the National Research Foundation of
Korea (NRF) grant funded by the Korea government(MSIP) (No. 2016R1E1A1A01941893).
}
\quad {\bf Renming Song\thanks{Research supported in part by a grant from
the Simons Foundation (\#429343, Renming Song).}} \quad and
\quad {\bf Zoran Vondra\v{c}ek}
\thanks{Research supported in part by the Croatian Science Foundation under the project 3526.}
}

\date{}

\maketitle

\begin{abstract}
Let $Z$ be a subordinate Brownian motion in $\R^d$,
$d\ge 2$, via a subordinator with Laplace exponent $\phi$. 
We kill
the process $Z$ upon exiting a bounded open set $D\subset \R^d$ to obtain
the killed process $Z^D$,
and then we subordinate the process $Z^D$ 
by a subordinator with Laplace exponent $\psi$.
The resulting process is denoted by $Y^D$. Both $\phi$ and $\psi$ are assumed to satisfy certain weak scaling conditions at infinity. 

We study the potential theory of  $Y^D$, in particular the boundary theory. 
First, in case that $D$ is a $\kappa$-fat bounded open set, we show that
the Harnack inequality holds. If, in addition, $D$ satisfies the local exterior volume condition, 
then we prove the Carleson estimate.  In case $D$ is a smooth open set and the lower weak scaling index of $\psi$ is strictly larger than $1/2$, we establish the boundary Harnack principle with explicit decay rate near the boundary of $D$. On the other hand, 
when $\psi(\lambda)=\lambda^{\gamma}$ with $\gamma\in (0,1/2]$, we show that the boundary Harnack principle near the boundary of $D$ fails for any 
bounded $C^{1,1}$ open set $D$. 
Our results give the first example where the Carleson estimate holds true, but the boundary Harnack principle does not.

One of the main ingredients in the proofs is 
the sharp two-sided estimates of the Green function of $Y^D$.  
Under an additional condition on $\psi$, we establish sharp 
two-sided estimates of the jumping kernel of $Y^D$ which
exhibit some unexpected boundary behavior.

We also prove a boundary Harnack principle for non-negative functions harmonic in  a smooth open set $E$ strictly contained in $D$, showing that the behavior of $Y^D$ in the interior of $D$ is determined by the composition $\psi\circ \phi$.
\end{abstract}



\section{Introduction}\label{s:intro}

The fractional Laplacian $(-\Delta)^{\alpha}$, $\alpha\in (0,1)$, in $\R^d$, $d\ge 1$, is a well-studied object in various branches of mathematics. There are many definitions of this operator as an operator on 
the Lebesgue spaces or the space of continuous functions. 
A detailed discussion of different definitions of 
the fractional Laplacian and their equivalence are given in
the recent paper \cite{Kwa}. From a probabilistic point of view, 
the fractional Laplacian is the infinitesimal generator of the 
semigroup of the isotropic $2\alpha$-stable process. 
To be more precise, let $X=(X_t, \P_x)$ be an isotropic $2\alpha$-stable process in $\R^d$, $\alpha\in (0,1)$. 
For any 
non-negative (or bounded) 
Borel function $f:\R^d\to \R$ 
and $t\ge 0$, let $Q_t f(x):=\E_x f(X_t)$. Then
\begin{equation}\label{e:intro1}
-(-\Delta)^{\alpha} f=\lim_{t\to 0}\frac{Q_t f-f}{t},
\end{equation}
provided the limit exists (in the appropriate function space).  

Another definition of the $\alpha$-fractional Laplacian with a probabilistic flavor is through Bochner's subordination of semigroups (or Bochner's functional calculus). Let $\Delta$ be the standard Laplacian defined on some function space, $(P_t)_{t\ge 0}$ the corresponding semigroup, and $
\chi(\lambda)=\lambda^{\alpha}$, $\lambda >0$. Then 
\begin{equation}\label{e:intro2}
-(-\Delta)^{\alpha}f=-
\chi(-\Delta)f=\frac{1}{|\Gamma(-\alpha)|}\int_0^{\infty}(P_t f-f)t^{-\alpha-1}\, dt\, .
\end{equation}
The probabilistic interpretation is as follows: Let $W=(W_t, \P_x)$ be a Brownian motion in $\R^d$ and $S=(S_t)_{t\ge 0}$ an independent 
$\alpha$-stable subordinator. Then the subordinate process 
$(W_{S_t})_{t\ge 0}$ is an isotropic $2\alpha$-stable process,
 and \eqref{e:intro1} and \eqref{e:intro2} are exactly the same. 

Once we move from the whole of $\R^d$ to an open subset $D\subset \R^d$, the question of defining a fractional Laplacian in $D$ becomes more delicate. From a probabilistic point of view, there are two obvious choices. The first, and the most common one, is to consider the isotropic 
$2\alpha$-stable process $X$ killed upon exiting $D$. 
More precisely,  
let $\tau^X_D:=\inf\{t>0: X_t\notin D\}$, and set $X^D_t:=X_t$ if $t<\tau^X_D$, and $X_t:=\partial$ if $t\ge \tau^X_D$ (here $\partial$ is an extra point usually called the cemetery). Define the corresponding semigroup $(Q^D_t)_{t\ge 0}$ by $Q_t^D f(x):=\E_x f(X^D_t)=\E_x(f(X_t), t<\tau^X_D)$. 
Let
\begin{equation}\label{e:intro3}
\LL_1 f:=\lim_{t\to 0}\frac{Q_t^D f-f}{t}
\end{equation}
be the infinitesimal generator of $(Q_t^D)_{t\ge 0}$. This is one possible definition of an $\alpha$-fractional Laplacian in $D$. It is usually called the fractional Laplacain in $D$ with zero exterior condition and can be denoted by 
$-(-\Delta)^{\alpha}_{\ |D}$. This definition corresponds to the definition in \eqref{e:intro1}.
 The second possible choice for an $\alpha$-fractional Laplacian in $D$ is to apply Bochner's functional calculus to the Dirichlet Laplacian 
$\Delta_{|D}$.
  To be more precise, let $(P_t^D)_{t\ge 0}$ be the semigroup corresponding to 
 the killed Brownian motion $W^D$, 
 and let 
 $\Delta_{|D}$
  be the infinitesimal generator of this semigroup. Set
\begin{equation}\label{e:intro4}
\LL_0 f:=
-\left(-\Delta_{|D}\right)^{\alpha}f=
-\chi\left(-\Delta_{|D}\right)f
=\frac{1}{|\Gamma(-\alpha)|}\int_0^{\infty}(P^D_t f-f)t^{-\alpha-1}\, dt\, .
\end{equation}
This definition corresponds to the one in \eqref{e:intro2}, but it yields an operator which is different from $\LL_1$. Probabilistically, $\LL_0$ is the infinitesimal generator of the 
semigroup $(\wt{Q}_t^D)_{t\ge 0}$ corresponding to the subordinate process 
$(W^D_{S_t})$ (where $S$ is still an $\alpha$-stable subordinator independent of $W$). The semigroup $(\wt{Q}_t^D)_{t\ge 0}$ is subordinate to the semigroup $(Q_t^D)_{t\ge 0}$ in the sense that $\wt{Q}_t^D f(x)\le Q_t^D f(x)$ for all non-negative $f$ and all $x\in D$.

Note that $\LL_0$ and $\LL_1$ are infinitesimal generators of 
processes that are obtained from the Brownian motion $W$ by
using two operations: killing and subordination. For $\LL_0$ we first kill 
the Brownian motion when it exits $D$ and 
then subordinate it via the subordinator $S$. For $\LL_1$ the order is reversed: we first subordinate $W$ to get $X$, and then kill $X$ when it exits $D$. This interpretation suggests that $\LL_0$ and $\LL_1$ are just two extremal possibilities for an infinite choice of fractional Laplacians in $D$. For example, let 
$\gamma, \delta \in (0,1)$
 be such that 
$\delta \gamma=\alpha$.
 Let $Z_t:=W_{S_t}$ be a subordinate Brownian motion where now $S$ is an independent $\delta$-stable subordinator. Note that $Z$ is an isotropic $2\delta$-stable process. Let $Z^D$ be the process $Z$ killed upon exiting $D$, and let $Y^D_t:=Z^D_{T_t}$ be the subordinate process via an independent $\gamma$-stable subordinator $T$. Denote by $X$ the 
(twicely)  subordinate Brownian motion: 
$X_t=Z_{T_t}=W(S_{T_t})$, and note that $X$ is a $2\alpha$-stable process.
If $(R_t^D)_{t\ge 0}$ denotes the semigroup of $Y^D$, then $\wt{Q}_t^D f(x)\le R_t^D f(x)\le Q_t^D f(x)$ for all non-negative $f$ and all $x\in D$. The infinitesimal generator $\LL$ of the semigroup $(R_t^D)$ can be written as
$$
\LL:=-\big((-\Delta)^{\delta}_{\ |D}\big)^{\gamma}
$$
and since $\delta \gamma=\alpha$, it also has the right to be called an $\alpha$-fractional Laplacian on $D$. 

The purpose of this paper is to study the potential theory of the operators 
as defined in the display above and see how their properties depend on 
$\delta$ and $\gamma$.
When $\delta=1$ and $\gamma \in (0,1)$ (so $\alpha=\gamma$),  $Y^D$ reduces to a subordinate killed Brownian motion. 
This case was recently studied in \cite{KSV16}
where it was shown that the boundary behavior of $Y^D$ is roughly the same as that of the killed Brownian motion in $D$, while in the interior of $D$, $Y^D$ behaves like a $2\alpha$-stable process. 
In the current case, namely $\delta\in (0,1)$, we will show that, in the interior of $D$, the process $Y^D$ still behaves like a $2\alpha$-stable process. Two potential-theoretic justifications of this are (i) the Green function interior estimates given in Proposition 
\ref{p:green-comparable}, and (ii) the scale invariant boundary Harnack principle, Theorem \ref{t:bhp-2}, which 
implies that the boundary behavior of non-negative functions which are harmonic
in an open set $E\subset \overline{E}\subset D$ is the same as for the $2\alpha$-stable process. 

On the other hand, the boundary potential theory of $Y^D$ is much more complicated and depends on the range of $\gamma$. In order to explain the intricacies involved, 
we first recall the statement of the boundary Harnack principle (BHP).
Let $D$ be a bounded open set in $\R^d$. We say that the BHP holds for $Y^D$ if there exists 
$\wh{R}>0$ such that for every $r\in (0, \wh{R}\, ]$, there exists a constant $c_r \ge 1$ such that for every $Q\in \partial D$ and any two non-negative functions $f$ and $g$ defined on $D$ which are harmonic in $D\cap B(Q,r)$ and vanish continuously on $\partial D \cap B(Q,r)$, it holds that
$$
\frac{f(x)}{g(x)}\le c_r \frac{f(y)}{g(y)},\qquad x,y \in D\cap B(Q,r/2)\, .
$$
If the constant $c_r$ above can be chosen independently of $r\in (0,\wh{R}\, ]$, we say that the scale invariant BHP holds. 

Usually, to prove the BHP one first establishes sharp two-sided Green function estimates and the Carleson estimate. This is the road we also take. As our first main result, we prove the Carleson estimate for rather rough open sets, see Theorem \ref{t:carleson} for the precise statement. 
The Green function estimates in 
an arbitrary bounded $C^{1,1}$ open set $D$
are given in Theorem \ref{t:green-function-estimate}. The Green function estimates provide indication for the decay rate of non-negative harmonic functions near the boundary of $D$. 

Our second main result concerns the case  $\gamma\in (1/2, 1)$. In this case we prove a scale invariant BHP with the explicit decay rate, cf.~Theorem \ref{t:main-bhp}, namely 
there exist constants $\wh{R}>0$ and $c\ge 1$,
such that for every $r\in (0,\wh{R}\, ]$, every $Q\in \partial D$ and any non-negative function $f$ defined on $D$ which is harmonic in $D\cap B(Q,r)$ and vanishes continuously on $\partial D \cap B(Q,r)$, it holds that
$$
\frac{f(x)}{\mathrm{dist}(x, \partial D)^{\delta}}\le c \frac{f(y)}{\mathrm{dist}(y, \partial D)^{\delta}},\qquad x,y \in D\cap B(Q,r/2)\, .
$$

Our last main result, and the most interesting one,
is about the case $\gamma \in (0,1/2]$. Here we show that even the non-scale invariant BHP does not hold in any bounded $C^{1,1}$ open set, cf.~Section \ref{s:counterexample}. We accomplish this by constructing a sequence of harmonic functions that, in the limit, have different decay rate at the boundary than that of the Green function. One of the possible reasons for the failure of the BHP is that the jumping kernel of $Y^D$ exhibits some unexpected boundary behavior which indicates a sort of phase transition at $\gamma=1/2$. This boundary behavior is very different from that of the jumping kernel of the
subordinate killed Brownian motion  studied in \cite{KSV16}.

To the best of our knowledge, 
our results give the first example where the Carleson estimate holds, 
but the BHP fails. 
It was shown in \cite[Theorem 4.3]{SV06}
 that when 
$D$ is a bounded $\kappa$-fat
open set,  the Martin boundary and the minimal Martin boundary of $Y^D$ both coincide
with the Euclidean boundary $\partial D$ for all
$\gamma \in (0,1)$.
This makes the failure of the BHP in the case
$\gamma \in (0,1/2]$ somewhat a surprise to us.

In \cite{KS} and \cite{KSV13}, it has been shown that the (non-scale invariant) BHP does not hold in a certain $C^{1,1}$ 
domain (connected open set) $D$ for truncated stable processes and subordinate Brownian motions with Gaussian components respectively. In both cases this was accomplished by showing 
the Carleson estimate fails in $D$. 
For an extensive survey of the literature on BHP, see the introduction of \cite{KSV-POT17}.

The paper is written for much more general subordinators than $\delta$ and  $\gamma$-stable ones.
The subordinators we consider are defined through their Laplace exponents $\phi$ and $\psi$ which are assumed to satisfy certain weak scaling conditions at infinity. The precise setting and the background material are given in Section \ref{s:prelim}. In Section \ref{s:cgs} we show that the Green functions of small sets away from the 
boundary of $D$ with respect to the processes
$Y^D$ and $X$ are comparable. The argument uses 
a non-local Feynman-Kac transform of $X$, 
which was already employed in \cite{KSV16}. 
 The main results of Section \ref{s:hi} are the Harnack inequality for non-negative functions harmonic 
 with respect to $Y^D$, Theorems \ref{uhp} and \ref{HP2}, 
 and the boundary Harnack principle for 
non-negative functions harmonic 
in an open set strictly contained in $D$,
Theorem \ref{t:bhp-2}, in case when $D$ is $\kappa$-fat.
The proof of the latter is only sketched, being similar to the proof of 
the corresponding result in \cite{KSV16}, while the proof of the former is different.  
In Section \ref{s:ce} we establish the Carleson estimate in bounded $\kappa$-fat open sets satisfying the local exterior volume condition
by using a parabolic version of the Carleson estimate for the process $Z$, cf.~Proposition \ref{p:pCeZ}. The main result of Section \ref{s:jkgfe} is  
the sharp two-sided estimates on the Green function of $Y^D$ 
in bounded $C^{1,1}$ open sets, cf.~Theorem \ref{t:green-function-estimate}.
In Section \ref{s:bhi} we prove the boundary Harnack principle with explicit decay rate in  
an arbitrary bounded $C^{1,1}$ open set $D$ for non-negative functions harmonic in a neighborhood of a boundary point of $D$ under the assumption that the lower weak scaling index of $\psi$ is strictly larger than $1/2$. 
In the proofs of the Carleson estimate and the boundary Harnack principle we follow the ideas from \cite{KSV16}, but contrary to that paper we do not use the explicit boundary behavior of the jumping kernel.
In Section \ref{s:j} we address the question of the boundary behavior of the jumping kernel of 
$Y^D$ in bounded $C^{1,1}$ open sets.
Here we need to distinguish between two cases, essentially corresponding to $\gamma$-stable subordinators with $\gamma<1/2$ and $\gamma>1/2$. 
In these two cases 
the jumping kernels exhibit different boundary behavior, and there seems to be some sort of phase transition going 
from one to the other case. 
The main result there is Theorem \ref{t:jumping-function-estimate}, while Example \ref{e:stable} gives the full picture for stable subordinators.
Finally, in the last section we show that, in case $\psi(\lambda)=\lambda^{\gamma}$ with $\gamma\in (0, 1/2]$, the BHP does not hold for any bounded $C^{1,1}$ open set.

In the remainder of this paper, unless indicated otherwise, whenever we mention the boundary Harnack principle, we mean the scale invariant one.

We end this introduction with a few words about notation used throughout the paper. 
For any two positive functions $f$ and $g$ and constant $c\geq 1$,
$f\asymp^c g$ means that $c^{-1}\, g \leq f \leq c\, g$ on their common domain of
definition.
We will use ``$:=$" to denote a
definition, which is read as ``is defined to be".
For $a, b\in \bR$,
$a\wedge b:=\min \{a, b\}$ and $a\vee b:=\max\{a, b\}$.
For any $x\in \R^d$ and $r>0$, we use $B(x, r)$ to denote the
open ball of radius $r$ centered at $x$.
For a Borel set $V$ in $\R^d$, 
$|V|$ denotes  the Lebesgue measure 
 of $V$.  
For any open set $U\subset \R^d$ and $x\in \R^d$, 
we use $\delta_U(x)$ to denote the distance between $x$ and the boundary $\partial U$.
For any process $(X_t)_{t \ge 0}$, we sometimes write $X(t)$ instead of $X_t$ for notational simplicity.
Upper case letters $C$'s without subscripts
denote  strictly positive constants in the statements of results and their values 
may change in each result.
Upper case letters with subscripts $C_i, i=0,1,2,  \dots$, denote constants 
that will be fixed throughout the paper. 
Lower case letters $c$'s without subscripts denote strictly positive
constants  whose values
are unimportant and which  may change from line to line,
while values of lower case letters with subscripts
$c_i, i=0,1,2,  \dots$, are fixed in each proof, 
and the labeling of these constants starts anew in each proof.
$c_i=c_i(a,b,c,\ldots)$, $i=0,1,2,  \dots$, denote  constants depending on $a, b, c, \ldots$.
The dependence on the dimension $d \ge 1$
may not be mentioned explicitly.
For the exit times from a Borel set $U$ of the three processes, 
we use notation
$\tau^Z_U:=\inf\{t>0: Z_t\notin U\}$, 
$\tau^X_U:=\inf\{t>0: X_t\notin U\}$, 
and $\tau_U:=\inf\{t>0:\, Y^D_t\notin U\}$.

%%%%%%%%%%%%%%%%%%%%%%%%%%%%%%%%%%%%%%%%%%%%%%%%%%%%%%%%%%%%%%%%%%
%%%%%%%%%%%%%                                          Preliminaries                             %%%%%%%%%%%%%%%%%%%%%%%%     
%%%%%%%%%%%%%%%%%%%%%%%%%%%%%%%%%%%%%%%%%%%%%%%%%%%%%%%%%%%%%%%%%%

\section{Preliminaries}\label{s:prelim}
Let $W=(W_t, \P_x)_{t\ge 0, x\in \R^d}$  be a Brownian motion in $\R^d$
running twice as fast as the standard Brownian motion,
$d\ge 2$, 
and $S=(S_t)_{t\ge 0}$ an independent subordinator with Laplace exponent $\phi$ and L\'evy measure $\mu$. We assume that $\phi$ is a complete Bernstein function 
satisfying the following weak scaling condition at infinity: 
There exist $a_1, a_2>0$ and $0<\delta_1\le \delta_2 <1$ such that
\begin{equation}\label{e:weak-scaling-phi}
a_1 \left(\frac{R}{r}\right)^{\delta_1}\le \frac{\phi(R)}{\phi(r)} \le a_2 \left(\frac{R}{r}\right)^{\delta_2}\, , \qquad 1<r\le R<\infty\, .
\end{equation}
It is clear that, for any $r_0\in(0, 1)$, \eqref{e:weak-scaling-phi} is still valid for 
$r_0<r\le R<\infty$ with constants $a_1, a_2$ 
depending on $r_0, \delta_1$ and $\delta_2$. 
We will implicitly use this throughout the paper
and will write $a_1(r_0)$ and $a_2(r_0)$ for the corresponding constants.
The same remark also applies to \eqref{e:weak-scaling-psi}, \eqref{e:weak-scaling-psi-phi} and 
\eqref{e:v-nu-asymp} below.
Without loss of generality we also assume that $\phi(1)=1$. 
Note that it follows from the right-hand side inequality in \eqref{e:weak-scaling-phi} that $\phi$ has no drift.

Define $Z=(Z_t)_{t\ge 0}$ by $Z_t:=W(S_t)$. Then $Z$ is an isotropic L\'evy process with characteristic exponent $\xi\mapsto \phi(|\xi|^2)$ and it is called a subordinate Brownian motion. 

Let $T=(T_t)_{t\ge 0}$ be another subordinator, independent of $Z$, with Laplace exponent $\psi$ and L\'evy measure $\nu$. We assume that $\psi$ is also a complete Bernstein function satisfying the following weak scaling condition at infinity: There exist $b_1, b_2>0$ and $0<\gamma_1\le \gamma_2 <1$ such that
\begin{equation}\label{e:weak-scaling-psi}
b_1 \left(\frac{R}{r}\right)^{\gamma_1}\le \frac{\psi(R)}{\psi(r)} \le b_2 \left(\frac{R}{r}\right)^{\gamma_2}\, , \qquad 1<r\le R<\infty\, .
\end{equation}
So $\psi$ also has  no drift.
Without loss of generality we also assume that $\psi(1)=1$.
 Let $X=(X_t)_{t\ge 0}$ be the process obtained by subordinating $Z$ via the subordinator $T$: $X_t:=Z(T_t)$. Since $X_t=W(S_{T_t})$, we see that $X$ is a subordinate Brownian motion via the subordinator $S\circ T$ with Laplace exponent $\psi\circ \phi$. It is straightforward to see that $\psi\circ \phi$ is a 
complete Bernstein function (see \cite[Corollary 7.9(iii)]{SSV})  satisfying the following weak scaling condition at infinity: 
\begin{equation}\label{e:weak-scaling-psi-phi}
b_1 a_1^{\gamma_1}\left(\frac{R}{r}\right)^{\gamma_1 \delta_1}\le \frac{(\psi\circ \phi)(R)}{(\psi\circ \phi)(r)} \le b_2 a_2^{\gamma_2}\left(\frac{R}{r}\right)^{\gamma_2 \delta_2}\, , \qquad 1<r\le R<\infty\, .
\end{equation}

Let $D\subset \R^d$ be an open set. 
We define the killed processes $Z^D$ and $X^D$ in the usual way. 
Note that both $Z^D$ and $X^D$ are killed subordinate Brownian motions, and that $X^D$ may also be regarded as the process obtained by first  subordinating  
the subordinate Brownian motion $Z$ via $T$  and then  killing  it upon exiting $D$.
We define another process by reversing the order of killing and subordination of $Z$: Let $Y^D=(Y^D_t)_{t\ge 0}$ be defined by $Y^D_t:=Z^D(T_t)$. Then $Y^D$ is 
the process obtained by first  killing  $Z$ upon exiting from $D$ and then  subordinating  the killed process via $T$.
 We will use $(Q_t^D)_{t\ge 0}$ to denote the semigroup of $X^D$ and $(R_t^D)_{t\ge 0}$ the semigroup of $Y^D$. 
 It was shown in \cite{SV08} that $Y^D$ can be realized as 
 $X^D$ killed at a terminal time
and that the semigroup $(R^D_t)_{t\ge 0}$ is subordinate to the semigroup $(Q^D_t)_{t\ge 0}$ in the sense that $R^D_t f(x)\le Q^D_t f(x)$ for all Borel $f:D\to [0,\infty)$, all $t\ge 0$ and all $x\in D$.  
As a consequence,
 cf.~\cite[Proposition 4.5.2]{FOT} and \cite[Proposition 3.2]{SV03}, 
we have the following relation between the killing function $\kappa^{Y^D}$ of $Y^D$ and the killing function $\kappa^{X^D}$  of $X^D$:
\begin{equation}\label{e:killing-functions-relation}
\kappa^{Y^D}(x) \ge \kappa^{X^D}(x)\, , \qquad x\in D\, .
\end{equation}

 Let $v(t)$ be the potential density of the subordinator $T$.
Since $\psi$ satisfies \eqref{e:weak-scaling-psi}, 
by \cite[Corollary 2.4 and Proposition 2.5]{KSV14}, there exists $c\ge 1$ such that 
\begin{equation}\label{e:v-nu-asymp}
v(t)\asymp^c \frac{1}{t\psi(t^{-1})}\, , \qquad \nu(t)\asymp^c \frac{\psi(t^{-1})}{t}\, , \qquad 0<t<1\, .
\end{equation}
Thus $\nu(t)$ satisfies the doubling property near zero:
For every $M>0$ there exists $c=c(M)>0$ such that
\begin{equation}\label{e:doubling-nu}
\nu(t)\le c \nu(2t)\, , \qquad 0<t\le M\, .
\end{equation}
Since $\psi$ is a complete Bernstein function, it follows from \cite[Lemma 2.1]{KSV12} that there exists $c>0$ such that
\begin{equation}\label{e:CBmu0}
\nu(t)\le c\nu(t+1), \qquad t\ge 1.
\end{equation}

\medskip
Now we list some auxiliary results. First, for simplicity, we let $\Phi(r):=\frac{1}{\phi(r^{-2})}$. 
By concavity and monotonicty of $\phi$, it is clear that 
\begin{equation}\label{e:phi(lambda-t)}
(1\wedge \lambda)\phi(t)\le \phi(\lambda t)\le (1\vee \lambda) \phi(t)\, , \qquad \lambda, t>0\, .
\end{equation}
As a consequence,
\begin{equation}\label{e:Phi(lambda-t)}
(1\wedge \lambda^2)\Phi(t)\le \Phi(\lambda t)\le (1\vee \lambda^2) \Phi(t)\, , \qquad \lambda, t>0\, .
\end{equation}
Let $\Phi^{-1}$ be the inverse function of $\Phi$. It is shown in \cite[(7.2)]{CKS} that $\Phi^{-1}$ satisfies the following scaling property: For each $T>0$, there exists 
$C(T)\ge 1$ such that
\begin{equation}\label{e:Phi-1-scaling}
C(T)^{-1}\left(\frac{r}{R}\right)^{1/2\delta_1} \le \frac{\Phi^{-1}(r)}{\Phi^{-1}(R)} \le C(T)\left(\frac{r}{R}\right)^{1/2\delta_2} \, ,\quad 0<r\le R\le T\, .
\end{equation}

 Next, since $\psi$, $\phi$ and $\psi\circ \phi$ are complete Bernstein functions, by \cite[Proposition 7.1 and (7.3)]{SSV}, 
\begin{equation}\label{e:t-psi-t-1}
 t/\phi(t), \,  t\psi(t^{-1}),\, t\phi(t^{-1})
 \text{ and } 
t(\psi\circ \phi)(t^{-1})
\quad \text{are  complete Bernstein functions}.
\end{equation} 
Hence, 
\begin{equation}\label{e:td-Phi-decreasing}
t\mapsto t^{-2}\Phi(t)\quad \text{is a decreasing function}, 
\end{equation}
and, since $d \ge 2$, 
\begin{equation}\label{e:psi-Phi-inv-increasing}
t\mapsto t^{d}\psi(\Phi(t)^{-1}) \quad \text{is an increasing function}.
\end{equation}

Let $p(t, x, y)$ be the transition density of $Z$.
 Note that $p(t,x,y)=p(t, |x-y|)$ where 
$$
p(t,r)=\int_0^{\infty}(4\pi s)^{-d/2}\exp\{-r^2/4s\}\P(S_t\in ds)
$$ 
is decreasing in $r$. 
We denote by $p^D(t,x,y)$ the transition density of $Z^D$, and by $(P_t^D)_{t\ge 0}$ the corresponding semigroup.
By the strong Markov property,
 $p^D(t,x,y)$ is given by the formula
\begin{align}
\label{e:pD}
p^D(t,x,y):=
p(t,x,y)-\E_x[p(t-\tau^Z_D, Z(\tau^Z_D), y), \tau^Z_D<t]\, , \quad t>0, \ x,y\in D\, .
\end{align}

Recall that 
$X$ can be regarded as a subordinate Brownian motion 
via $S\circ T$. 
Let $q(t, x, y)$ be the transition density of $X$, and $q^D(t,x,y)$   the transition density of $X^D$.
We will need the upper bound for the Green function $G^X_D$ of the process $X^D$ in case $D$ is bounded. When the process $X$ is transient (which is always true when $d\ge 3$), the required upper bound is the standard upper bound for the Green function of $X$. The lemma below covers also the case $d=2$.
\begin{lemma}\label{l:upper-bound-wo-transience}
Let $D\subset \R^d$ be a bounded open set, $d\ge 2$. There exists a constant $c=c(\text{diam}(D))\ge 1$ such that for all $x,y\in D$, 
\begin{equation}\label{e:GX-upper-estimate}
G^X_D(x,y)\le c \frac1{|x-y|^{d}\psi(\Phi(|x-y|)^{-1})}\, .
\end{equation}
\end{lemma}
\pf 
It follows from \cite[Lemma 48.3]{Sato} that, for any bounded open set $U\subset\R^d$, there
exists $t_1>0$ such that $\sup_{x\in \R^d}\P_x(X_{t_1}\in U)<1$.  Using this and 
the Chapman-Kolmogorov equation, 
one can easily show (see, for instance, the proof of \cite[Lemma 3.7]{CKS-12}) that, 
for every $R>0$, there exist $c_1, c_2>0$ such that
for every $x_0\in \R^d$,
$$
q^{B(x_0, R)} (t,x,y) \le c_2 e^{-c_1 t} \quad \text{ for all } (t,x,y) 
\in (t_1, \infty) \times B(x_0, R) \times B(x_0, R).
$$
Using this and the upper bound of $q(t,x,y)$ for $t \in(0,1]$ in \cite[Theorem 3]{KSz} and \cite[Proposition 3.2]{KSV17},  we immediately get that,  
for every $R, T>0$, 
there exist $c_3, c_4>0$ such that for every $x_0\in \R^d$ and  $x,y \in B(x_0, R)$,
\begin{align}
\label{e:pDIU}
q^{B(x_0, R)} (t,x,y) \le c_4 \begin{cases} 
e^{-c_3 t}& \text{ for all } t \in [T, \infty) \\ \displaystyle
[( \psi \circ\phi)^{-1}(t^{-1})]^{d/2} \wedge \frac{t (\psi \circ\phi) (|x-y|^{-2})} {|x-y|^{d}}& \text{ for all } t 
\in (0, T). 
\end{cases}
\end{align}
Let $T_R=2/(\psi \circ\phi) ((2R)^{-2})$.
Note that using \eqref{e:weak-scaling-psi-phi}, the analog of \eqref{e:Phi-1-scaling} for $\psi\circ \phi$, and the fact $d\ge 2$, by the  change of variable $u=1/t(\psi \circ\phi) (|x-y|^{-2})$ we have that for  $x,y \in B(x_0, R)$,
\begin{align*}
&|x-y|^{d}\psi(\Phi(|x-y|)^{-1}) \int_0^{T_R} [( \psi \circ\phi)^{-1}(t^{-1})]^{d/2} \wedge \frac{t (\psi \circ\phi) (|x-y|^{-2})} {|x-y|^{d}}dt\\
=&\int_1^\infty
  u^{-3}du+
  \int_{1/T_R(\psi \circ\phi) (|x-y|^{-2})}^1
  u^{-2}  \left(\frac{ ( \psi \circ\phi)^{-1}(u(\psi \circ\phi) (|x-y|^{-2}))}
  { ( \psi \circ\phi)^{-1}((\psi \circ\phi) (|x-y|^{-2}))}  \right)^{d/2}du
 \\
  \le & c_5+
c_6 \int_{0}^1
  u^{-2+d/(2\gamma_2 \delta_2)}  
 du =c_7(R)<\infty.
\end{align*}
Using this and \eqref{e:pDIU} we see that
\begin{align*}
&G_{D}^X(x,y) \le G_{B(x_0, \text{diam}(D) )}^X(x,y)\nn\\
=&\int_0^{T_{\text{diam}(D)}}  q^{B(x_0, \text{diam}(D) )} (t,x,y)dt+\int_{T_{\text{diam}(D)}}^{\infty} q^{B(x_0, \text{diam}(D) )} (t,x,y)dt\nn\\
&\le c_8 \int_0^{T_{\text{diam}(D)}}
[( \psi \circ\phi)^{-1}(t^{-1})]^{d/2} \wedge \frac{t (\psi \circ\phi) (|x-y|^{-2})} {|x-y|^{d}} dt
+c_8\int_{T_{\text{diam}(D)}}^{\infty} e^{-c_9 t}dt\nn\\
&\le  c \frac1{|x-y|^{d}\psi(\Phi(|x-y|)^{-1})}\, ,\qquad \text{for all } x,y\in D\, ,
\end{align*}
thus proving \eqref{e:GX-upper-estimate}. \qed

Recall  that  $Y^D_t=Z^D(T_t)$ is the process obtained by subordinating the process
$Z^D$ via the independent subordinator $T$.  
The transition semigroup  $(R^D_t)_{t\ge 0}$ of $Y^D$ admits a transition density given by
\begin{equation}\label{e:transition-YD}
r^D(t,x,y):=\int_0^{\infty}p^D(s,x,y)\P(T_t \in ds)\, .
\end{equation}

Since
the semigroup $(R^D_t)_{t\ge 0}$ is subordinate to the semigroup $(Q^D_t)_{t\ge 0}$, when $D$ is a bounded open set, 
the process $Y^D$ is transient and  admits a Green function
\begin{equation}\label{e:green-YD}
G^{Y^D}(x,y):=\int_0^{\infty}r^D(t,x,y)\, dt=\int_0^{\infty}p^D(t,x,y)v(t)\, dt\, .
\end{equation}
Moreover, for every bounded open set $D$,
\begin{equation}\label{e:greenup}
G^{Y^D}(x,y)\le   G_D^X(x,y),  \quad x,y\in D.  
\end{equation}

Let $J^{Y^D}(x,y)$ be the jumping kernel of $Y^D$ 
given by
\begin{align}
\label{e:JY}
J^{Y^D}(x,y) :=\int_0^{\infty}p^D(t,x,y)\nu(t)\, dt\, ,
\end{align}
and let $J^X(x, y):=j^X(|x-y|)$ be the L\'evy density of $X$ 
given  by
\begin{align}
\label{e:JX}
j^X(|x-y|):=\int^\infty_0p(t, x, y)\nu(t)\, dt.
\end{align}
Clearly,
\begin{align}
\label{e:JYDJX}
J^{Y^D}(x,y) \le j^X(|x-y|),   \quad x, y \in D.
\end{align}
 Furthermore, by
\cite[Lemma 3.2]{KSV14}, for any $M>0$ there exists $c(M)\ge 1$ such that
\begin{equation}\label{e:JD-estimate}
c(M)^{-1} \frac{\psi(\Phi(|x-y|)^{-1})}{|x-y|^{d}}\le 
j^X(|x-y|)\le c(M) \frac{\psi(\Phi(|x-y|)^{-1})}{|x-y|^{d}}\, ,\qquad |x-y|\le M\, .
\end{equation}

Recall that, for any Borel set $U\subset D$, we use the notation $\tau_U=\inf\{t>0:\, Y^D_t\notin U\}$ for the exit time from $U$ of $Y^D$. 
 
 Since $Y^D$ can be realized as $X^D$ killed at a terminal time, it follows from \cite{Sztonyk} that if $U\subset D$ is a Lipschitz open set, then 
 \begin{equation}\label{e:exit-time-YD}
\P_x(Y^D_{\tau_U}\in \partial U)=0.
\end{equation}
 We will use $\zeta$ to denote the lifetime of $Y^D$.
Then it follows from \cite[Corollary 4.2(i)]{SV08} that the process $Y^D$ dies inside $D$ almost surely, 
i.e., 
\begin{equation}\label{e:YD-killed-inside}
\P_x(Y^D_{\zeta-}\in D)=1\, \qquad \text{for all }x\in D\, .
\end{equation}

For any open set $U\subset D$,  let $Y^{D,U}$ be  the subprocess of $Y^D$ killed upon exiting $U$ and
\begin{equation}\label{e:hkos1} 
r^{D,U}(t,x,y)  :=  
r^D(t,x,y) -  \E_x [ r^D(t - \tau_U,Y^D_{\tau_U},y) : \tau_U < t]\quad  
t>0,\  x,y \in U.
\end{equation}
By the strong Markov property, $r^{D, U}(t,x,y)$ is the transition density of $Y^{D, U}$.
Then the Green function of $Y^{D, U}$ is given by $G^{Y^D}_U(x,y):= \int_0^\infty  r^{D,U}(t, x,y) dt$.
Further, let $Y^U_t:=Z^U(T_t)$ be the process obtained by subordinating the killed process $Z^U$ via the subordinator $T$. The Green function of $Y^U$ will be denoted by  $G^{Y^U}(x,y)$, $x,y\in U$. 
Since the semigroup of $Y^U$ is subordinate to the semigroup of 
$Y^{D,U}$, cf. \cite[Proposition 3.1]{SV08}, we have 
\begin{equation}\label{e:UB-smaller-UDB}
G^{Y^U}(x,y)\le G^{Y^D}_U(x,y)\le G^{Y^D}(x,y)\, \qquad \text{for all }x,y\in U\, .
\end{equation}

Note that one can follow the argument in \cite[Section 2]{KSV16} and see that $(R^D_t)_{t\ge 0}$  satisfies the strong Feller property.

%%%%%%%%%%%%%%%%%%%%%%%%%%%%%%%%%%%%%%%%%%%%%%%%%%%%%%%%
%%%%%%%%%                           Harnack inequality                               %%%%%%%%%%%%%%%%%%%%
%%%%%%%%%%%%%%%%%%%%%%%%%%%%%%%%%%%%%%%%%%%%%%%%%%%%%%%%

\section{Comparability of Green  functions of $Y^D$ and $X$  on small open sets}\label{s:cgs}

In the remainder of the paper, we always assume that  $d\ge 2$, 
\eqref{e:weak-scaling-phi} and \eqref{e:weak-scaling-psi} hold true, 
and $D$ is a bounded open set in $\R^d$.   
 \begin{defn}\label{d:C11} Let $U\subset \R^d$ be an open set and let $Q\in \partial U$. We say that $U$ is $C^{1,1}$ near $Q$ if there exist 
a localization radius $R>0$, 
a $C^{1,1}$-function 
$\varphi_Q:\R^{d-1}\to \R$ satisfying $\varphi_Q(0)=0$, $\nabla \varphi_Q(0)=(0,\dots, 0)$, $\| \nabla \varphi_Q \|_{\infty}\le \Lambda$, $|\nabla \varphi_Q(z)-\nabla \varphi_Q(w)|\le \Lambda |z-w|$, and an orthonormal coordinate system $CS_Q$ with its origin at $Q$ such that
$$
B(Q,R)\cap U=\{y=(\wt{y},y_d)\in B(0,R) \textrm{ in } CS_Q:\, y_d>\varphi_Q(\wt{y})\}\, ,
$$
where $\wt{y}:= (y_1, \dots, y_{d-1})$.
The pair $(R,\Lambda)$ is called the $C^{1, 1}$ 
characteristics of $U$ near $Q$.
An open set $U\subset \R^d$ is said to be a (uniform) $C^{1,1}$ open set with characteristics $(R,\Lambda)$ if it is $C^{1,1}$ with characteristics $(R,\Lambda)$ near every boundary point $Q\in \partial U$.
\end{defn}

Recall that, for any open set  $U\subset D$, $X^U$ (respectively, $Y^{D,U}$)
is the process $X$ (respectively, $Y^D$)  killed upon exiting $U$.
The approach and proofs in this section  are  very similar to  
that of \cite[Section 7]{KSV16}; we first show that  
  when $U$ is  a  relatively compact  open  subset of $D$,
 the process $Y^{D,U}$ can be thought of as a non-local Feynman-Kac transform of $X^U$. 
  Second, if $U$ is a certain $C^{1,1}$  open set,  
 the Green functions of $X^U$ and $Y^{D,U}$ are comparable because
the conditional gauge function related to this transform is bounded 
between two  positive constants. 
We do not give 
full proofs of the results which are almost identical to 
the ones in \cite[Section 7]{KSV16}. 
For the convenience of our readers
we will provide the exact references for the proofs that we omit.

\medskip

We denote by $(P_t)_{t\ge 0}$ the semigroup of $Z$. Let $(\EE^{X^U}, \DD(\EE^{X^U}))$ be the Dirichlet form of $X^U$. Then, cf.~\cite[Section 13.4]{SSV},
\begin{eqnarray}
\EE^{X^U}(f,f)&=&\int_0^{\infty}\int_U f(x)(f(x)-P_s f(x))\, dx\,  \nu(s)ds\, , \label{e:df-of-XU}\\
\DD(\EE^{X^U})&=&\{f\in L^2(U, dx):\, \EE^{X^U}(f,f)<\infty\}\, .\label{e:dom-df-of-XU}
\end{eqnarray} 
Furthermore, for $f\in \DD(\EE^{X^U})$,
\begin{equation}\label{e:df-of-XU-alt}
\EE^{X^U}(f,f)=\frac12 \int_U \int_U (f(x)-f(y))^2 J^X(x,y)dy dx+\int_U f(x)^2 \kappa^X_U(x)dx,
\end{equation}
where $J^X$ is defined in \eqref{e:JX} and
$$
\kappa^X_U(x):=\int_{\R^d\setminus U}J^X(x,y)dy\, .
$$

Recall that  $(P_t^D)_{t\ge 0}$  is the semigroup of $Z^D$.
The Dirichlet form $(\EE^{Y^D}, \DD(\EE^{Y^D}))$ of $Y^D$ is given by
$$
\EE^{Y^D}(f,f)=\int_0^{\infty}\int_D f(x)(f(x)-P_s^D f(x))\, dx\,  \nu(s) ds
$$
and $\DD(\EE^{Y^D})=\{f\in L^2(D,dx):\, \EE^{Y^D}(f,f)<\infty\}$. Moreover, for $f\in \DD(\EE^{Y^D})$, 
$$
\EE^{Y^D}(f,f)=\frac12 \int_D \int_D (f(x)-f(y))^2 J^{Y^D}(x,y)dydx+\int_D f(x)^2 \kappa^{Y^D}(x) dx,
$$
where $J^{Y^D}$ is defined in \eqref{e:JY} and
\begin{align}
\label{e:kms}
\kappa^{Y^D}(x):=\int_0^{\infty}(1-P_t^D 1(x))\nu(t)\, dt\, .
\end{align}
Hence, it follows that the Dirichlet form $(\EE^{Y^{D,U}}, \DD(\EE^{Y^{D,U}}))$ of $Y^{D,U}$ 
is given by
\begin{eqnarray}
\EE^{Y^{D,U}}(f,f)&=&\int_0^{\infty}\int_U f(x)(f(x)-P^D_s f(x))\, dx\,  \nu(s)ds\, , \label{e:df-of-YDU}\\
\DD(\EE^{Y^{D,U}})&=&\{f\in L^2(U, dx):\, \EE^{Y^{D,U}}(f,f)<\infty\}\, .\label{e:dom-df-of-YDU}
\end{eqnarray}

We will need the following simple result.
Recall that $\delta_U(x)$ is the distance between $x$ and the boundary $\partial U$.
\begin{lemma}\label{l:J-difference}
For $x,y\in D$,
\begin{equation}\label{e:J-difference}
J^X(x,y)-J^{Y^D}(x,y)\le j^X(\delta_D(y))\, .
\end{equation}
\end{lemma}
\pf
By \eqref{e:pD},  \eqref{e:JY} and \eqref{e:JX}, 
we have that  for $x,y\in D$,
\begin{eqnarray}\label{e:J-difference0}
J^X(x,y)-J^{Y^D}(x,y)&=&\int_0^{\infty} \E_x[p(t-\tau^Z_D, Z(\tau^Z_D), y), \tau^Z_D<t]\nu(t)\, dt\nn\\
&=&\E_x\left[\int_{\tau^Z_D}^{\infty}p(t-\tau^Z_D, Z(\tau^Z_D), y)\nu(t)dt\right]\, .
\end{eqnarray}
Since, $p(s,z-y) \le p(s, \delta_D(y))$ for every $s>0$ and $z\in  D^c$,  we have that for every $s<t$ and 
$(y,z) \in D \times D^c$,
\begin{eqnarray*}
\int_s^{\infty}p(t-s,z-y)\nu(t)dt& \le &\int_0^{\infty}p( u , \delta_D(y)) \nu(u+s)du\\
&\le &\int_0^{\infty}p( u , \delta_D(y)) \nu(u)du=j^X(\delta_D(y)).
\end{eqnarray*} 
This and \eqref{e:J-difference0} imply the lemma.
\qed

Using Lemma \ref{l:J-difference}, the proof of the next result is the same as 
 that of  \cite[Lemma 7.2]{KSV16}.
\begin{lemma}\label{l:domains-equal}
Let $U$ be a relatively compact  open subset of $D$.
Then $\DD(\EE^{X^U})=\DD(\EE^{Y^{D,U}})$.
\end{lemma}

For $x \in D$, let 
$$
q_U(x):=
\int_U(J^X(x,y)-J^{Y^D}(x,y))dy.
$$
For $x,y\in D$, $x\ne y$, let 
$$
F(x,y):=\frac{J^{Y^D}(x,y)}{J^X(x,y)}-1=\frac{J^{Y^D}(x,y)-J^X(x,y)}{J^X(x,y)}\, ,
$$
and $F(x,x):=0$.  We also define $F(x,\partial):=0$, 
where $\partial$ is the cemetery point.
Then $-1<F(x,y)\le 0$.  Note that
$$
\int_U F(x,y) J^X(x,y)dy=-q_U(x)\, .
$$

Using Lemma \ref{l:J-difference} and \eqref{e:JD-estimate}, the proof of the next result 
is very similar to (and even simpler than) that of \cite[Lemma 7.3]{KSV16}. So we omit the proof.
\begin{lemma}\label{l:estimate-of-F}
There exists $b=b(\phi,d)>2$  such that for all $x_0\in D$ and all $r\in (0,1/b]$ satisfying $B(x_0, (b+1)r)\subset D$, it holds that
$$
\sup_{x,y\in B(x_0,r)}|F(x,y)| \le \frac12\, .
$$
\end{lemma}

Let $b>2$ be the constant from Lemma \ref{l:estimate-of-F}. For $r<1/b$,  
let $U\subset D$ be an open set such that $\mathrm{diam}(U)\le r$ and $\mathrm{dist}(U,\partial D)\ge (b+2)r$. 
Then there exists a ball $B(x_0,r)$ such that $U\subset B(x_0,r)$ and $B(x_0, (b+1)r)\subset D$. Since by Lemma \ref{l:estimate-of-F},
\begin{align}
\label{e:Fbound}
|F(x,y)|\le 1/2 \qquad \text{ for all }x,y\in U,
\end{align}
using the following non-local multiplicative functional
\begin{align*}
K_t^U
&:=\exp\left(\sum_{0<s\le t}\log(1+F(X_{s-}^U,X_s^U))\right)\, ,
\end{align*}
we define 
$$
T^U_tf(x):=\E_x[K^U_t f(X_t^U)]\, .
$$
By \cite[(4.5) and Theorem 4.8]{CS03a} and the proof of \cite[Lemma 7.4]{KSV16}, $(T^U_t)_{t\ge0}$ is a strongly continuous semigroup on $L^2(U,dx)$ with associated quadratic form $(\EE^{Y^{D,U}}, \DD(\EE^{Y^{D,U}}))$.

Let $\E_x^y$ be the expectation with respect to the conditional probability $\P_x^y$ defined via Doob's $h$-transform, with 
$h(\cdot)=G^X_U(\cdot, y)$, starting from $x\in D$.
For $x,y\in U$, $x\neq y$, we define  the conditional gauge function for $K^U_t$
$$
u^U(x,y):=\E_x^y [K^U_{\tau_U^X}],
$$
which is less than or equal to $1$ because $F\le 0$. 
By \cite[Lemma 3.9]{che02}, we have that
\begin{equation}\label{e:gf-XU-YDU-comparable}
G^{Y^D}_U(x,y)=u^U(x,y)G_U^X(x,y), \quad x,y\in U.
\end{equation}

\medskip

Let $\chi(r)= (\psi\circ \phi) (r)$.
For $r>0$, define
$$
\chi^r(\lambda):=\frac{\chi(\lambda r^{-2})}{\chi(r^{-2})}\, ,\quad \lambda >0\, .
$$
Then $\chi^{r}$ is a complete Bernstein function with $\chi^{r}(1)=1$ and $\chi^{1}(s)=\chi(s)$ .
 It follows from  \eqref{e:weak-scaling-psi-phi}
that for each $M>0$, $(\chi^{r})_{r \le M}$ satisfies the following weak scaling condition at infinity, uniformly in $r\in (0,M]$: there exists 
$c=c(M)>1$ such that  
\begin{equation}\label{e:weak-scaling-psi-phir}
c^{-1}\left(\frac{S}{s}\right)^{\gamma_1 \delta_1}\le \frac{\chi^r(S)}{\chi^r(s)} \le
c\left(\frac{S}{s}\right)^{\gamma_2 \delta_2}\, , \qquad 1<s\le S<\infty\, .
\end{equation}

Let $K^r=( K^r_t)_{t\ge 0}$ be a subordinator with Laplace exponent $\chi^r$ independent of 
the Brownian motion  $W$. Let $X^r=(X^r_t)_{t\ge 0}$ be defined by $X^r_t:=W_{ K^r_t}$
Then $X^r$ is an isotropic L\'evy process with characteristic exponent $\chi^r(|\xi|^2)=\chi(|\xi|^2 r^{-2}) / \chi(r^{-2})$, $\xi\in \R^d$, which shows that $X^r$ is identical in law to the process 
$\{r^{-1}X_{t/ \chi(r^{-2})}\}_{t\ge 0}$.

Let $V\subset \R^d$ be a bounded $C^{1,1}$ open set.  
For $r\in (0,1]$, let $V^r:=\{rx:\, x\in V\}$. 
Denote by $G_V^{X^r}$ (respectively $G_{V^r}^X$) the Green function of $V$ with respect to  $X^r$ (respectively the Green function of $V^r$ with respect to  $X$). Then by scaling,
\begin{equation}\label{e:scaling-for-G}
G_{V^r}^X(x,y)=r^{-d}\chi(r^{-2})^{-1}G_V^{X^r}(x/r,y/r)\, , \quad x,y\in V^r\, .
\end{equation}
For any open set $U\subset \R^d$, we let
\begin{align}
\label{e:gU}
 g^{r}_U(x,y):=
\left(1\wedge \frac{\chi^r(|x-y|^{-2})}{\sqrt{\chi^r(\delta_U(x)^{-2})\chi^r(\delta_U(y)^{-2})}}\right)\, \frac{1}{|x-y|^{d}\chi^r(|x-y|^{-2})}\, , \quad x,y\in U\,,
\end{align}
and $g_U(x,y):=g^1_U(x,y)$.
\begin{prop}\label{p:scale-inv-gfe}
Let $V\subset \R^d$ be a bounded $C^{1,1}$ open set with characteristics $(R,\Lambda)$  and 
$\mathrm{diam}(V)\le 1$. There exists a constant 
$C=C(R,\Lambda,\phi,\psi, d)\ge 1$ 
such that for every $r\in (0,1]$,
$$
C^{-1}g_{V^r}(x,y)\le G_{V^r}^X(x,y)\le C g_{V^r}(x,y)\, , \qquad x,y\in V^r\, .
$$
The dependence of $c$ on $\phi$ and $\psi$  
is only through the constants in assumptions 
 \eqref{e:weak-scaling-phi} and \eqref{e:weak-scaling-psi}.  
\end{prop}
\pf 
The proof is similar to that of \cite[Proposition 7.5]{KSV16}.
 In fact, by \cite[Theorem 1.2]{KM14} and \eqref{e:weak-scaling-psi-phir},
there exists a constant $c>0$ such that 
for all $r \in (0,1]$, 
\begin{equation}\label{e:scale-inv-gfe1}
c^{-1}g^r_V(x,y)\le G_V^{X^r}(x,y) \le c g^r_V(x,y)\, ,\quad x,y\in V\, .
\end{equation}
Since $g^r_V(x/r,y/r)=r^d \phi(r^{-2})g_{V^r}(x,y)$, 
the claim of the proposition follows by combining this,  \eqref{e:scaling-for-G} and \eqref{e:scale-inv-gfe1}.
 \qed

Note that the fact that we scale around the origin is irrelevant.
For any $z\in \R^d$ we could use the scaling 
$V^r:=\{r(x-z)+z;\, x\in V\}$ and obtain the same result.

Using Proposition \ref{p:scale-inv-gfe} and \eqref{e:Fbound}, the proof of the next result 
is very similar to (and even simpler than) that of \cite[Lemmas 7.6--7.8]{KSV16}. So we omit the proof.
\begin{lemma}\label{l:bounded-gauge}
Let $R>0$ and $\Lambda >0$. There exists $C=C(R,\Lambda, \phi,\psi, d)\in (0,1)$ 
such that for every $r\in (0,1/b]$ and every $C^{1,1}$ open set $U\subset D$ with characteristics $(rR,\Lambda/r)$ and $\mathrm{diam}(U)\le r$ satisfying $\mathrm{dist}(U,\partial D)\ge (b+2)r$, we have
$$
C\le u^U(x,y) \le 1\, , \quad x,y\in U,\,  x\neq y.
$$
\end{lemma}

Combining this lemma with \eqref{e:gf-XU-YDU-comparable} we arrive at
\begin{prop}\label{p:green-comparable}
Let $R>0$ and $\Lambda >0$. There exists $C=C(R,\Lambda, \phi,\psi, d)\in (0,1)$ 
such that for every $r\in (0,1/b]$ and every $C^{1,1}$ open set $U\subset D$ with characteristics $(rR,\Lambda/r)$ and $\mathrm{diam}(U)\le r$ satisfying $\mathrm{dist}(U,\partial D)\ge (b+2)r$, we have 
$$
CG_U^X(x,y)\le G^{Y^D}_U(x,y)\le G_U^X(x,y)\, ,\quad x,y\in U\, .
$$
\end{prop}

\section{Behavior of harmonic functions in the interior of $D$} 
\label{s:hi}

We first recall the 
definitions of $\kappa$-fat  open sets in $\R^d$ and harmonic functions.
\begin{defn}\label{def:UB}
Let $0<\kappa< 1$.
  We say that an open set $D\subset \R^d$ is $\kappa$-fat if
  there is $R_1>0$ such that for all $x\in \overline{D}$ and all $r\in
  (0,R_1]$, there is a ball
  $B(A_r(x), \kappa r)
  \subset D\cap B(x,r)$. The pair $(R_1, \kappa)$ is called the characteristics of
the $\kappa$-fat open set $D$.
\end{defn}
\begin{defn}\label{def:har}
Suppose that $D\subset \R^d$ is an open set.  
(1) 
A non-negative Borel function $u$ on $D$ is said to be
harmonic in an open set $U\subset D$ with respect to $Y^D$ if
for every open set $B$ whose closure is a compact subset of $U$,
\begin{equation}\label{e:har}
u(x)= \E_x \left[ u(Y^D_{\tau_{B}})\right], \qquad
\hbox{for every } x\in B.
\end{equation}

\noindent
(2) 
A non-negative Borel function $u$ on $D$ is said to be regular
harmonic in an open set $U\subset D$ with respect to $Y^D$ if
$$
u(x)= \E_x \left[ u(Y^D_{\tau_{U}})\right],
\qquad \hbox{for every } x\in U.
$$
\end{defn}
Clearly, a regular harmonic function in $U$ is harmonic in $U$.

 The first goal of this section is to prove 
a scale invariant Harnack inequality for
non-negative functions harmonic with respect to $Y^D$ 
when $D$ is a bounded $\kappa$-fat open set in $\R^d$.
Then we study the decay rate of non-negative functions in $D$ which are (regular) harmonic near 
a portion of the boundary, strictly contained in $D$, of an open set $E\subset D$
and vanish locally on $E^c$.

Let $j^Z$ be the L\'evy density of $Z$.
Recall that $\Phi(r)=\frac{1}{\phi(r^{-2})}$.
 Analogously to \eqref{e:JD-estimate}, it follows from \cite[Lemma 3.2]{KSV14} that for any $M>0$ there exists $c(M)\ge 1$ such that 
\begin{equation}\label{e:J-estimate-Z}
c(M)^{-1} \frac1{|x-y|^{d}\Phi(|x-y|)}\le 
j^Z(|x-y|)\le c(M) \frac1{|x-y|^{d}\Phi(|x-y|)}\, ,\qquad |x-y|\le M\, .
\end{equation}

The following result is a consequence of
\cite[Proposition 3.6]{CKS} and it is valid for any bounded open set $D$.
\begin{lemma}\label{l:ilb4BMheatkernel}
Let $D\subset\R^d$ be a bounded open set 
and $a, T$ be positive constants.
There exists $C=C(a, T, \phi, {\rm diam}(D))>0$
such that
\begin{equation}\label{e:ilb4BMheatkernel}
p^D(t, x, y)\ge C\left(\Phi^{-1}(t)^{-d}\wedge \frac{t}{|x-y|^d \Phi(|x-y|)}\right), \qquad t\le 
\Phi(a(\delta_D(x)\wedge\delta_D(y)))
\wedge T.
\end{equation}
\end{lemma}

Note that 
\begin{equation}\label{e:Phi-less-x-y}
\Phi^{-1}(t)^{-d} \le  \frac{t}{|x-y|^d \Phi(|x-y|)} \quad \text{if and only if }\quad t\ge \Phi(|x-y|)\, .
\end{equation}
\begin{lemma}\label{l:ilb4jk}
Let  $D\subset\R^d$ be a bounded open set and $\varepsilon_0\in (0, 1]$ be
a constant.  There exists a constant 
$C=C(\varepsilon_0, \phi, \mathrm{diam}(D))\in (0, 1)$ such that for
every $x_0\in D$ 
and $r\le 1/2$ satisfying $B(x_0, (1+\varepsilon_0)r)\subset D$, we have
\begin{equation}\label{e:ilb4jk}
CJ^X(x, y)\le J^{Y^D}(x, y)\le J^X(x, y), \qquad x, y\in B(x_0, r).
\end{equation}
\end{lemma}

\pf  The second inequality in \eqref{e:ilb4jk} is \eqref{e:JYDJX}. 
 So we only need to prove the first inequality. 
It follows from Lemma \ref{l:ilb4BMheatkernel} that there exists 
$c_1= c_1(\varepsilon_0, \phi, \mathrm{diam}(D))  >0$ 
 such that for any $t\le \Phi((2/\varepsilon_0) (\delta_D(x)\wedge\delta_D(y) ) \wedge \Phi($diam$(D))$,
$$
p^D(t, x, y)\ge c_1\left(\Phi^{-1}(t)^{-d}\wedge \frac{t}{|x-y|^d \Phi(|x-y|)}\right).
$$
Thus using the fact that $s \mapsto \nu(s)$ is decreasing, \eqref{e:Phi-less-x-y} and \eqref{e:v-nu-asymp}, we have that for $x, y\in B(x_0, r)$,
\begin{align*}
&J^{Y^D}(x, y)\ge\int^{\Phi(|x-y|)}_0p^D(t, x, y)\nu(t)dt
\ge \nu(\Phi(|x-y|))\int^{\Phi(|x-y|)}_0p^D(t, x, y)dt\\
&\ge  c_2  \frac{\psi(\Phi(|x-y|)^{-1})}{\Phi(|x-y|)} \int^{\Phi(|x-y|)}_0\frac{t}{|x-y|^d \Phi(|x-y|)}dt
\ge  c_3  \frac{\psi(\Phi(|x-y|)^{-1})}{|x-y|^{d}}\\
&\ge c_4J^X(x, y),
\end{align*}
where in the last inequality we used \eqref{e:JD-estimate}.
\qed

In the remainder of this section we assume that $D$ is a bounded 
$\kappa$-fat open set in $\R^d$.

 Combining \cite[Corollary 1.4]{CKS} with \eqref{e:J-estimate-Z} we get 
that for every $T>0$ there exists $C(T)\ge 1$ such that
\begin{align}\label{e:p-r-comparable-kappa}
&C(T)^{-1}\P_x(\tau^Z_D>t)\P_y(\tau^Z_D>t)\left(\Phi^{-1}(t)^{-d}
\wedge \frac{t}{|x-y|^d \Phi(|x-y|)}\right) \nn\\
&\le p^D(t, x, y) \le  C(T)
\P_x(\tau^Z_D>t)\P_y(\tau^Z_D>t)\left(\Phi^{-1}(t)^{-d}\wedge \frac{t}{|x-y|^d \Phi(|x-y|)}\right)
\end{align}
for all $(t,x,y)\in (0,T]\times D \times D$.
\begin{prop}\label{p:estimate-of-J-away}
Suppose that $D\subset\R^d$ is a bounded $\kappa$-fat open set.
For every $\eps_0 \in (0,1]$, 
there exists a constant $C\ge 1$ 
 such that for all $x_0\in D$ and all $r \le 1$ satisfying $B(x_0, (1+\eps_0)r)\subset D$, it holds that
\begin{equation}\label{e:estimate-of-J-away}
J^{Y^D}(z,x_1)\le C J^{Y^D}(z,x_2)\, ,\quad x_1,x_2\in B(x_0,r), \ \ z\in D\setminus B(x_0, (1+\eps_0)r)\, .
\end{equation}
\end{prop}

\pf 
We closely follow the proof of \cite[Proposition 3.5]{KSV16}. 
Suppose that $r, \eps_0 \le 1$, $B(x_0, (1+\eps_0)r)\subset D$ and $x_1, x_2\in B(x_0, r)$.  
Then it follows from \cite[Lemma 3.2]{CKS} that for $t<\Phi(\eps_0r)$, we have
$$
\P_{x_2}(t<\tau^Z_D)\ge \P_{x_2}(t<\tau^Z_{B(x_2, \Phi^{-1}(t))})\ge c_1
$$
for some constant $c_1=c_1(\eps_0)>0$ independent of $x_2$ and $t<\Phi(\eps_0r)$.
By combining this with \eqref{e:p-r-comparable-kappa} we have that 
there exists $c_2>1$ such that for $z\in D$ and 
 $t < \Phi(\eps_0 r)$,  
\begin{align*}
&p^D(t, z, x_1)\le c_2 \P_z(t<\tau^Z_D)\left(\Phi^{-1}(t)^{-d}\wedge \frac{t}{|x_1-z|^d \Phi(|x_1-z|)}\right),\\
&p^D(t, z, x_2)\ge c_2^{-1}\P_z(t<\tau^Z_D)\left(\Phi^{-1}(t)^{-d}\wedge \frac{t}{|x_2-z|^d \Phi(|x_2-z|)}\right).
\end{align*}
Now suppose that $z\in D\setminus B(x_0,(1+\eps_0)r)$ so that 
$$
\frac{ \eps_0 }{1+\eps_0}|z-x_0|\le |z-x_i|\le \left(1+\frac1{1+\eps_0}\right)|z-x_0|,
\qquad i=1, 2. 
$$
Then
\begin{align}
&\int^{\Phi(\eps_0r)}_0p^D(t, z, x_1)\nu(t)dt\nonumber\\
&\le c_3\int^{\Phi(\eps_0r)}_0\P_z(t<\tau^Z_D)\left(\Phi^{-1}(t)^{-d}\wedge \frac{t}{|z-x_2|^d \Phi(|z-x_2|)}\right)\nu(t)dt \nonumber \\
& \le c_3 c_2 \int_0^{\Phi(\eps_0 r)} p^D(t,z,x_2) \nu(t)\, dt  .\label{e:rs1}
\end{align}

Using the parabolic Harnack principle (see, for instance, \cite[Theorem 1.4]{CKK}),
we get that there exists $c_4>1$ such that
\begin{align*}
&\int^\infty_{\Phi(\eps_0r)}p^D(t, z, x_1)\nu(t)dt=\sum^\infty_{n=1}\int^{(n+1)\Phi(\eps_0r)}_{n\Phi(\eps_0r)}p^D(t, z, x_1)\nu(t)dt\\
&\le c_4\sum^\infty_{n=1}\int^{(n+1)\Phi(\eps_0r)}_{n\Phi(\eps_0r)}p^D(t+\frac{\Phi(\eps_0r)}2, z, x_2)\nu(t)dt\\
&=c_4\sum^\infty_{n=1}\int^{(n+\frac32)\Phi(\eps_0r)}_{(n+\frac12)\Phi(\eps_0r)}p^D(t, z, x_2)\nu(t-\frac{\Phi(\eps_0r)}2)dt\, .
\end{align*}
If 
$t \in (\Phi(\eps_0 r), 1)$, 
 then by \eqref{e:doubling-nu} we have
$\nu(t-\Phi(\eps_0r)/2)\le \nu(t/2)\le c_5 \nu(t)$  with $c_5\ge 1$. 
If $t\ge 1$, then $\nu(t-\Phi(\eps_0r)/2)\le \nu(t-1/2)$. 
By \eqref{e:CBmu0}, there exists $c_6 \ge 1$ such that $\nu(s)\le c_6 \nu(s+1/2)$ for all $s>1/2$. Hence, $\nu(t-1/2)\le c_6\nu(t)$. With $c_7=c_5\vee c_6$,  we conclude that $\nu(t-\Phi(\eps_0r)/2)\le c_7\nu(t)$ for all $t\ge 3\Phi(\eps_0r)/2$. Hence,
\begin{align}
&\int^\infty_{\Phi(\eps_0r)}p^D(t, z, x_1)\nu(t)dt\le c_4 c_7\sum^\infty_{n=1}\int^{(n+\frac32)\Phi(\eps_0r)}_{(n+\frac12)\Phi(\eps_0r)}p^D(t, z, x_2)\nu(t)dt\nonumber \\
&\le c_8\int^\infty_{(3\Phi(\eps_0r))/2}p^D(t, z, x_2)\nu(t)dt\, .
\label{e:rs2}
\end{align}

Combining \eqref{e:JY}, \eqref{e:rs1} and \eqref{e:rs2}, we get that there exists $c_9>1$ such that
$$
J^{Y^D} (z, x_1)=\int^\infty_0p^D(t, z, x_1)\nu(t)dt\le c_9\int^\infty_0p^D(t, z, x_2)\nu(t)dt=
c_9J^{Y^D} (z, x_2)\, ,
$$
which finishes the proof. \qed

Recall that $b \ge 2$ is the constant in  Lemma \ref{l:estimate-of-F}.
\begin{lemma}\label{l:tauB}
Suppose that $D\subset\R^d$ is a bounded $\kappa$-fat open set.
There exists $C>0$ such that for every $x_0\in D$ and every $r\in (0,1)$ 
with $B(x_0,r)\subset D$,
$$
\E_x \tau_{B(x_0, r)}  \ge 
\E_x \tau_{B(x_0,4a r)} 
\ge C \frac1{\psi(\Phi(r)^{-1})}\, ,\quad x\in B(x_0,ar)\, ,
$$
where $a:=2^{-5}b^{-1} \in (0, 2^{-6}]$. 
\end{lemma}
\pf 
The lemma follows from \eqref{e:weak-scaling-psi-phi}, Propositions \ref{p:scale-inv-gfe} and \ref{p:green-comparable}. 
In fact, by Propositions \ref{p:scale-inv-gfe} and \ref{p:green-comparable}, we have that 
for  $x,y\in B(x_0,2ar)$,
$$
G^{Y^D}_{B(x_0,4ar)}(x,y) \ge c_1 G^{X}_{B(x_0,4ar)}(x,y) \ge c_2 \frac1{|x-y|^{d}\psi(\Phi(|x-y|)^{-1})}.
$$
Thus for $x\in B(x_0,ar)$,
\begin{align*}
&\E_x \tau_{B(x_0,4ar)} =
\int_{B(x_0,4ar)}
 G^{Y^D}_{B(x_0,4ar)}(x,y)dy \ge
c_2 \int_{B(x_0,2ar)}  \frac1{|x-y|^{d}\psi(\Phi(|x-y|)^{-1})}dy\\
&\ge c_2 \int_{B(x,ar)}  \frac1{|x-y|^{d}\psi(\Phi(|x-y|)^{-1})}dy \ge c_3
\int_{0}^{ar}\psi(\Phi(s)^{-1}) s^{-1}ds.
\end{align*}
Since \eqref{e:weak-scaling-psi-phi} implies that 
\begin{align*}
&\int_{0}^{ar}\frac{ds}{\psi(\Phi(s)^{-1}) s}= 
\frac1{\psi(\Phi(ar)^{-1})} \int_{0}^{ar}\frac{\psi(\Phi(ar)^{-1})}{\psi(\Phi(s)^{-1}) s}ds\\
&\ge c_4 \frac1{\psi(\Phi(r)^{-1})} \int_{0}^{ar}\left(\frac{s}{ a r}\right)^{2\gamma_2\delta_2} s^{-1}ds 
\ge c_5 \frac1{\psi(\Phi(r)^{-1})} ,
\end{align*}
we have proved the lemma. 
\qed

By using  
\eqref{e:greenup}, \eqref{e:GX-upper-estimate} and \eqref{e:weak-scaling-psi-phi}, one gets that  there exists $c>0$ such that for every $x_0\in D$ 
and every $r\in (0,1)$ with $B(x_0,r)\subset D$,  we have
\begin{equation}\label{e:exit-time-est}
\E_x \tau_{B(x_0,r)} \le c \frac1{\psi(\Phi(r)^{-1})}\, ,\quad x\in B(x_0,r)\, .
\end{equation}

For any open set $U\subset D$, let
$$
K^{D,U}(x,z):=\int_U G^{Y^D}_U(x,y)J^{Y^D}(y,z)\, dy\, ,\quad x\in U, z\in \overline{U}^c\cap D
$$
be the Poisson kernel of $Y^D$ on $U$.
\begin{thm}
[Harnack inequality]\label{uhp}
Suppose that $D\subset\R^d$ is a bounded $\kappa$-fat open set.
There exists a constant $C>0$ such that for any
$r\in (0,1]$ and $B(x_0, r) \subset D$ and 
any Borel function 
$f$ which is non-negative in
$D$ and harmonic in $B(x_0, r)$ with respect to $Y^D$, we have
$$
f(x)\le C f(y), \qquad \text{ for all } x, y\in B(x_0, r/2).
$$
\end{thm}
\pf
Let $b_1=a/2$ and $b_2=a/4$, where $a$ is the constant from Lemma \ref{l:tauB}.
We will first show that there exists a constant 
$c_1>0$, independent of $x_0$ and $r$,
such that for all $y_0 \in B(x_0,r/2)$, $x_1,x_2\in B(y_0, b_2r/4)$ and $z\in \overline{B(y_0,b_1r)}^c\cap D$,
\begin{equation}\label{e:pk-est}
K^{D,B(y_0,b_1r)}(x_1,z)\le c_1 K^{D,B(y_0,b_1r)}(x_2,z)\, .
\end{equation}
Let $A(x, r_1, r_2):=\{y\in \R^d: r_1\le |y-x|<r_2\}$. Note that  
\begin{align}
&K^{D,B(y_0,b_1r)}(x_1,z)\nn\\
&=\int_{B(y_0,b_2 r)}G^{Y^D}_{B(y_0,b_1r)}(x_1,y)J^{Y^D}(y,z)\, dy+\int_{A(y_0,b_2 r,b_1r)}G^{Y^D}_{B(y_0,b_1r)}(x_1,y)J^{Y^D}(y,z)\, dy \nn\\
&=:I_1+I_2\, .\label{e:aaa0}
\end{align}
In order to estimate $I_2$ we use 
Propositions \ref{p:scale-inv-gfe} and \ref{p:green-comparable}  together with 
$|x_1-y|\asymp^{c_2} |x_2-y|$ for all $y \in A(y_0,b_2r, b_1r)$ 
 and $\delta_{B(y_0, b_1 r)}(x_1) \asymp^{c_3} \delta_{B(y_0, b_1 r)}(x_2)$ 
to get 
\begin{align}\label{e:aaa1}
I_2 
&\le  c_4 \int_{A(y_0,b_2r, b_1r)}G^{Y^D}_{B(y_0,b_1r)}(x_2,y)J^{Y^D}(y,z)\, dy \le c_5 K^{D,B(y_0,b_1r)}(x_2,z)\, .
\end{align}
By Proposition \ref{p:estimate-of-J-away} (with $\eps_0=1$),   $J^{Y^D}(y_0,z)\asymp^{c_6} J^{Y^D}(y,z)$
for $z\in B(y_0,b_1r)^c= B(y_0,2b_2 r)^c$ and $y\in B(y_0,b_2 r)$.
Hence, by \eqref{e:exit-time-est} and  Lemma \ref{l:tauB}, 
\begin{eqnarray}
I_1&\le & c_7 J^{Y^D}(y_0,z)\int_{B(y_0,b_2 r)}G^{Y^D}_{B(y_0,b_1r)}(x_1,y)\, dy \le c_7 J^{Y^D}(y_0,z) \E_{x_1}\tau_{B(y_0,b_1r)}\nn\\
&\le &c_8 J^{Y^D}(y_0,z)\psi(\Phi(r)^{-1})^{-1}\le  c_9 J^{Y^D}(y_0,z) \E_{x_2}\tau_{B(y_0, b_2 r)}\nn\\
&=&c_9 J^{Y^D}(y_0,z)\int_{B(y_0,b_2 r)}G^{Y^D}_{B(y_0,b_2 r)}(x_2,y)\, dy\nn\\
&\le & c_9 J^{Y^D}(y_0,z)\int_{B(y_0,b_2 r)}G^{Y^D}_{B(y_0,b_1r)}(x_2,y)\, dy\nn\\
&\le & c_{10} \int_{B(y_0,b_2 r)}G^{Y^D}_{B(y_0,b_1r)}(x_2,y)J^{Y^D}(y,z)\, dy \le c_{10} K^{D,B(y_0,b_1r)}(x_2,z)\, .
\label{e:aaa2}
\end{eqnarray}
Combining \eqref{e:aaa0}--\eqref{e:aaa2}, 
we have proved \eqref{e:pk-est}. Now, let $f$ be a 
non-negative function in $D$ which is harmonic with respect to $Y^D$ in $B(x_0,r)$. 
Then, by the L\'evy system formula and \eqref{e:exit-time-YD},
for all $y_0 \in B(x_0,r/2)$ and $x_1,x_2 \in B(y_0,b_1 r/8)$ we have
\begin{align*}
&f(x_1)=\int_{\overline{B(y_0,b_1r)}^c}K^{D,B(y_0,b_1r)}(x_1,z)f(z)\, dz\\
&\le c_{11} \int_{\overline{B(y_0,b_1r)}^c}
K^{D,B(y_0,b_1 r)}
(x_2,z)f(z)\, dz=c_{11} f(x_2)\, .
\end{align*}
For $x_1, x_2\in B(x_0,r/2)$, the theorem follows by a standard chain argument.
\qed

Now we prove the following version of the Harnack inequality.
\begin{thm}\label{HP2}
Suppose that $D\subset\R^d$ is a bounded $\kappa$-fat open set.
There exists a constant $C=C(\phi, \psi, \mathrm{diam}(D))>1$ such that the
following is true: If $L>0$ and $x_{1}, x_{2}\in D$ and $r\in (0,1)$ are such
that $|x_{1}-x_{2}|< Lr$ and $B(x_{1}, r)\cup B(x_{2},r) \subset D$, then for any Borel function 
$f$ which is non-negative in
$D$ and harmonic  in $B(x_{1}, r)\cup B(x_{2},r)$ with respect to $Y^D$, we
have
$$
C^{-1}  (L \vee 1)^{-d-\delta_2}  f(x_{2})\,\leq \,f(x_{1}) \,\leq\, C (L \vee 1)^{d+\delta_2}  f(x_{2}).
$$
\end{thm}

\pf
Let $r \in (0,1)$,  $x_{1}, x_{2}\in D$ be such that $|x_1-x_2|< Lr$ and
$B(x_{1}, r)\cup B(x_{2},r) \subset D$. Let
$f$ be a non-negative function which is harmonic in
$B(x_{1}, r)\cup B(x_{2},r)$ with respect to $Y^D$.
If $ |x_{1}- x_{2}|  <  \frac14 r$, then since $r<1$, 
the claim is true by Theorem \ref{uhp}.
Thus we only need to consider the case when $\frac14 r \le |x_{1}- x_{2}| \le L r$ with $L > \frac14$.

By Lemma \ref{e:ilb4BMheatkernel}, \eqref{e:weak-scaling-phi} and \eqref{e:v-nu-asymp}, for every
$(x,y) \in B(x_2, \frac{r}{16}) \times B(x_1, \frac{r}{16})$,
\begin{align}
&J^{Y^D}(x, y)\ge\int^{\Phi(\delta_D(x)\wedge\delta_D(y))}_0p^D(t, x, y)\nu(t)dt \nn\\
&\ge c_1\int^{\Phi(\delta_D(x)\wedge\delta_D(y))}_0
 \left(\Phi^{-1}(t)^{-d} \wedge   \frac{t }{|x-y|^d \Phi(|x-y|)} \right)
\nu(t) dt \nn\\\
&\ge c_2 \int^{\Phi(r/8)}_0\frac{t \nu(t)}{|x-y|^d \Phi(|x-y|)}
dt \ge c_3 \nu(\Phi(r/8) ) \int^{\Phi(r/8)}_0\frac{t }{|x-y|^d \Phi(|x-y|)} dt \nn \\
&\ge  c_4  \nu(\Phi(r/8) )\int^{\Phi(r/8)}_0\frac{t}{(2Lr)^d \Phi(2Lr)}dt
\ge  c_5  {\psi(\Phi(r/8)^{-1})} \frac{\Phi(r/8)}{(2Lr)^d \Phi(Lr)} \nn\\
&\ge  c_6 L^{-d-\delta_2} r^{-d}  {\psi(\Phi(r/8)^{-1})}. \label{e:JYDbigger}
\end{align}

Note that, by Proposition \ref{p:estimate-of-J-away},  for every
$y \in B(x_1, \frac{r}{16})$,
 it holds that
$$ 
K^{D, B(x_2, \frac{r}{16})}(x_2, y)=
\int_{B(x_2,\frac{r}{16})} G^{Y^D}_{B(x_2,\frac{r}{16})}(x_2,z)J^{Y^D}(z, y) dz\,
\ge\, c_7\, J^{Y^D}(x_2, y) \E_{x_2} [\tau_{B(x_2, \frac{r}{16})}].
$$
Thus using this, \eqref{e:JYDbigger} and Lemma \ref{l:tauB}, we have that 
for every $y \in B(x_1, \frac{r}{16})$,
\begin{equation}\label{e:asdads1}
K^{D, B(x_2,\frac{r}{16})}(x_2, y) \,\ge\, 
c_7c_6 L^{-d-\delta_2} r^{-d}  \frac{\psi(\Phi(r/8)^{-1})}{\psi(\Phi( {r}/(16))^{-1})} 
\ge c_8  L^{-d-\delta_2} r^{-d}.
\end{equation}

For any $y\in B(x_1, \frac{r}{16})$, $f$ is regular harmonic in $B(y,\frac{15r}{16})\cup B(x_1, \frac{15r}{16})$. 
Since $|y-x_1|< \frac{r}{16}$, by Theorem \ref{uhp},
\begin{equation}\label{e:asdads2} 
f(y)\ge c_9f(x_1), \quad y\in B(x_1, \frac{r}{16}),
\end{equation}
for some constant $c_9>0$. Therefore, by \eqref{e:asdads1},
\begin{eqnarray*}
f(x_2) &=& \E_{x_2}\left[f(Y^D_{\tau_{B(x_2,\frac{r}{16}) }})\right]
\ge \E_{x_2}\left[f(Y^D_{\tau_{B(x_2,\frac{r}{16})}});
Y^D_{\tau_{B(x_2,\frac{r}{16})}} \in B(x_1,\frac{r}{16}) \right]\\
&\ge &c_{10} \,f(x_1)\, \P_{x_2}\left(Y^D_{\tau_{B(x_2,\frac{r}{16})}}
\in B(x_1,\frac{r}{16}) \right)
= c_{10}\, f(x_1) \int_{B(x_1,\frac{r}{16})} K^{D, B(x_2,\frac{r}{16})} (x_2,w)\, dw \\
&\ge  &c_{11}L^{-d-\delta_2} f(x_1)  |B(x_1,\frac{r}{16})| r^{-d}   = c_{12}\,L^{-d-\delta_2} f(x_1).
\end{eqnarray*}
Thus we have proved the theorem.
\qed
\begin{remark}{\rm
Suppose that $d=1$ and that $D$ is a 
bounded open set which is a union of disjoint open
intervals so that the infimum of the lengths of all these intervals is at least $R_0$ and the infimum of the distances between these intervals is at least $R_0$.

Let $\widetilde D=D\times [0,1]$
 and let $\widetilde{Y}^{\widetilde{D}}$ be the subordinate killed Brownian motion in $\widetilde D\subset
\R^{2}$. 
By checking the definitions, we can see that if $h$ is harmonic in $U\subset D$ with
respect to $Y^D$, then $\widetilde h(x, y)=h(x)$, $x\in D, y\in 
[0,1]$, is harmonic in $U\times [0,1]^2$ with respect 
to $\widetilde{Y}^{\widetilde{D}}$. 
Note that $U\times [0,1]$ is a $\kappa$-fat open set in $\R^2$.
Thus, in fact, 
using Theorems \ref{uhp} and \ref{HP2} for $d= 2$, 
Theorems \ref{uhp} and \ref{HP2} also hold 
for $d= 1$. 
}
\end{remark}
\begin{thm}\label{t:bhp-2}
Suppose $d\ge 2$. 
Let $D\subset \R^d$ be a bounded $\kappa$-fat open set. 
There exists a constant 
$b=b(\phi, \psi, d)>2$ such that, for 
every open set $E\subset D$ and every $Q\in \partial E\cap D$  such that $E$ is $C^{1,1}$ near $Q$ with characteristics 
$(R, \Lambda)$, the following holds: 
There exists a constant 
$C=C(\delta_D(Q)\wedge R,\Lambda, \psi, \phi,d)>0$ 
such that for every 
$r\le (\delta_D(Q)\wedge 1)/(b+2)$ and every non-negative function $f$ on $D$ which is regular harmonic in $E\cap B(Q,r)$ with respect to $Y^D$ and vanishes on $E^c\cap B(Q,r)$, we have 
$$
\sqrt{(\psi\circ \phi)(\delta_E(x)^{-2})}f(x)\le C \sqrt{(\psi\circ \phi)(\delta_E(y)^{-2})} f(y)\, , \qquad x,y\in E\cap B(Q,2^{-6}\kappa_0^4 r)\, ,
$$
where $\kappa_0=(1+(1+\Lambda)^2)^{-1/2}$.
\end{thm}

\pf
Using  \eqref{e:weak-scaling-psi-phi}, \cite[Lemma 2.2]{So}, Proposition \ref{p:green-comparable}, Proposition \ref{p:estimate-of-J-away}
and 
 the factorization from either \cite[Lemma 5.5]{KSV12} or \cite[Lemma 5.4]{KM14}, the proof is the same as that of 
  \cite[Theorem 1.3]{KSV16} 
 with $\psi\circ \phi$ instead of $\phi$. So we omit the 
  details. 
\qed

%%%%%%%%%%%%%%%%%%%%%%%%%%%%%%%%%%%%%%%%%%%%%%%%%%%%%%%%%%%%%%%%%
%%%%%%%%%%%%%%                     Carleson estimate in $C^{1,1}$ domain                 %%%%%%%%%%%%%%%%%%%%
%%%%%%%%%%%%%%%%%%%%%%%%%%%%%%%%%%%%%%%%%%%%%%%%%%%%%%%%%%%%%%%%%

\section{Carleson estimate}\label{s:ce}

 We will 
establish the Carleson estimate for $Y^D$. 
Unlike \cite{KSV16}, neither the explicit boundary behavior of the jumping kernel nor that of
the Green function is used in the proof of the Carleson estimate. 
The Carleson estimate of $Y^D$
is established for a large class of non-smooth open sets.

Using  \eqref{e:YD-killed-inside},  the proof of the next lemma is the same as that of \cite[Lemma 5.1]{KSV16}. 
\begin{lemma}\label{l:regularity}
Suppose that $D\subset \R^d$ is an open set. Let $x_0\in \R^d$, and $r_1<r_2$ be two positive numbers such that $ D \cap B(x_0, r_1) \neq \emptyset$.
Suppose $f$ is a non-negative function in $D$ that is
harmonic in $D\cap B(x_0, r_2)$ with respect to $Y^D$ and vanishes
continuously on $\partial D\cap B(x_0, r_2)$. Then $f$ is regular harmonic in
$D\cap B(x_0, r_1)$ with respect to $Y^D$, i.e.,
\begin{equation}\label{e:regularity}
f(x)=\E_x\left[ f\big(Y^D(\tau_{D\cap B(x_0, r_1)})\big)\right]
 \qquad \hbox{
for all }x\in D\cap B(x_0, r_1)\, .
\end{equation}
\end{lemma}

For $x\in D$, let $z_x$ be a point on $\partial D$ 
such that $|z_x-x|=\delta_D(x)$.
We say $D\subset\R^d$ satisfies the local exterior volume condition with 
characteristics $(R_0, C_0)$
if  for every $z \in \partial D$ and  $x \in B(z, R_0) \cap D$,
 $|D^c\cap B(z_x,\delta_D(x))| \ge C_0 \delta_D(x)^{d}$. It is easy to see that, if $D^c$ is $\kappa$-fat, then $D$ satisfies the local exterior volume condition.
 
We recall that $\zeta$ is  the lifetime of $Y^D$.
Let
\begin{align}
\label{e:defg}
g(r):=\frac{1}{r^d \psi(\Phi(r)^{-1})}\, ,\qquad r>0\, .
\end{align}
\begin{lemma}\label{lower bound} 
Suppose that $D\subset\R^d$ satisfies the local exterior volume condition with 
characteristics $(R_0, C_0)$.
Then there exists a constant
$\delta_*=\delta_*(R_0, C_0)>0$ such that for all
$x\in D$ with $\delta_D(x) < R_0/2$,
$$
\P_x\left(  \tau(x)=\zeta     \right)\ge \delta_*\, ,
$$
where $\tau(x):=\tau_{
B(x,\delta_D(x)/2)}=\inf\{t>0:\,
Y^D_t\notin
B(x,\delta_D(x)/2)\}$.
\end{lemma}

\pf
By \cite[Theorem 4.5.4(1)]{FOT},
$$
\P_x\left(  \tau(x)=\zeta     \right)=\P_x(Y^D_{\zeta-}\in 
B(x,\delta_D(x)/2))
= \int_{
B(x,\delta_D(x)/2)} G^{Y^D}(x,y) \kappa^{Y^D} (y)dy,
$$
where  $\kappa^{Y^D}$ 
is the density of the killing measure of $Y^D$ given in \eqref{e:kms}. 
Since $D$ satisfies the local exterior volume condition,
we have (see the proof of \cite[Proposition 5.12]{KSV11}) 
\begin{equation}\label{e:estimate-kappa-D}
\kappa^{Y^D} (y) \ge \kappa^{X^D}(y)\ge c_1 \psi(\Phi (\delta_D(y))^{-1}), \qquad y \in 
B(x,\delta_D(x)/2).
\end{equation}
Here \eqref{e:killing-functions-relation} is used in the first inequality. 
Thus, using \eqref{e:weak-scaling-psi-phi}, \eqref{e:estimate-kappa-D}, 
\eqref{e:UB-smaller-UDB}, and 
Propositions \ref{p:scale-inv-gfe} and \ref{p:green-comparable},
\begin{align*}
&\P_x\left(  \tau(x)=\zeta     \right) = \int_{B(x,\delta_D(x)/2)} G^{Y^D}(x,y) \kappa^{Y^D} (y)dy\\
&\ge \int_{B(x,\delta_D(x)/(4b))} G_{B(x,\delta_D(x)/(8b))}^{Y^D}(x,y) \kappa^{Y^D} (y)dy\\
&\ge c_2 \int_{ B(x,\delta_D(x)/(4b))} g(|x-y|)\psi(\Phi (\delta_D(y))^{-1})dy\\
& \ge c_3 \psi(\Phi (\delta_D(x))^{-1})   
\int_{ B(x, \delta_D(x)/(4b))} \frac{1}{|x-y|^d \psi(\Phi(|x-y|)^{-1})} dy \\
&\ge c_4 \psi(\Phi (\delta_D(x))^{-1})
\frac{1}{\psi(\Phi (\delta_D(x)/(4b))^{-1})} \ge c_5,
\end{align*}
where $b>2$ is the constant in Lemma \ref{l:estimate-of-F}.
 \qed

In the remainder of this section we will assume 
$D\subset\R^d$ is a bounded $\kappa$-fat open set with characteristics $(R_1, \kappa)$.
Combining \eqref{e:JY} and 
\cite[Theorem 1.3(iii) and Corollary 1.4]{CKS},
 we immediately get the following
\begin{prop}\label{p:est-JYkappa} 
For any $T>0$, there exists 
$C=C(R_1, \kappa, T)\ge 1$ such that for all $x, y\in D$,
$$
C^{-1}\wt{J}^D(x, y)\le J^{Y^D}(x, y)\le C\wt{J}^D(x, y), 
$$
where
\begin{align*}
\wt{J}^D(x, y)
=&\int_0^T \P_x(\tau^Z_D>t)\P_y(\tau^Z_D>t)\left(\Phi^{-1}(t)^{-d}\wedge \frac{t}{|x-y|^d \Phi(|x-y|)}\right) 
\nu(t)dt\\
&\qquad+\P_x(\tau^Z_D>1)\P_y(\tau^Z_D>1)\, .
\end{align*}
\end{prop}

Before we prove the Carleson estimate for $Y^D$,
we first show the following form of parabolic Carleson type estimate for $Z$.
\begin{prop}\label{p:pCeZ}
For any $T>0$ and $c_0 \in (0,1)$, there exists $C=C(R_1, \kappa, c_0, T)\ge 1$ such that for all $t \in (0,T]$, $r \le R_1/2$, $Q\in \partial D$ and 
$x, x_0\in D \cap B(Q, r)$ with $\delta_D(x_0) \ge c_0r$,
\begin{align}\label{e:pCeZ}
\P_x(\tau^Z_D>t) \le  C  \P_{x_0}(\tau^Z_D>t).
\end{align}
\end{prop}
\pf For simplicity, without loss of generality we assume 
$T=R_1=1$.
In this proof, we always assume that $t, r \in (0,1]$ and 
$x, x_0 \in D \cap B(Q, r)$ 
with $\delta_D(x_0) \ge c_0r$. 

\noindent
{\it Case 1. $r \ge 2^{-4} \kappa \Phi^{-1}(t)/3$:}
In this case by  \cite[Lemma 3.2]{CKS} 
and \eqref{e:Phi(lambda-t)},
\begin{align*}
&\P_{x_0}(\tau^Z_D>t) \ge    \P_{x_0} (\tau^Z_{B(x_0, c_0 r)}  > t) \ge
\P_{x_0} (\tau^Z_{B(x_0, c_0 r)}  > \Phi(3\cdot 2^4  \kappa^{-1}r)) \\
& \ge  \P_{x_0} (\tau^Z_{B(x_0, c_0 r)}  > c_2 \Phi(c_0r)) \ge c_3>0.
\end{align*}
Thus $\P_x(\tau^Z_D>t) \le 1 \le  c_3^{-1}  \P_{x_0}(\tau^Z_D>t).$

\noindent
{\it Case 2. $r \le 2^{-4} \kappa \Phi^{-1}(t)/3$:}
In this case, we will use  \cite[Lemma 4.1]{CKS}. 
Let $A:=A_{\Phi^{-1}(t)/2}(Q)$ 
so that $B(A, \kappa \Phi^{-1}(t)/2) \subset B(Q, \Phi^{-1}(t)/2) \cap D$.
Since $B(Q, \Phi^{-1}(t)/2) \subset B(y, \Phi^{-1}(t))$ for all $y \in B(Q, r)$, 
we have 
\begin{align}\label{e:pCeZ1}
B(A, 2^{-2}\kappa \Phi^{-1}(t))\subset B(A, 2^{-1}\kappa \Phi^{-1}(t))\subset B(y, \Phi^{-1}(t)) \cap D,
\quad \forall y \in B(Q, r)  \cap D.
\end{align}
Define
\begin{align}
U(x):= D\cap B \big(x,|x-A|+ 2^{-2}\kappa \Phi^{-1}(t)/3 \big), \label{e:UV}\\
\wh U(x_0):= D\cap B \big(x_0,|x_0-A|+2^{-1} \kappa \Phi^{-1}(t)/3 \big).
\label{e:UV2}
\end{align}
Note that $U(x) \subset  \wh U(x_0)$. In fact, 
by assumption, we have $|x-x_0| \le 2^{-3} \kappa \Phi^{-1}(t)/3$,  so 
for $y \in U(x)$,
\begin{align*}
|x_0-y| &\le |x-x_0|+ |x-A|+ 2^{-2}\kappa \Phi^{-1}(t)/3 \\
&\le |x_0-A|+ 2|x-x_0|+ 2^{-2}\kappa \Phi^{-1}(t)/3\\
& \le |x_0-A|+2^{-1} \kappa \Phi^{-1}(t)/3 .
\end{align*}
Since $D$ is $2^{-2}\kappa$-fat and $A=A_{2^{-2}\kappa \Phi^{-1}(t)}(x)$ by \eqref{e:pCeZ1},
 using  \cite[Lemma 4.1]{CKS}, we have 
 $$ \P_x(\tau^Z_D>t) \le c_4  \P_x(Z_{\tau_{U(x)}}\in D).$$
 Since $y \to \P_y(Z_{\tau_{U(x)}}\in D)$  is  regular harmonic in 
  $D\cap B(Q, 2r)$
  with respect to $Z$ and vanishes in $D^c \cap B(Q, 2r)$, 
by \cite[Lemma 5.5]{KSV12},
$$
\P_x(Z_{\tau_{U(x)}}\in D) \le c_5 \P_{x_0}(Z_{\tau_{U(x)}}\in D).
$$ 
Since $x_0\in B\left(x, 2^{-3}\kappa \Phi^{-1}(t)/3\right)\subset
B\left(x, 2^{-1}(|x-A|+2^{-2}\kappa \Phi^{-1}(t)/3)\right)$,
by \cite[Lemma 2.4]{CK},
$$
\P_{x_0}(Z_{\tau_{U(x)}}\in D) \le \P_{x_0}\left(Z_{\tau_{U(x)}}\in
B\left(x, |x-A|+ 2^{-2}\kappa \Phi^{-1}(t)/3 \right)^c\right)\le c_6  t^{-1} \E_{x_0}[\tau^Z_{U(x)}].
$$
Since  $U(x) \subset  \wh U(x_0)$, combining the above inequalities we get
$$
 \P_x(\tau^Z_D>t) \le c_4c_5c_6  t^{-1}
    \E_{x_0}[\tau^Z_{\wh U(x_0)}].
$$
Finally, since $D$ is $2^{-1}\kappa$-fat and $A=A_{2^{-1}\kappa \Phi^{-1}(t)}(x_0)$ by \eqref{e:pCeZ1},
using  \cite[Lemma 4.1]{CKS}, we have 
$$
 \P_x(\tau^Z_D>t) \le c_4c_5c_6  t^{-1}
  \E_{x_0}[\tau^Z_{\wh U(x_0)}] 
 \le c_7  \P_{x_0}(\tau^Z_D>t).
$$
\qed
\begin{thm}[Carleson estimate]\label{t:carleson}
Suppose that $D\subset\R^d$ is a bounded 
$\kappa$-fat open set 
with  characteristics $(R_1, \kappa)$
satisfying the local exterior volume condition with  characteristics $(R_0, C_0)$.
There exists a constant
$C=C(R_1,\kappa, R_0, C_0)>0$
 such that for every $Q\in \partial D$, 
$r\in (0, (R_0 \wedge R_1)/2)$, 
and every non-negative function
$f$ in $D$ that is harmonic in $D \cap B(Q, r)$ with respect to
$Y^D$ and vanishes continuously on $ \partial D \cap B(Q, r)$, we have
\begin{equation}\label{e:carleson}
f(x)\le C f(x_0) \qquad \hbox{for }  x\in D\cap B(Q,r/2),
\end{equation}
where $x_0\in D \cap B(Q,r)$ with $\delta_D(x_0)\ge \kappa r/2$.
\end{thm}

\pf
In this proof, the constants $\delta_*, \nu, \gamma, \beta_1, \eta$ and $c_i$'s are always
independent of $r$. Without loss of generality,
we assume that $R_0=R_1 \le 1$ and diam($D$) $\le 1$. 
By
Theorem \ref{HP2},
it suffices to prove (\ref{e:carleson}) for 
$x\in D\cap B(Q, \kappa r/(24))$.

Choose $0<\gamma < \frac{2\gamma_1 \delta_1}{d+2}\wedge  \frac12$. 
For any $x\in D\cap B(Q,\kappa r/(12))$, define
$$
D_0(x)=D\cap B(x,2\delta_D(x))\, ,\qquad B_1(x)=B(x,r^{1-\gamma}
\delta_D(x)^{\gamma})\,
$$
and
$$
B_2=B(x_0,\kappa\delta_D(x_0)/3)\, ,\qquad 
B_3=B(x_0, 2\kappa\delta_D(x_0)/3).
$$
Since $x\in B(Q,\kappa r/(12))$, we have $\delta_D(x)<r/(12)$. By 
the choice of $\gamma<1/2$, 
we have that $D_0(x)\subset B_1(x)$.
By Lemma \ref{lower bound}, there exists $\delta_*=\delta_*
(R_0, C_0)>0$ such that
\begin{equation}\label{e:c:1}
\P_x(\tau_{D_0(x)}=\zeta)\ge
\P_x(\tau_{B(x,\delta_D(x)/2)}=\zeta)  \ge \delta_*\, ,
\quad x\in D\cap B(Q,\kappa r/(12))\, .
\end{equation}
Further, by \eqref{e:weak-scaling-psi-phi},
\eqref{e:GX-upper-estimate} and \eqref{e:greenup},
\begin{align}\label{e:exit-from-D_0}
\E_x \tau_{D_0(x)}^{Y^D}&\le \int_{D_0(x)}G^{Y^D}(x,y)\, dy 
\le c_1\int_{B(x,2\delta_D(x))}g(|x-y|)\, dy \nonumber \\
&\le c_2  \int_{B(x,2\delta_D(x))} \frac{1}{|x-y|^d \psi(\Phi(|x-y|)^{-1})}dy
\le \frac{c_3}{ \psi(\Phi(\delta_D(x))^{-1})}\, ,
\end{align}
where
$g$ was defined in \eqref{e:defg}.
By Theorem \ref{HP2},  we have 
\begin{equation}\label{e:c:2}
f(x)<c_4(\delta_D(x)/r)^{-\beta_1} f(x_0)\, ,
\quad x\in D\cap B(Q,\kappa r/(12))\,,
\end{equation}
where $\beta_1:=d+\delta_2>0$.
Since $f$ is regular harmonic in
$D_0(x)$ with respect to $Y^D$ by Lemma \ref{l:regularity}, 
for every $x\in D\cap B(Q,\kappa r/(12)))$,
\begin{align}
f(x)=&\E_x\big[f\big(Y^D(\tau_{D_0(x)})\big); Y^D(\tau_{D_0(x)})\in
B_1(x)\big]\nonumber\\
&+ \E_x\big[f\big(Y^D(\tau_{D_0(x)})\big);
Y^D(\tau_{D_0(x)})\notin B_1(x)\big]. \label{e:c:3}
\end{align}
We first show that there exists 
$\eta \in (0, 2^{-4})$
 such that
for all $x\in D \cap B(Q, \kappa r/(12)) $ 
with $\delta_D(x) < \eta r$,
\begin{equation}\label{e:c:4}
\E_x\big[f\big(Y^D(\tau_{D_0(x)})\big); Y^D(\tau_{D_0(x)})\notin
B_1(x)\big]\le f(x_0). 
\end{equation}
 Since $\gamma<\frac12$, we have $2^{-4}<4^{-(1-\gamma)^{-1}}$.
Thus
for $\delta_D(x)< 2^{-4} r$,
$$
2\delta_D(x) \le r^{1-\gamma} \delta_D(x)^{\gamma} - 2\delta_D(x).
$$
Hence, if $x\in D \cap B(Q, \kappa r/(12))$
with $\delta_D(x) < 2^{-4}r$,  then
$|x-y|\le 2|z-y|$ for $z\in D_0(x)$, $y\notin B_1(x)$. 

If $z\in B_2$ and $y\in D\setminus B_3$, then it follows from 
Proposition \ref{p:estimate-of-J-away} 
(with $\eps_0=1$ and $r=\kappa \delta_D(x_0)/3$) 
that $J^{Y^D}(z,y)\ge c_4 J^{Y^D}(x_0,y)$. 
By using this estimate in the third line below, Lemma \ref{l:tauB} in the fourth and $\delta_D(x_0) \ge \kappa r/2$ in the fifth, we get that
\begin{align}\label{e:c:6-new}
&f(x_0)\ge \E_{x_0}\left[f(Y^D(\tau_{B_2})); Y^D(\tau_{B_2})\notin B_3\right]\nonumber \\
&=    \E_{x_0} \int_0^{\tau_{B_2}}\left( \int_{D\setminus B_3} J^{Y^D}(Y^D_t, y)f(y)\, dy \right) dt\nonumber \\
&\ge c_5 
\E_{x_0}[\tau_{B_2}]
\int_{D\setminus B_3}  J^{Y^D}(x_0, y)f(y)\, dy \nonumber \\
&\ge \frac{c_6}{ \psi(\Phi(\delta_D(x_0))^{-1})}\int_{D\setminus B_3}  J^{Y^D}(x_0, y)f(y)\, dy \nonumber \\
&\ge\frac{ c_7}{\psi(\Phi(r)^{-1})}\int_{D\setminus B_3}  J^{Y^D}(x_0, y)f(y)\, dy\, .
\end{align}
Next,
\begin{align}\label{e:c:5-new}
&\E_x\big[f\big(Y^D(\tau_{D_0(x)})\big); Y^D(\tau_{D_0(x)}) \notin B_1(x)\big]\nonumber\\
&=\E_x \int_0^{\tau_{D_0(x)}}    \int_{D\setminus B_1(x)} J^{Y^D}(Y^D_t,y)f(y)\, dy\, dt \nonumber \\
&=\E_x \int_0^{\tau_{D_0(x)}}    \int_{(D\setminus B_1(x))\cap B_3^c} J^{Y^D}(Y^D_t,y)f(y)\, dy\, dt \nonumber \\
&\quad +\E_x \int_0^{\tau_{D_0(x)}}    \int_{(D\setminus B_1(x))\cap B_3} J^{Y^D}(Y^D_t,y)f(y)\, dy\, dt =: I_1+I_2\, .
\end{align}

In order to estimate $I_2$ we first use 
 Theorem \ref{HP2},
\eqref{e:JYDJX} and  \eqref{e:JD-estimate} to get 
\begin{align}\label{e:I2-first-part}
&I_2\le c_8  f(x_0)\E_x \int_0^{\tau_{D_0(x)}^{Y^D}} \int_{(D\setminus B_1(x))\cap B_3} \frac{\psi(\Phi(|Y_t^D-y|)^{-1})}{|Y_t^D-y|^d}\, dy \, dt. \end{align}
Since 
$|z-y|\ge \frac12 |x-y|$ for $(z, y)\in D_0(x) \times (D\setminus B_1(x))$ and 
$|y-x|\ge |x_0-Q|-|x-Q|-|y-x_0|>\kappa\delta_D(x_0)/6 \ge \kappa^2 r/12$ for $y\in B_3$, 
we have 
$$
{\psi(\Phi(|z-y|)^{-1})}|z-y|^{-d} \le c_7 {\psi(\Phi(r)^{-1})}{r^{-d}}\quad  \text{for }(z, y) \in D_0(x) \times 
(B_3\cap (D\setminus B_1(x))).
$$
 Therefore, by using \eqref{e:exit-from-D_0} in the second inequality and the fact that 
  $\delta_D(x_0) < r$ in the third, we have
\begin{align}\label{e:I2-second-part}
&I_2\le c_{9} f(x_0)\E_x [\tau_{D_0(x)}^{Y^D}] \int_{(D\setminus B_1(x))\cap B_3} \frac{\psi(\Phi(r)^{-1})}{r^d}\, dy\nonumber \\
&\le c_{10} f(x_0)\frac1{\psi(\Phi(\delta_D(x))^{-1})} |B_3| \frac{\psi(\Phi(r)^{-1})}{r^d}\nonumber \\
&\le c_{11} f(x_0)\frac{\psi(\Phi(r)^{-1})}{\psi(\Phi(\delta_D(x))^{-1})} \, .
\end{align}

In order to estimate $I_1$ 
we use Propositions \ref{p:est-JYkappa} and \ref{p:pCeZ}.
If $w\in D_0(x)$, then $\delta_D(w)\le r$. 
If we further assume $y\in D\setminus B_1(x)$, 
then $|w-y|\ge \frac12 |x-y|$, implying that $|w-y|^d\Phi(|w-y|)\ge 2^{-d-2}|x-y|^d\Phi(|x-y|)$.
Therefore, using Propositions \ref{p:est-JYkappa}--\ref{p:pCeZ}
and the fact that $a\wedge (cb)\le c(a\wedge b)$ for $a,b>0$ 
and $c\ge 1$,
\begin{align}\label{e:I1-first-part}
&I_1\le 
c_{12}\E_x\int_0^{\tau_{D_0(x)}^{Y^D}}\int_{(D\setminus B_1(x))\cap B_3^c}\left( \int_0^1  \P_{x_0}(\tau^Z_D>s)\P_y(\tau^Z_D>s)
\right. \nonumber \\
& \quad \left. \times \left(\Phi^{-1}(s)^{-d}\wedge \frac{2^{d+2}s}{|x-y|^d \Phi(|x-y|)}\right)\nu(s)ds+\P_{x_0}(\tau^Z_D>1)\P_y(\tau^Z_D>1)\right) f(y)dy\,  dt \nonumber \\
&\le 2^{d+2}c_{12}\E_x [\tau_{D_0(x)}^{Y^D}]\int_{(D\setminus B_1(x))\cap B_3^c}\left( \int_0^1 
\P_{x_0}(\tau^Z_D>s)\P_y(\tau^Z_D>s)
\right. \nonumber \\
&\quad \left. \times \left(\Phi^{-1}(s)^{-d}\wedge \frac{s}{|x-y|^d \Phi(|x-y|)}\right)\nu(s)ds +\P_{x_0}(\tau^Z_D>1)\P_y(\tau^Z_D>1)\right) f(y)dy\, .
\end{align}
Recall that $x\in B(Q, \kappa r/(24))$. For $y\in D\setminus B_1(x)$ we have $|y-x|\ge r^{1-\gamma}\delta_D(x)^{\gamma}$, and therefore 
$$
|y-x_0|\le |y-x|+r\le
|y-x|+r^{\gamma}\delta_D(x)^{-\gamma}|y-x|\le 2r^{\gamma}\delta_D(x)^{-\gamma}|y-x|\, .
$$
Thus by \eqref{e:Phi(lambda-t)},
we have $\Phi(|y-x_0|)\le  (2r^{\gamma}\delta_D(x)^{-\gamma})^2 \Phi(|x-y|)$.
Hence, by using 
\eqref{e:exit-from-D_0} in the first
inequality below, the fact that $a\wedge (cb)\le c(a\wedge b)$ for $a,b>0$ 
and $c\ge 1$ in the second, 
Proposition \ref{p:est-JYkappa} in the third,
and  \eqref{e:c:6-new} in the last inequality, we have
\begin{align}\label{e:I1-second-part}
& I_1 \le 
\frac{c_{13}}{\psi(\Phi(\delta_D(x))^{-1})}
\int_{(D\setminus B_1(x))\cap B_3^c}\left( \int_0^1  \P_{x_0}(\tau^Z_D>s)\P_y(\tau^Z_D>s)
\right. \nonumber \\
&\quad \left. \times \left(\Phi^{-1}(s)^{-d}\wedge \frac{(2r/\delta_D(x))^{\gamma(d+2)}s}{|x_0-y|^d \Phi(|x_0-y|)}\right)\nu(s)ds+\P_{x_0}(\tau^Z_D>1)\P_y(\tau^Z_D>1)\right) f(y)dy \nonumber \\
& \le 
\frac{c_{13} 2^{\gamma(d+2)}}{ \psi(\Phi(\delta_D(x))^{-1})} 
 \left(\frac{\delta_D(x)}{r}\right)^{-\gamma(d+2)}   \int_{D\setminus B_3}
 \left(\int_0^1\P_{x_0}(\tau^Z_D>s)\P_y(\tau^Z_D>s)\right. \nonumber \\
& \quad \left. \times \left(\Phi^{-1}(s)^{-d}\wedge \frac{s}{|x_0-y|^d \Phi(|x_0-y|)}\right)\nu(s)ds+\P_{x_0}(\tau^Z_D>1)\P_y(\tau^Z_D>1) \right)f(y)dy \nonumber \\
&\le 
\frac{c_{14}} { \psi(\Phi(\delta_D(x))^{-1})} \left(\frac{\delta_D(x)}{r}\right)^{-\gamma(d+2)}  \int_{D\setminus B_3} J^{Y^D}(x_0,y) f(y)\, dy\nonumber \\
&\le c_{15} \frac{\psi(\Phi(r)^{-1})}{\psi(\Phi(\delta_D(x))^{-1})} \left(\frac{\delta_D(x)}{r}\right)^{-\gamma(d+2)} f(x_0)\, .
\end{align}

Combining \eqref{e:c:5-new}, \eqref{e:I2-second-part} and \eqref{e:I1-second-part}, we obtain
\begin{align}\label{e:c:9-new}
&\E_x[f(Y^D(\tau_{D_0(x)}));\, Y^D(\tau_{D_0(x)})\notin B_1(x)]\nonumber\\
&\le c_{16} f(x_0)\frac{\psi(\Phi(r)^{-1})}{\psi(\Phi(\delta_D(x))^{-1})}\left( \left(\frac{\delta_D(x)}{r}\right)^{-\gamma(d+2)} +1\right) \nn \\
&\le c_{17} f(x_0)\left(\frac{\delta_D(x)}{r}\right)^{2\delta_1 \gamma_1}
\left( \left(\frac{\delta_D(x)}{r}\right)^{-\gamma(d+2)}+1\right)\, .
\end{align}
Since $2\gamma_1 \delta_1 -(d+2)\gamma>0$, 
we can choose $\eta\in (0, 2^{-4})$ so that
$$
c_{17}\,
\left(\eta^{2\gamma_1 \delta_1 -(d+2)\gamma}+\eta^{2\gamma_1 \delta_1}\right)
\,\le\, 1\, .
$$
Then  for $x\in  D \cap B(Q, \kappa r/(12))$ with $\delta_D(x) < \eta r$, we
have by \eqref{e:c:9-new},
\begin{eqnarray*}
\E_x\left[f(Y^D(\tau_{D_0(x)}));\, Y^D(\tau_{D_0(x)})\notin B_1(x)\right] &\le & c_{17}\, f(x_0)
\left(\eta^{2\gamma_1 \delta_1 -(d+2)\gamma}+\eta^{2\gamma_1 \delta_1} \right)
\le f(x_0)\, .
\end{eqnarray*}
This completes the proof of \eqref{e:c:4}.

With \eqref{e:c:4}, one can prove the Carleson estimate \eqref{e:carleson} for 
$x\in D\cap B(Q, \kappa r/(24))$ 
by a method of contradiction. Since this part of the proof is the same as the corresponding part in the proof of 
\cite[Theorem 5.4]{KSV16}, we omit 
the  details.   \qed

%%%%%%%%%%%%%%%%%%%%%%%%%%%%%%%%%%%%%%%%%%%%%%%%%%%%%%%%%%%%%%%%%%%
%%%%%%%%%%%%%%%              Jumping kernel and Green function estimates                  %%%%%%%%%%%%%%%%%%%
%%%%%%%%%%%%%%%%%%%%%%%%%%%%%%%%%%%%%%%%%%%%%%%%%%%%%%%%%%%%%%%%%%%

\section{Green function and exit time estimates}
\label{s:jkgfe}

In this section, we assume that $D$ is a bounded $C^{1,1}$ open set in  $\R^d$,  $d\ge 2$,
with $C^{1, 1}$ characteristics $(R, \Lambda)$.
The first goal of this section is to derive sharp two-sided estimates for  $G^{Y^D}$. See \cite[Theorem 3.1]{KSV16m}
for the corresponding result for killed subordinate Brownian motion.

Recall that $\Phi(r)=\frac{1}{\phi(r^{-2})}$.
For $t>0$ and $x,y\in D$, let
\begin{equation}\label{e:rtxy}
r(t,x,y):=\left(\frac{\Phi(\delta_D(x))^{1/2}}{t^{1/2}}\wedge 1\right)
\left(\frac{\Phi(\delta_D(y))^{1/2}}{t^{1/2}}\wedge 1\right)\left(\Phi^{-1}(t)^{-d}\wedge \frac{t}{|x-y|^d \Phi(|x-y|)}\right)\, .
\end{equation}
 Combining \cite[Corollary 1.6]{CKS} with \eqref{e:J-estimate-Z} we get 
that for every $T>0$ there exist $C_1=C_1(T, R, \Lambda, a_2, a_2, \delta_1, \delta_2)\ge 1$  and $C_2=C_2(T, R, \Lambda, a_2, a_2, \delta_1, \delta_2, {\rm diam}(D))\ge 1$
such that 
\begin{equation}\label{e:p-r-comparable}
C_1^{-1}r(t,x,y)\le p^D(t,x,y)\le C_1 r(t,x,y)
\end{equation}
for all $(t,x,y)\in (0,T]\times D \times D$,  and
\begin{equation}\label{e:p-large-time}
C_2^{-1}e^{-\lambda_1 t} \Phi(\delta_D(x))^{1/2}\Phi(\delta_D(y))^{1/2} \le p^D(t,x,y) \le C_2 e^{-\lambda_1 t} \Phi(\delta_D(x))^{1/2}\Phi(\delta_D(y))^{1/2}
\end{equation}
for all $(t,x,y)\in (T,\infty)\times D \times D$. Here $-\lambda_1<0$ is the largest eigenvalue of the infinitesimal generator of $Z^D$.

It follows easily from \eqref{e:weak-scaling-psi} and \eqref{e:v-nu-asymp} that
\begin{equation}\label{e:weak-scaling-v}
b_2^{-1} \left(\frac{s}{t}\right)^{1-\gamma_2}\le \frac{v(t)}{v(s)} \le b_1^{-1} \left(\frac{s}{t}\right)^{1-\gamma_1}\, ,\quad 0<t\le s \le 1\, .
\end{equation}

 It is shown in \cite[Lemma 7.1]{CKS} that 
\begin{equation}\label{e:a(x,y)-comparable}
\left(\frac{\Phi(\delta_D(x))^{1/2}\Phi(\delta_D(y))^{1/2}}{\Phi(|x-y|)}\wedge 1\right) \asymp^c \left(\frac{\Phi(\delta_D(x))^{1/2}}{\Phi(|x-y|)^{1/2}}\wedge 1\right) \left(\frac{\Phi(\delta_D(y))^{1/2}}{\Phi(|x-y|)^{1/2}}\wedge 1\right)\, .
\end{equation}
\begin{lemma}\label{l:third-integral}
Let $f:[0,\infty)\to [0,\infty)$ be a decreasing function.
For any $T>0$,
there exists $C=C(f,\phi, \mathrm{diam}(D),T)>0$ such that for all $x,y \in D$,
$$
\int_T^{\infty}p^D(t,x,y)f(t)\, dt \le C
\left(\frac{\Phi(\delta_D(x))^{1/2}}{\Phi(|x-y|)^{1/2}}\wedge 1\right) \left(\frac{\Phi(\delta_D(y))^{1/2}}{\Phi(|x-y|)^{1/2}}\wedge 1\right)\frac{\Phi(|x-y|)f(\Phi(|x-y|))}{|x-y|^d}\, .
$$
\end{lemma}
\pf 
For $x,y\in D$, let 
$a(x,y):=\Phi(\delta_D(x))^{1/2}\Phi(\delta_D(y))^{1/2}$
and 
$$b(x,y):=\left(\frac{a(x,y)}{\Phi(|x-y|)}\wedge 1\right) \frac{\Phi(|x-y|)f(\Phi(|x-y|))}{|x-y|^d}.$$
If $a(x,y)\le \Phi(|x-y|)$, then
\begin{eqnarray}\label{e:bbb1}
b(x,y)=a(x,y)\frac{f(\Phi(|x-y|))}{|x-y|^d}\ge \frac{f(\Phi(\mathrm{diam}(D)))}{\mathrm{diam}(D)^d} a(x,y)\, .
\end{eqnarray}
If $a(x,y)> \Phi(|x-y|)$, then by \eqref{e:td-Phi-decreasing},
\begin{align}\label{e:bbb2}
&b(x,y)=f(\Phi(|x-y|))\frac{\Phi(|x-y|)}{|x-y|^d} \nonumber\\
&\ge f(\Phi(\mathrm{diam}(D)))\frac{\Phi(\mathrm{diam}(D))}{\mathrm{diam}(D)^d}\ge \frac{f(\Phi(\mathrm{diam}(D)))}{\mathrm{diam}(D)^d} a(x,y)\, .
\end{align}
Since by \eqref{e:p-large-time},
$$
\int_T^{\infty}p^D(t,x,y)f(t)\, dt \le  c_3 a(x,y)f(T)\int_T^{\infty} e^{-\lambda_1 t}\, dt \le c_4 a(x,y)\,,
$$
the claim now follows from \eqref{e:bbb1}, \eqref{e:bbb2}
 and \eqref{e:a(x,y)-comparable}.
\qed
\begin{lemma} \label{l:second-integral}
Let $T>\Phi(\mathrm{diam}(D))$.
There exists a constant
$C=C(T,  a_2, \delta_2 )\ge 1$ 
such that for all $x,y\in D$,
\begin{eqnarray*}
\int_{\Phi(|x-y|)}^T r(t,x,y) dt
\le C\left(\frac{\Phi(\delta_D(x))^{1/2}}{\Phi(|x-y|)^{1/2}}\wedge 1\right)
\left(\frac{\Phi(\delta_D(y))^{1/2}}{\Phi(|x-y|)^{1/2}}\wedge 1\right)
\frac{\Phi(|x-y|)}{|x-y|^d}\, .
\end{eqnarray*}
\end{lemma}
\pf First note that by \eqref{e:Phi-1-scaling}, for $\Phi(|x-y|) \le t \le T$,
$$
\Phi^{-1}(t)^{-d}\le C_T^d |x-y|^{-d}\Phi(|x-y|)^{d/2\delta_2} t^{-d/2\delta_2}\, .
$$
Therefore, using \eqref{e:Phi-less-x-y} we have
\begin{eqnarray*}
\lefteqn{\int_{\Phi(|x-y|)}^T r(t,x,y) dt}\\
&\le &\int_{\Phi(|x-y|)}^T \left(\frac{\Phi(\delta_D(x))^{1/2}}{\Phi(|x-y|)^{1/2}}\wedge 1\right) \left(\frac{\Phi(\delta_D(y))^{1/2}}{\Phi(|x-y|)^{1/2}}\wedge 1\right) \Phi^{-1}(t)^{-d}\, dt \\
&\le & c_1\left(\frac{\Phi(\delta_D(x))^{1/2}}{\Phi(|x-y|)^{1/2}}\wedge 1\right) \left(\frac{\Phi(\delta_D(y))^{1/2}}{\Phi(|x-y|)^{1/2}}\wedge 1\right) 
\frac{\Phi(|x-y|)^{d/2\delta_2}}{|x-y|^d}
\int_{\Phi(|x-y|)}^{\infty} t^{-d/2\delta_2} \, dt \\
&=&c_2 \left(\frac{\Phi(\delta_D(x))^{1/2}}{\Phi(|x-y|)^{1/2}}\wedge 1\right) \left(\frac{\Phi(\delta_D(y))^{1/2}}{\Phi(|x-y|)^{1/2}}\wedge 1\right) 
\frac{\Phi(|x-y|)}{|x-y|^d}\, .
\end{eqnarray*}
\qed

Lemmas \ref{l:third-integral} and \ref{l:second-integral} will also be used in Section \ref{s:j}. 
\begin{lemma}\label{l:green-main-integral}
There exists a constant $C\ge 1$ such that for all $x,y\in D$,
\begin{align*}
&\int_0^{\Phi(|x-y|)} r(t,x,y) v(t)\, dt \\
\asymp^{C} 
&\left(\frac{\Phi(\delta_D(x))^{1/2}}{\Phi(|x-y|)^{1/2}}\wedge 1\right) \left(\frac{\Phi(\delta_D(y))^{1/2}}{\Phi(|x-y|)^{1/2}}\wedge 1\right)\frac{\Phi(|x-y|)v(\Phi(|x-y|))}{|x-y|^d}\, .
\end{align*}
\end{lemma}
\pf 
 By using \eqref{e:rtxy}, \eqref{e:Phi-less-x-y} and \eqref{e:weak-scaling-v} in the second line, the change of variable $s=\Phi(|x-y|)t^{-1}$ in the third line, the fact that  $(sa)\wedge  1\le s(a\wedge 1)$ for every $a>0$ and $s\ge 1$  in the fourth line, we get
\begin{align}
&\int_0^{\Phi(|x-y|)} r(t,x,y) v(t)\, dt\nn\\
&\le c_1 \frac{v(\Phi(|x-y|))\Phi(|x-y|)^{-\gamma_1+1}}{|x-y|^d\Phi(|x-y|)}
 \int_0^{\Phi(|x-y|)} \left(\frac{\Phi(\delta_D(x))^{1/2}}{t^{1/2}}\wedge 1\right) \left(\frac{\Phi(\delta_D(y))^{1/2}}{t^{1/2}}\wedge 1\right) t^{\gamma_1}\, dt \nn\\ 
&=c_1 \frac{v(\Phi(|x-y|))\Phi(|x-y|)}{|x-y|^d} \int_1^{\infty} 
 \left(\frac{s^{1/2}\Phi(\delta_D(x))^{1/2}}{\Phi(|x-y|)^{1/2}}\wedge 1\right) \left(\frac{s^{1/2}\Phi(\delta_D(y))^{1/2}}{\Phi(|x-y|)^{1/2}}\wedge 1\right) s^{-\gamma_1-2}\, ds\nn\\
  &\le c_1
   \left(\frac{\Phi(\delta_D(x))^{1/2}}{\Phi(|x-y|)^{1/2}}\wedge 1\right)
 \left(\frac{\Phi(\delta_D(y))^{1/2}}{\Phi(|x-y|)^{1/2}}\wedge 1\right)\frac{v(\Phi(|x-y|))\Phi(|x-y|)}{|x-y|^d} \int_1^{\infty} s^{-\gamma_1-1}\, ds\, . \label{e:ddd1}
\end{align}
Since $\int_1^{\infty} s^{-\gamma_1-1}\, ds <\infty$, this proves the upper bound.

On the other hand, by using  \eqref{e:rtxy}, \eqref{e:Phi-less-x-y}  and the fact that $v$ is decreasing 
in the second line, we get
\begin{eqnarray}
\lefteqn{\int_0^{\Phi(|x-y|)} r(t,x,y) v(t)\, dt} \nn\\
 &\ge & \frac{v(\Phi(|x-y|))}{|x-y|^d\Phi(|x-y|)}  \int_0^{\Phi(|x-y|)} \left(\frac{\Phi(\delta_D(x))^{1/2}}{t^{1/2}}\wedge 1\right) \left(\frac{\Phi(\delta_D(y))^{1/2}}{t^{1/2}}\wedge 1\right) t\, dt \nn\\ 
 &\ge &\left(\frac{\Phi(\delta_D(x))^{1/2}}{\Phi(|x-y|)^{1/2}}\wedge 1\right) \left(\frac{\Phi(\delta_D(y))^{1/2}}{\Phi(|x-y|)^{1/2}}\wedge 1\right) \frac{v(\Phi(|x-y|))}{|x-y|^d\Phi(|x-y|)} \int_0^{\Phi(|x-y|)} t\, dt\nn\\
  &=&\frac12 \left(\frac{\Phi(\delta_D(x))^{1/2}}{\Phi(|x-y|)^{1/2}}\wedge 1\right) \left(\frac{\Phi(\delta_D(y))^{1/2}}{\Phi(|x-y|)^{1/2}}\wedge 1\right)\frac{v(\Phi(|x-y|))\Phi(|x-y|)}{|x-y|^d}\, .\label{e:ddd2}
\end{eqnarray}
\qed

Recall the function $g$ defined in \eqref{e:defg}.
Note that $g$ satisfies the doubling property near 0 by \eqref{e:weak-scaling-psi-phi}
 and, by \eqref{e:psi-Phi-inv-increasing}, $g$ is a decreasing function.
\begin{thm}\label{t:green-function-estimate}
There exists a constant $C\ge 1$ 
such that for all $x,y\in D$,
$$
G^{Y^D}(x,y)\asymp^C \left(\frac{\Phi(\delta_D(x))^{1/2}}{\Phi(|x-y|)^{1/2}}\wedge 1\right) \left(\frac{\Phi(\delta_D(y))^{1/2}}{\Phi(|x-y|)^{1/2}}\wedge 1\right) g(|x-y|)\, .
$$
\end{thm}
\pf Let $T=2\Phi(\mathrm{diam}(D))$. By \eqref{e:green-YD}, \eqref{e:rtxy}, \eqref{e:p-r-comparable} and \eqref{e:p-large-time},
\begin{align*}
&G^{Y^D}(x,y)=\int_0^{\infty}p^D(t,x,y)v(t)\, dt \\
&\le  c_1 \left(\int_0^{\Phi(|x-y|)} r(t,x,y)v(t)\, dt +
\int_{\Phi(|x-y|)}^Tr(t,x,y)v(t)\, dt +\int_T^{\infty} p^D(t,x,y)v(t)\, dt\right) \\
&=:c_1(I_1+I_2+I_3)\, .
\end{align*}
By Lemmas \ref{l:third-integral}--\ref{l:green-main-integral} and the fact that $v$ is decreasing, for each $j=1,2,3$,
\begin{eqnarray*}
I_j &\le & c_2 \left(\frac{\Phi(\delta_D(x))^{1/2}}{\Phi(|x-y|)^{1/2}}\wedge 1\right) \left(\frac{\Phi(\delta_D(y))^{1/2}}{\Phi(|x-y|)^{1/2}}\wedge 1\right)\frac{\Phi(|x-y|)v(\Phi(|x-y|))}{|x-y|^d}\\
&\le & c_3 \left(\frac{\Phi(\delta_D(x))^{1/2}}{\Phi(|x-y|)^{1/2}}\wedge 1\right) \left(\frac{\Phi(\delta_D(y))^{1/2}}{\Phi(|x-y|)^{1/2}}\wedge 1\right)g(|x-y|)\, .
\end{eqnarray*}
Here the last line follows from \eqref{e:v-nu-asymp} and the definition of $g$.
 For the lower bound we 
use $G^{Y^D}(x,y)\ge c_4 I_1$ and  Lemma \ref{l:green-main-integral}.
\qed 

 Let $Q\in \partial D$ and choose $\varphi_Q$ as in Definition \ref{d:C11}. 
Define $\rho_Q (x) := x_d -  \varphi_Q(\wt x),$ where $(\wt x,
x_d)$ are the coordinates of $x$ in $CS_Q$. Note that for every $Q
\in \partial D$ and $ x \in B(Q, R)\cap D$, we have
\begin{equation}\label{e:d_com}
(1+\Lambda^2)^{-1/2} \,\rho_Q (x) \,\le\, \delta_D(x)  \,\le\,
\rho_Q(x).
\end{equation}
We define for $r_1, r_2>0$,
\begin{align}
\label{e:DQ}
D_Q( r_1, r_2) :=\left\{ y\in D: r_1 >\rho_Q(y) >0,\, |\wt y | < r_2
\right\}.
\end{align}

 Let $\kappa_0=(1+(1+ \Lambda)^2)^{-1/2}$.  
It is well known (see, for instance, \cite[Lemma
2.2]{So}) that there exists $L=L(R, \Lambda, d)>0$ such that for
every $z \in\partial D$ and $r \le \kappa_0 R$, one can find a $C^{1,1}$
domain $V_Q(r)$ with characteristics $(rR/L, \Lambda L/r)$
such that $D_Q( r/2, r/2) \subset V_Q(r)  \subset  D_Q( r, r) $. 
In this and the following two sections,
given a $C^{1, 1}$ 
open set $D$, 
$V_Q(r)$ always refers to the $C^{1, 1}$ 
domain above.

It is easy to see that for every $Q \in\partial D$ and $r \le \kappa_0 R$,
\begin{align}
\label{e:Uzr}
V_Q(r) \subset  D_Q( r, r) \subset D\cap B(Q,  r /\kappa_0).
\end{align}
In fact, for all $y \in D_Q(r,r)$, 
 \bee\label{e:lsd}
|y|^2 = |\wt y|^2+ |y_d|^2 <r^2 +(|y_d- \varphi_Q(\wt y)|+ |\varphi_Q(\wt
y)|)^2 < (1+(1+ \Lambda)^2) r^2.
\eee

For any $r\le 1$, let $\phi^r$ be defined by 
$\phi^r(\lambda):=\phi(\lambda r^{-2})/\phi(r^{-2})$. Then $\phi^r$ is also a
complete Bernstein function. Let $S^r$
be a subordinator independent of the Brownian motion $W$ and let $Z^r$ be defined by 
$Z^r_t:=W_{S^r_t}$. 
Then, cf. \cite[p.~247]{KSV14},
$Z^r$ is identical in law to the process $\{r^{-1}Z_{t/\phi(r^{-2})}\}_{t\ge 0}$. 
Let $ p^r (t, x, y)$ be the transition density of $Z^r$. 
For any open set $U$, let $p^{r, U}(t, x, y)$ be the transition density of $Z^{r, U}$, the subprocess of $Z^r$
killed upon exiting $U$.

By the fact that $r^{-1}V_Q(r)$ is a $C^{1, 1}$ open set with $C^{1, 1}$ characteristics $(R/L, \Lambda L)$, 
we get from  \eqref{e:p-r-comparable}--\eqref{e:p-large-time} 
 that there exists $C_3\ge1$   such that 
$$
 C_3^{-1} \widetilde r(t,x,y)\le p^{r, r^{-1}V_Q(r)}(t,x,y)\le  C_3  \widetilde r(t,x,y)
$$
for all $(t,x,y)\in (0, a_2(2/\kappa_0)^{2 \delta_2}]\times r^{-1}V_Q(r) \times r^{-1}V_Q(r)$, where
\begin{align*}
&\widetilde r(t,x,y)\\
&:=\left(\frac{\Phi(\delta_{r^{-1}V_Q(r)}(x))^{1/2}}{t^{1/2}}\wedge 1\right)
\left(\frac{\Phi(\delta_{r^{-1}V_Q(r)}(y))^{1/2}}{t^{1/2}}\wedge 1\right)\left(\Phi^{-1}(t)^{-d}\wedge \frac{t}{|x-y|^d \Phi(|x-y|)}\right)\, .
\end{align*}
By the scaling property mentioned in the paragraph above we get that
$$
p^{V_Q(r)}(t, x, y)=r^{-d}p^{r, r^{-1}V_Q(r)}(\phi(r^{-2})t, r^{-1}x, r^{-1}y)=r^{-d}p^{r, r^{-1}V_Q(r)}(\Phi(r)^{-1}t, r^{-1}x, r^{-1}y).
$$
Thus, since $\Phi(2r/\kappa_0) \le  a_2(2/\kappa_0)^{2 \delta_2} \Phi(r)$   by \eqref{e:weak-scaling-phi}, 
we have
\begin{equation}\label{e:r-hat}
 C_3^{-1} \widehat r(t,x,y)\le  p^{V_Q(r)}  (t,x,y)\le  C_3  \widehat r(t,x,y)
\end{equation}
for all $(t,x,y)\in (0,\Phi(2r/\kappa_0)]\times V_Q(r) \times V_Q(r)$, where
\begin{align*}
&\widehat r(t,x,y)\\
&:=\left(\frac{\Phi(\delta_{V_Q(r)}(x))^{1/2}}{t^{1/2}}\wedge 1\right)
\left(\frac{\Phi(\delta_{V_Q(r)}(y))^{1/2}}{t^{1/2}}\wedge 1\right)\left(\Phi^{-1}(t)^{-d}\wedge \frac{t}{|x-y|^d \Phi(|x-y|)}\right)\, .
\end{align*}

Since diam$(V_Q(r)) \le 2r/\kappa_0$, using the lower bound in  \eqref{e:r-hat} and 
following the argument in \eqref{e:ddd2}, 
one can easily prove the following 
\begin{prop}\label{p:gfcnlb-scaling}
There exists $C=C(R, \Lambda)\ge 1$ such that for all $Q\in \partial D$,
 $r\le \kappa_0R$   
and $x, y\in V_Q(r)$,
$$
G^{Y^{V_Q(r)}}(x, y)\ge C^{-1}
\left(\frac{\Phi(\delta_{V_Q(r)}(x))^{1/2}}{\Phi(|x-y|)^{1/2}}\wedge 1\right) \left(\frac{\Phi(\delta_{V_Q(r)}(y))^{1/2}}{\Phi(|x-y|)^{1/2}}\wedge 1\right) g(|x-y|)\, .
$$
\end{prop}

For simplicity, let $\tau_D:=\zeta$. 
Then $\E_x \tau_D=\int_D G^{Y^D}(x,y)dy$.
The final goal of the section is to give 
sharp two-sided estimates on  $\E_x \tau_D$.
Lemmas \ref{l:gamma-small} and \ref{l:gamma-1/2} below
will be used in 
Section \ref{s:counterexample}.
\begin{lemma}\label{l:gamma-large}
If $\gamma_1>1/2$, then there exists 
$C=C(D, \phi,\psi)>0$ such that
$$
\E_x \tau_D 
\asymp^C
 \Phi(\delta_D(x))^{1/2}\, ,\qquad x\in D\, .
$$
\end{lemma}
\pf 
Let $T=\mathrm{diam}(D)$. 
By using Theorem \ref{t:green-function-estimate} in the first inequality, 
\eqref{e:weak-scaling-psi} in the penultimate inequality, and \eqref{e:weak-scaling-phi}
 and $\gamma_1>1/2$ in the last inequality, we get that
\begin{align}
& \E_x \tau_D= \int_DG^{Y^D}(x,y)\, dy \nn \\
&\le c_{1} \Phi(\delta_D(x))^{1/2} \int_D \frac{1}{\Phi(|x-y|)^{1/2} |x-y|^d \psi(\Phi(|x-y|)^{-1})}\, dy \nn \\
&\le c_{2} \Phi(\delta_D(x))^{1/2}
 \int_{B(x,T)} \frac{1}{\Phi(|x-y|)^{1/2} |x-y|^d \psi(\Phi(|x-y|)^{-1})}\, dy \nn \\
 &\le c_{3} \Phi(\delta_D(x))^{1/2}
 \int_{0}^{T} \frac{1}{\Phi(s)^{1/2} s\psi(\Phi(s)^{-1})}\, ds \nn \\
 &\le c_{4} \frac{\Phi(\delta_D(x))^{1/2}}{\Phi(T)^{1/2}\psi(\Phi(T
 )^{-1})}\int_0^{T} \frac{ds}{s}\left(\frac{\Phi(T)}{\Phi(s)}\right)^{\frac12-\gamma_1} \nn \\
& \le c_{5} \Phi(\delta_D(x))^{1/2}\, .
\end{align}

 For the lower bound, recall that any $C^{1,1,}$ open set satisfies the interior ball condition with some radius $\wh{r}$. 
 Let $a:=(\mathrm{diam}(D)/10)\wedge \wh{r}$ 
 and $D_a:=\{y\in D:\, \delta_D(y)>a\}$.

\noindent Case (i): $\delta_D(x)<a$. Then $B(x,\delta_D(x)/2)\cap D_{2a}=\emptyset$. For $y\in D\setminus B(x, \delta_D(x)/2)$, we have
$$
2|x-y|\ge \delta_D(x), \qquad 3|x-y|\ge \delta_D(x)+|x-y|\ge \delta_D(y),
$$
hence by Theorem \ref{t:green-function-estimate},
$$
G^{Y^D}(x,y)\ge c_6 \frac{\Phi(\delta_D(x))^{1/2}\Phi(\delta_D(y))^{1/2}}{\Phi(|x-y|)}g(|x-y|)\, .
$$
Recall from (2.11) that $\theta(t):=\Phi(t)\psi(\Phi(t)^{-1})$ is increasing, hence $t\mapsto t^d \theta(t)$ is also increasing. Thus we have
\begin{align*}
& \E_x \tau_D \ge c_7\int_{D_{2a}}\Phi(\delta_D(x))^{1/2}\Phi(\delta_D(y))^{1/2} \frac{1}{\Phi(|x-y|)|x-y|^d \psi(\Phi(|x-y|)^{-1})}\\
& \ge c_7\Phi(\delta_D(x))^{1/2}\Phi(2a)^{1/2}\int_{D_{2a}}\frac{dy}{|x-y|^d \theta(|x-y|)}\\
&\ge c_8\Phi(\delta_D(x))^{1/2} \frac{|D_{2a}|}{(\mathrm{diam}(D))^d \theta(\mathrm{diam}(D))}\\
&= c_9 \Phi(\delta_D(x))^{1/2}\, .
\end{align*}

\noindent Case (ii): $\delta_D(x)\ge a$. 
Since $T \ge \delta_D(x)\ge a$, we have
\begin{align*}
& \E_x \tau_D \ge c_{10}\int_{B(x,\delta_D(x)/2)} g(|x-y|) 
\ge c_{11} \int_0^{a/2} \frac{ds}{s\psi(\Phi(s)^{-1})} \ge c_{12}\Phi(\delta_D(x))^{1/2}\, . 
\end{align*}
\qed
\begin{lemma}\label{l:gamma-small}
If $\gamma_2<1/2$, then there exists 
$C=C(D, \phi,\psi)>0$ such that
$$
\E_x \tau_D 
\asymp^C
 \frac{1}{\psi(\Phi(\delta_D(x))^{-1})}\, ,\qquad x\in D\, .
$$
\end{lemma}
\pf Let $T=\mathrm{diam}(D)$. By using 
\eqref{e:weak-scaling-psi}  and \eqref{e:weak-scaling-phi}, we have
\begin{align*}
& \E_x \tau_D\le c_1 \int_D \left(\frac{\Phi(\delta_D(x))^{1/2}}{\Phi(|x-y|)^{1/2}}\wedge 1\right) g(|x-y|)\, dy \\
&= c_1  \int_{D\cap \{|x-y|<\delta_D(x)\}} g(|x-y|)\, dy \\
& \ \ \ +
c_1 \int_{D\cap \{|x-y|\ge \delta_D(x)\}}\frac{\Phi(\delta_D(x))^{1/2}}{\Phi(|x-y|)^{1/2} |x-y|^d \psi(\Phi(|x-y|)^{-1})}\, dy \\
& \le c_2 \frac{1}{\psi(\Phi(\delta_D(x))^{-1})}+ 
\frac{c_3}{\psi(\Phi(\delta_D(x))^{-1})}
\int_{\delta_D(x)}^T \frac{\Phi(\delta_D(x))^{1/2}\psi(\Phi(\delta_D(x))^{-1})}{\Phi(s)^{1/2} s\psi(\Phi(s)^{-1})}\, ds \\
& \le c_2 \frac{1}{\psi(\Phi(\delta_D(x))^{-1})}+ \frac{c_4}{\psi(\Phi(\delta_D(x))^{-1})} \int_{\delta_D(x)}^T \frac1s
\left(\frac{\Phi(\delta_D(x))}{\Phi(s)}\right)^{\frac12-\gamma_2}\, ds \\
&\le c_5\frac{1}{\psi(\Phi(\delta_D(x))^{-1})}\, .
\end{align*}
The lower bound is an immediate consequence of Lemma 4.6:
$$
\E_x \tau_D\ge \E_x \tau_{B(x,\delta_D(x))}\ge c_6 
\frac{1}{\psi(\Phi(\delta_D(x))^{-1})}\, .
$$
\qed
\begin{lemma}\label{l:gamma-1/2}
If $\psi(\lambda)=\lambda^{1/2}$, then there exists 
$C=C(D, \phi)>0$ such that
$$
\E_x \tau_D
\asymp^C
 \Phi(\delta_D(x))^{1/2}\log(1/\delta_D(x))\, ,\qquad x\in D\, .
$$
\end{lemma}
\pf By following the proof of the upper bound in Lemma \ref{l:gamma-small}, we obtain
\begin{align*}
& \E_x \tau_D\le c_1  
\frac{1}{\psi(\Phi(\delta_D(x))^{-1})}
+ c_2\Phi(\delta_D(x))^{1/2}\int_{\delta_D(x)}^{\mathrm{diam}(D)} \frac{1}{\Phi(s)^{1/2} s(\Phi(s)^{-1})^{1/2}}\, ds \\
&=c_1  \frac{1}{\psi(\Phi(\delta_D(x))^{-1})}
+c_2\Phi(\delta_D(x))^{1/2}\log(\mathrm{diam}(D)/\delta_D(x))\\
&\le c_3  \Phi(\delta_D(x))^{1/2}\log(1/\delta_D(x))\, .
\end{align*}

For the lower bound, as one can see from the end of the proof of 
 Lemma \ref{l:gamma-small}, we only need to consider $x$ close to the boundary.
 Since $D$ is $C^{1,1}$, there exists a constant $L>0$  such that 
for every $x \in D$ with $\delta_D(x) \le L/2$, one can find a cone 
$\CC$ with vertex in $x$, pointing inward, 
with height  $L$, and aperture not depending on $x$. 
Moreover, such a cone $\CC$ can be chosen 
so that $\CC\subset \{y\in D: |x-y|\le \delta_D(y)\}$. Then for $y\in \CC\cap \{y\in D:\, \delta_D(x)\le |x-y|\}$, we have that
$$
G^{Y^D}(x,y)\ge c_4 \frac{\Phi(\delta_D(x))^{1/2}}{\Phi(|x-y|)^{1/2}}g(|x-y|)\, .
$$
Hence
\begin{align*}
&\E_x \tau_D \ge c_5 \Phi(\delta_D(x))^{1/2}\int_{\CC\cap \{y\in D:\, \delta_D(x)\le |x-y|\}} \frac{1}{\Phi(|x-y|)^{1/2}|x-y|^d \Phi(|x-y|)^{-1/2}}\\
&\ge c_6 \Phi(\delta_D(x))^{1/2}\int_{\delta_D(x)}^{L}\frac{ds}{s}\ge c_7 \Phi(\delta_D(x))^{1/2}\log(1/\delta_D(x))\, .
\end{align*}

\qed

%%%%%%%%%%%%%%%%%%%%%%%%%%%%%%%%%%%%%%%%%%%%%%%%%%%%%%%%%%%%%%%%%%
\section{Boundary Harnack principle in $C^{1,1}$  open set }\label{s:bhi}
%%%%%%%%%%%%%%%%%%%%%%%%%%%%%%%%%%%%%%%%%%%%%%%%%%%%%%%%%%%%%%%%%%

In this section we assume that
the lower weak scaling index $\gamma_1$ of $\psi$ is strictly larger than $1/2$.
We continue assuming that  $d\ge 2$ and that 
$D\subset \R^d$ is a bounded $C^{1,1}$ open set.
 Again, let $(R, \Lambda)$ be the $C^{1, 1}$ characteristics of $D$.
Without loss of generality we assume that $R \le 1$.

Combining \eqref{e:JY} and \eqref{e:p-r-comparable}--\eqref{e:p-large-time}, we immediately get the following
\begin{prop}\label{p:est-JY} For any $T>0$, there exists $C\ge 1$ such that for all $x, y\in D$,
$$
C^{-1}J^D(x, y)\le J^{Y^D}(x, y)\le CJ^D(x, y), 
$$
where
\begin{align*}
J^D(x, y)
:=&\int_0^T \left(\frac{\Phi(\delta_D(x))^{1/2}}{s^{1/2}}\wedge 1\right)
\left(\frac{\Phi(\delta_D(y))^{1/2}}{s^{1/2}}\wedge 1\right) \\
&\quad \times \left(\Phi^{-1}(s)^{-d}\wedge \frac{s}{|x-y|^d \Phi(|x-y|)}\right)\nu(s)ds+\Phi(\delta_D(x))^{1/2}\Phi(\delta_D(y))^{1/2}\, .
\end{align*}
\end{prop}

Recall $\kappa_0=(1+(1+ \Lambda)^2)^{-1/2} $, and for $Q \in \partial D$,  $V_Q(r)$ is  a $C^{1,1}$ 
domain  with characteristic
$(rR/L, \Lambda L/r)$ such that $D_Q( r/2, r/2) \subset V_Q(r)  \subset  D_Q( r, r) $ where  $L=L(R, \Lambda, d)>0$.
Recall that $g$ is defined in \eqref{e:defg}. 
\begin{lemma}\label{l:ltwdnw} There exists  $C>0$ such that for every $r \le  \kappa_0 R/2$, $Q \in \partial D$ and $x \in D_Q ( r/4 , r/4 )$, 
\begin{equation}\label{e:tww}
\E_x\tau_{V_Q(r)}\le  \int_{V_Q(r)}G^{Y^D}(x,y)\, dy \le C\frac{\Phi(\delta_D(x))^{1/2}}{\Phi(r)^{1/2}\psi(\Phi(r
)^{-1})}\, .
\end{equation}
\end{lemma}
\pf By using Theorem \ref{t:green-function-estimate} in the 
second inequality, 
\eqref{e:Uzr}
 in the third inequality, 
\eqref{e:weak-scaling-psi} in the penultimate inequality, and 
\eqref{e:weak-scaling-phi}
 and $\gamma_1>1/2$ in the last inequality, we get that
\begin{align*}
& \E_x \tau_{V_Q(r)} =\int_{V_Q(r)}G^{Y^D}_{\tau_{V_Q(r)}}(x,y)\, dy \le \int_{V_Q(r)}G^{Y^D}(x,y)\, dy \nn \\
&\le c_{1} \Phi(\delta_D(x))^{1/2} \int_{V_Q(r)} \frac{1}{\Phi(|x-y|)^{1/2} |x-y|^d \psi(\Phi(|x-y|)^{-1})}\, dy \nn \\
&\le 
c_1 \Phi(\delta_D(x))^{1/2}
 \int_{B(x,2r/\kappa_0)} \frac{1}{\Phi(|x-y|)^{1/2} |x-y|^d \psi(\Phi(|x-y|)^{-1})}\, dy \nn \\
 &\le 
 c_2 \Phi(\delta_D(x))^{1/2}
 \int_{0}^{2r/\kappa_0} \frac{1}{\Phi(s)^{1/2} s\psi(\Phi(s)^{-1})}\, ds \nn \\
 &\le 
 c_3 \frac{\Phi(\delta_D(x))^{1/2}}{\Phi(r)^{1/2}\psi(\Phi(r
 )^{-1})}\int_0^{2r/\kappa_0} \left(\frac{\Phi(r)}{\Phi(s)}\right)^{1/2-\gamma_1} \frac{ds}{s}\nn \\
& \le 
c_4 \frac{\Phi(\delta_D(x))^{1/2}}{\Phi(r)^{1/2}\psi(\Phi(r
)^{-1})}\, . 
\end{align*}
\qed
 
 \begin{lemma}\label{L:ggggg}
There exists
 $C>0$ such that 
 for all $r \le  \kappa_0 R/2$, $Q \in D$ and $x \in D_Q(2^{-3}r, 2^{-3}r)$,
 \begin{align*}
&\P_x\Big( Y^D(\tau_{ V_Q(r)}) \in D_Q (  2 r , r)\setminus D_Q (  3 r/2 , r)\Big) 
\ge C \frac{\Phi(\delta_D(x))^{1/2}}{\Phi(r)^{1/2}}.
\end{align*}
 \end{lemma}
 \pf  Without loss of generality, we assume $Q=0$.  Recall that 
$V_0(r) \subset  D_0( r, r) \subset D\cap B(0,  r /\kappa)$. 
Note that, for $y\in D_0 (  2 r , r)\setminus D_0 (  3 r/2 , r)$ and $z\in V_0(r)$, it holds that
$$
|z-y|\asymp^{c_1} r, \quad \delta_D(y) \ge c_2 r \quad \text{and} \quad \delta_D(z) \le c_3 r\, .
$$
Using this and  \eqref{e:Phi(lambda-t)}, 
\eqref{e:v-nu-asymp} and
Proposition \ref{p:est-JY},  
 we have that  for $z \in V_0(r)$, 
 \begin{align}\label{e:Jlow22}
&\int_{D_0 (  2 r , r)\setminus D_0 (  3 r/2 , r)}     J^{Y^D}(z,y)dy\nn\\
&  \ge c_4 \int_{D_0 (  2 r , r)\setminus D_0 (  3 r/2 , r)} 
\int_{\Phi(|y-z|)}^{2\Phi(|y-z|)} \left(\frac{\Phi(\delta_D(z))^{1/2}}{t^{1/2}}\wedge 1\right)\left(\frac{\Phi(\delta_D(y))^{1/2}}{t^{1/2}}\wedge 1\right)\Phi^{-1}(t)^{-d} \nu(t)dt  dy \nn \\
& \asymp^{c_5}  \int_{D_0 (  2 r , r)\setminus D_0 (  3 r/2 , r)}   \left(\frac{\Phi(\delta_D(z))^{1/2}}{\Phi(|z-y|)^{1/2}}\wedge 1\right) \frac{\psi(\Phi(|z-y|)^{-1})}{|z-y|^d} dy \nn\\
& \ge c_6\frac{\Phi(\delta_D(z))^{1/2}}{\Phi(r)^{1/2}}{\psi(\Phi(r)^{-1})}\, .
\end{align}

For $a, b>0$, we define the cone $\sC(x,a ,b)$ above $x$ by
$$
\sC(x,a ,b):=\{y=(\wt y,y_d)\in B(x,a) \mbox{ in }  CS : y_d>x_d ,|\wt x-\wt y|<b (y_d -x_d)\}.
$$
Since $D$ is $C^{1,1}$ and $x \in D_0(2^{-3}r, 2^{-3}r)$, it is easy to see that   
there exists $\varepsilon \in (0, (2(1+ \Lambda))^{-2})$ such that 
\begin{align}
\label{e:Ceta}
\sC(x,2^{-4} r ,\varepsilon)\subset D_0(2^{-2}r, 2^{-2}r). \end{align}
By \eqref{e:Ceta} and the fact that $D_0( r/2, r/2) \subset V_0(r)  \subset  D_0( r, r)$, we have that $\delta_D(z) =\delta_{V_0(r)}(z)$  for all 
$z \in \overline{ \sC(x,2^{-4}  r ,\varepsilon)}$.

It is easy to see that there exists $c_5\in (0,1]$ such that 
\begin{equation}\label{e:D2_421}
c_5|x-z|\le 
\delta_{V_0(r)}(z)=\delta_{D} (z),
\qquad z \in 
\overline{\sC(x,2^{-5} r, 2^{-1}\eps)}.
\end{equation}
We claim that for $z \in \sC(x,2^{-5} r, 2^{-1}\varepsilon) $,
 \begin{align}
\label{e:D2_420}
\left(\frac{\Phi(\delta_{D}(x))}{\Phi(|x-z|)}\wedge 1
\right)\Phi(\delta_{D}(z))  \ge c_6 \Phi(\delta_{D}(x)).
\end{align}
If $z \in \sC(x,2^{-5}r, 2^{-1}\eps)  \setminus  B(x,   \delta_D(x)/2)$,
then $|x-z| \ge \delta_D(x)/2$, 
so by \eqref{e:D2_421} and \eqref{e:Phi(lambda-t)},
\begin{align*}
\left(\frac{\Phi(\delta_{D}(x))}{\Phi(|x-z|)}\wedge 1 \right) \Phi(\delta_{D}(z))  \ge 
c_7\left(\frac{\Phi(\delta_{D}(x))}{\Phi(|x-z|)}\wedge 1\right)\Phi(|x-z|)\ge c_8 \Phi(\delta_{D}(x))\, .
\end{align*}
If $z \in \sC(x,2^{-5}r , 2^{-1}\eps)\cap  B(x,   \delta_D(x)/2)$
then $|x-z| < \delta_D(x)/2$ and so by  \eqref{e:Phi(lambda-t)},
$$ 
\Phi(\delta_{D}(z)) \ge \Phi(\delta_{D}(x) -|x-z|) > \Phi(\delta_{D}(x)/2) \ge 2^{-2}\Phi(\delta_{D}(x)).
$$
Thus,
 \begin{align*}
\left(\frac{\Phi(\delta_{D}(x))}{\Phi(|x-z|)}\wedge 1\right) \Phi(\delta_{D}(z))  \ge \Phi(\delta_{D}(z)) \ge
2^{-2} \Phi(\delta_{D}(x))\, .
\end{align*}
We have proved \eqref{e:D2_420}.

Since $G^{Y^D}_{V_0(r)}\ge G^{Y^{V_0(r)}}$ by \eqref{e:UB-smaller-UDB},  
we have by using Proposition \ref{p:gfcnlb-scaling} in the second below,
\eqref{e:D2_421} in the third,
\eqref{e:D2_420} in the fourth and
\eqref{e:weak-scaling-psi-phi} in the fifth,
\begin{align}
\label{e:D2_41}
& \E_{x}
\int_0^{\tau_{ V_0(r)}} \Phi(\delta_D(Y^D_t))^{1/2} dt =\int_{V_0(r)} G^{Y^D}_{V_0(r)}(x,z) \Phi(\delta_D(z))^{1/2}\, dz \nn \\
&\ge  c_9 \int_{V_0(r)} \left(\frac{\Phi(\delta_{V_0(r)}(x))^{1/2}}{\Phi(|x-z|)^{1/2}}\wedge 1\right) \left(\frac{\Phi(\delta_{V_0(r)}(z))^{1/2}}{\Phi(|x-z|)^{1/2}}\wedge 1\right) g(|x-z|) \Phi(\delta_D(z))^{1/2}  dz \nn \\
&\ge  c_{10} \int_{\sC(x,2^{-5} r ,
\varepsilon/2)} \left(\frac{\Phi(\delta_{D}(x))^{1/2}}{\Phi(|x-z|)^{1/2}}\wedge 1\right) g(|x-z|) \Phi(\delta_D(z))^{1/2}   dz \nn \\
&\ge  c_{11} \Phi(\delta_D(x))^{1/2}  \int_{\sC(x,2^{-5} r ,
\varepsilon/2)} g(|x-z|)   dz \nn \\
&\ge  c_{12}\Phi(\delta_D(x))^{1/2} \int_{0}^{2^{-5} r}
\frac1{s\psi(\Phi(s)^{-1})}ds\ge c_{13}
\frac{\Phi(\delta_D(x))^{1/2}}{\psi(\Phi(r)^{-1})} .
\end{align}
Combining \eqref{e:Jlow22} and \eqref{e:D2_41}, we get
 \begin{align*}
&\P_x\Big( Y^D(\tau_{ V_0(r)}) \in D_0 (  2 r , r)\setminus D_0 (  3 r/2 , r)\Big)  =\E_{x} \int_0^{\tau_{ V_0(r)}} \int_{D_0 (  2 r , r)\setminus D_0 (  3 r/2 , r)}     J^{Y^D}(Y^D_t,y)\, dy\, dt \nn \\
&\ge c_{14} \frac{1}{\Phi(r)^{1/2}} {\psi(\Phi(r)^{-1})} 
\E_{x} \int_0^{\tau_{ V_0(r)}} \Phi(\delta_D(Y^D_t))^{1/2} \, dt \nn \\
&\ge c_{15}
  \frac{1}{\Phi(r)^{1/2}}
  {\psi(\Phi(r)^{-1})}
\frac{\Phi(\delta_D(x))^{1/2}}{\psi(\Phi(r)^{-1})}
=
c_{15} \frac{\Phi(\delta_D(x))^{1/2}}{\Phi(r)^{1/2}} \, .
\end{align*}
 \qed
 
  Let
 \begin{align}
 \label{e:defj}
 j(r):=\frac{\psi(\Phi(r)^{-1})}{r^d}\, , \qquad r>0\, .
 \end{align}
Note that $j$ is a decreasing function and it satisfies the doubling property near 0.
 
\begin{lemma}\label{L:2} Assume that $\gamma_1>1/2$.
There exists  $C>0$ such that for every 
$r \le  \kappa_0 R/2$, $Q \in \partial D$ 
and $x \in D_Q ( 2^{-3} r, 2^{-4} r )$,
  \bee\label{e:L:2}
\P_{x}\left(Y^D( \tau_{ V_Q (r)}) \in
D_Q (  2 r , 2r)\right) \,\le\, C \,  \, 
 \frac{\Phi(\delta_D(x))^{1/2} }{\Phi(r)^{1/2}}.   
 \eee
\end{lemma}

\pf  
Without loss of generality, we assume $Q=0$ and fix $r \le  \kappa^{-1} R/2$.  
Note that, by \eqref{e:Uzr}, 
$V_0(r) \subset  D_0( r, r) \subset D\cap B(0,  r /\kappa_0)$.
Using the L\'evy system formula and \eqref{e:JYDJX},  \eqref{e:exit-time-YD}, \eqref{e:UB-smaller-UDB} we get
\begin{align}
\label{e:ub1}
&\P_{x}\left(Y^D( \tau_{ V_0 (r)}) \in
D_0 (  2 r , r)\setminus D_0 (  3 r/2 , r)\right)\nn\\
=&
\E_{x}
\int_0^{\tau_{ V_0(r)}}
\int_{D_0 (  2 r , r)\setminus D_0 (  3 r/2 , r)}    J^{Y^D}(Y^D_t,z)dzdt\nn\\
=&
\int_{ V_0(r)} G^{Y^D}_{V_0(r)}(x,y)
\int_{D_0 (  2 r , r)\setminus D_0 (  3 r/2 , r)}    J^{Y^D}(y,z)dzdy\nn\\
 \le& 
\int_{D_0 (  2 r , r)\setminus D_0 (  3 r/2 , r)}  \int_{ V_0(r)} G^{Y^D}(x,y)
   J^{X}(y,z)dy dz.
\end{align}
We first note that 
\begin{align}
\label{e:xzr}
c_{0} r <|y-z| \le 4(1+\kappa_0^{-1})r  \normal \quad \text{ for } (y,z) \in (V_0(r) \times( D_0 (  2 r , r)\setminus D_0 (  3 r/2 , r)).
\end{align}
Thus by \eqref{e:JD-estimate}  and Lemma \ref{l:ltwdnw}, 
for $x \in D_0 ( r/4 , r/4 )$ and  $z \in D_0 (  2 r , r)\setminus D_0 (  3 r/2 , r)$,
\begin{align}
\label{e:ub2}
&\int_{ V_0(r)} G^{Y^D}(x,y) J^{X}(y,z)dy \nn\\
&\le c_1  \int_{ V_0(r)}G^{Y^D}(x,y) j(|y-z|) \, dy \nn \\
&\le  c_2  \int_{ V_0(r)}G^{Y^D}(x,y) \, dy \frac{\psi(\Phi(r)^{-1})}{r^d } \le c_3  
\frac{\Phi(\delta_D(x))^{1/2} }{\Phi(r)^{1/2}} r^{-d} .
\end{align}
By \eqref{e:ub1} and   \eqref{e:ub2} we obtain that  for $x \in D_0 ( r/4 , r/4 )$
\begin{align}
\label{e:D2_43nu}
&\P_x\left( Y^D(\tau_{ V_0(r)}) \in D_0 (  2 r , r)\setminus D_0 (  3 r/2 , r)\right) 
 \le c_4  \frac{\Phi(\delta_D(x))^{1/2} }{\Phi(r)^{1/2}}. 
\end{align}
Note that, by the same argument in the proof of  Lemma \ref{l:ltwdnw}, we get 
that for $w \in D_0(r/2, r/2)$ and each $\tilde{r}\in (0,r/2)$, 
\begin{align}
\label{e:D2_4n1}
& \int_{D \cap B(w, \tilde{r})} G^{Y^D}(w,z)dz\le c_5  \frac{\Phi(\delta_D(w))^{1/2}}{\Phi(r)^{1/2}\psi(\Phi(\tilde{r})^{-1})}.
\end{align}
Thus, using \eqref{e:D2_4n1}, we have that for $w \in D_0(r/2, r/2)$ and  $0 <4R_1 \le R_2 <r $,
\begin{align}\label{e:ggggg1}
&\P_w \left(Y^D_{\tau_{D \cap B(w, R_1)}} \in D\setminus B(w, R_2)\right)\nn\\
=& \int_{D \cap B(w, R_1)} G_{D \cap B(w, R_1)}^{Y^D}(w,z) 
\int_{D\setminus B(w, R_2)} J^{Y^D}(y,z) dydz\nn
\\
\le&c_6
\int_{D \cap B(w, R_1)} G^{Y^D}(w,z) \int_{\R^d\setminus B(0, R_2)}  |y|^{-d}\psi(\Phi(|y|)^{-1}) dydz\nn
\\
\le & c_7 \frac{\Phi(\delta_D(w))^{1/2}}{\Phi(R_1)^{1/2}\psi(\Phi(R_1
)^{-1})}
\int_{R_2}^{\infty}
 s^{-1}\psi(\Phi(s)^{-1}) ds 
 \le c_8\frac{\Phi(\delta_D(w))^{1/2}}{\Phi(R_1)^{1/2}\psi(\Phi(R_1)^{-1})}\psi(\Phi(R_2)^{-1}). 
\end{align}

Let 
$$H_2=\{Y^D( \tau_{ V_0(r)}) \in
D_0 (  2 r , 2r)  \}, \quad
H_1=\{Y^D( \tau_{ V_0(r)}) \in
D_0 (  2 r , r)\setminus D_0 (  3 r/2 , r)\}. 
$$
We claim that $\P_x( H_2)\leq c_{9} \P_x( H_1)$ 
for all $r \le  \kappa^{-1} R/2$ and $x \in D_0 ( 2^{-3} r , 2^{-4} r )$. 
Combining this claim with \eqref{e:D2_43nu},
we arrive at the conclusion of the lemma:
$$
\P_x\left( Y^D( \tau_{ V_0(r)}) \in
D_0 (  2 r , 2r) \right)\leq c_{9} \P_x( H_1) \le c_{10}\frac{\Phi(\delta_D(x))^{1/2} }{\Phi(r)^{1/2}}, \quad x \in D_0 ( 2^{-3} r , 2^{-4} r )  .
$$

Now we give the proof of the claim, which is inspired by the proof of \cite[Lemma 5.3]{G} 
(see also \cite[Lemma 6.2]{KSV16}).
Note that, by Lemma \ref{L:ggggg},  for $w \in D_0(r/8, r/8)$,
\begin{align}\label{e:D2_43n2}
&\P_w ( H_1 ) \ge c_{11}  \frac{\Phi(\delta_D(w))^{1/2}}{\Phi(r)^{1/2}}.
\end{align}

For $i\geq 1$, set
$$J_i=D_0(2^{-i-2}r,s_i)\setminus  D_0(2^{-i-3}r,s_i),\ \ \ \ s_i=\frac{1}{4}\left(\frac{1}{2}-\frac{1}{50}\sum_{j=1}^i\frac{1}{j^2}\right)r,$$
and $s_0=s_1$.  
Note that $r/(10)<s_i <r/8$.
Define for $i\geq 1$ ,
\begin{align}\label{q1q}
d_i=d_i(r)=
\sup_{z\in J_i}\frac{\P_z( H_2)}{\P_z( H_1)},\ \ \
\widetilde{J}_i=D_0(2^{-i-2}r,s_{i-1}),\ \ \ \
\tau_i=\tau_{\widetilde{J}_i}.
\end{align}
By \eqref{e:D2_43n2}, $\sup_{r \le \kappa^{-1} R/2} d_i(r)$ is finite for all $i \ge 1$.
Repeating the argument leading to \cite[(6.29)]{KSV16}, we get that 
 for $z\in J_i$ and $i\geq 2$,
\begin{align}\label{q12q}  \P_z( H_2) 
\leq    \left ( \sup_{1\leq k\leq i-1}  d_k \right) \P_{z}(H_1) +\P_z\left( Y^{D}_{\tau_i}\in
D_0 (  2 r , 2r) \setminus \cup_{k=1}^{i-1}J_k \right).
\end{align}

For  $i\ge 2$, define $\sigma_{i,0}=0, \sigma_{i,1}=\inf\{t>0: |Y^{D}_t-Y^{D}_0|\geq 2^{-i-2}r\} $ and
$\sigma_{i,m+1}=\sigma_{i,m}+\sigma_{i,1}\circ\theta_{\sigma_{i,m}}$
for $m\geq 1$. By Lemma \ref{lower bound},  we have that there exists $k_1 \in (0,1)$ such that 
\begin{align}\label{e:rs}
\P_{w}(Y^{D}_{\sigma_{i,1}}\in \widetilde{J}_i) \le 1-
\P_{w}( \sigma_{i,1}=\zeta ) \le 1-\P_{w}(
\tau_{B(w,\delta_D(w)/2)}
=\zeta )  <k_1,\ \ \ w\in  \widetilde{J}_i.
\end{align}
For the purpose of further estimates, we now choose a
positive integer $l\ge b_2a_2^{\gamma_2}$ 
such that $k_1^l\le 2^{-2\gamma_2\delta_2}$, where $a_2, \delta_2$ and $b_2, \gamma_2$ 
are the constants from 
\eqref{e:weak-scaling-phi}  and \eqref{e:weak-scaling-psi} respectively.
Next we choose 
$i_0 \ge 2$ large enough so that
$2^{-i}<1/(200 l i^3)$ for all $i\ge i_0$.
Now we assume $i\ge i_0$.
Using \eqref{e:rs} and the strong Markov property we have
that  for   $z\in J_i$,
\begin{align}\label{q132q}    
&\P_z( \tau_{i}>\sigma_{i,li})\leq \P_z(Y^{D}_{\sigma_{i,k}}\in \widetilde{J}_i, 1\leq k\leq
li
 )\nonumber\\
 &=
 \E_z \left[ \P_{Y^{D}_{\sigma_{i,li-1}}} (Y^{D}_{\sigma_{i,1}}\in  \widetilde{J}_i) : Y^{D}_{\sigma_{i,li-1}} \in  \widetilde{J}_i,  Y^{D}_{\sigma_{i,k}}\in \widetilde{J}_i, 1\leq k\leq
li-2
  \right]
 \nonumber\\      
 &\leq \P_z\left(Y^{D}_{\sigma_{i,k}}\in \widetilde{J}_i, 1\leq k\leq
li-1
 \right)k_1\leq k_1^{li}.
\end{align}
Note that if $z\in J_i$ and $y\in D_0 (  2 r , 2r) \setminus[ \widetilde{J}_i \cup(\cup_{k=1}^{i-1}J_k)]$, 
then $|y-z|\ge (s_{i-1}-s_i) \wedge (2^{-3}-2^{-i-2}) r = r/(200 i^2)$.
Furthermore, 
since  $2^{-i-2} r< r/(200 i^2)$ (recall that $i\ge  i_0$), 
if $Y^D_{\tau_i}(\omega)\in D_0(2r,2r)\setminus \cup_{k=1}^{i-1}J_k$ and $\tau_i(\omega)\le \sigma_{i,li}(\omega)$, then $\tau_i(\omega)=\sigma_{i,k}(\omega)$ for some $k=k(\omega)\le li$.  
Dependence of $k$ on $\omega$ will be omitted  in the next few lines. 
Hence on
$\{Y^{D}_{\tau_{i}}\in D_0 (  2 r , 2r) \setminus \cup_{k=1}^{i-1}J_k,\ \ \tau_{i}\leq
\sigma_{i,li}\}$ with $Y^{D}_0=z\in J_i$,  we have
$|Y^D_{\sigma_{i,k}}-Y^D_{\sigma_{i,0}}|=|Y^D_{\tau_i}-Y^D_0|>
\frac{1}{200i^2}r$ for some $1\leq k\leq li$.
Thus  for some $1\leq k\leq li$,
$$
\sum_{j=0}^k|Y^D_{\sigma_{i,j}}-Y^D_{\sigma_{i,j-1}}|>
\frac{1}{200i^2}r
$$
which implies for some $1\leq j\leq k\le  li$,
$$ 
\big|Y^{D}_{\sigma_{i,j}}-Y^{D}_{\sigma_{i,j-1}}\big|\geq
\frac1k\frac{1}{200i^2}r\ge \frac{1}{ li} \frac{1}{200i^2}r\, .
$$
Thus, we have 
\begin{align}& \{Y^{D}_{\tau_{i}}\in D_0 (  2 r , 2r) \setminus
\cup_{k=1}^{i-1}J_k,\ \ \tau_{i}\leq \sigma_{i,li}\}\nonumber\\
  \subset &   \cup_{j=1}^{li}\{|Y^D_{\sigma_{i,j}}- Y^D_{\sigma_{i,j-1}}|\geq r/(400li^3),Y^D_{\sigma_{i,j-1}}\in \widetilde{J}_{i}
\}.
\end{align}
Now, using  \eqref{e:ggggg1} 
 (noting that  $4 \cdot 2^{-i-2}<1/(400 l i^3)$ for all $i\ge i_0$), 
 and repeating the argument leading to \cite[(6.34)]{KSV16}, we get
\begin{align}\label{dsa} 
&\P_z \left(Y^{D}_{\tau_{i}}\in D_0 (  2 r , 2r ) \setminus
\cup_{k=1}^{i-1}J_k,\ \ \tau_{i}\leq \sigma_{i,li} \right)\nonumber\\
 \leq  &li \sup_{z\in \widetilde{J}_{i}} 
  \P_z\left(|Y^D_{\sigma_{i,1}}-z  |\geq r/(400li^3)\right)
\le c_{12}li  
\frac{\psi(\Phi( r/(li^3))^{-1}) }{\psi(\Phi(2^{-i}r)^{-1})}.
\end{align}

By (\ref{q132q}) and  (\ref{dsa}),  we have for 
$z\in J_i$, 
\begin{align}\label{as}
&\P_z\left( Y^D_{\tau_{i}}\in D_0 (  2 r , 2r) \setminus \cup_{k=1}^{i-1}J_k \right) \leq k_1^{li} 
+c_{12} li  \frac{\psi(\Phi( r/(li^3))^{-1})}{\psi(\Phi(2^{-i}r)^{-1})}. 
\end{align} 
By our choice of $l$, using \eqref{e:weak-scaling-psi-phi} we have 
\begin{align}\label{e:ggggg2}
&li    \frac{\psi(\Phi( r/(li^3))^{-1}) }{\psi(\Phi(2^{-i}r)^{-1})}  
\ge a_2^{-\gamma_2} b_2^{-1}li (li^3/2^i)^{2\gamma_2\delta_2}
= a_2^{-\gamma_2} b_2^{-1}l^{1+2\gamma_2 \delta_2} i^{1+6\gamma_2 \delta_2}  (2^{-2\gamma_2 \delta_2})^{i} 
 \ge (k_1^{l})^{i}. 
\end{align}
Thus 
combining \eqref{e:ggggg2} with \eqref{as}, and then 
 using \eqref{e:D2_43n2}, \eqref{e:weak-scaling-phi}, \eqref{e:weak-scaling-psi} and the fact that $ \gamma_1>1/2$, we get that  for
$z\in J_i$,
\begin{align}\label{as1}
&\frac{\P_z( Y^{D}_{\tau_i}\in
D_0 (  2 r , 2r) \setminus \cup_{k=1}^{i-1}J_k)}{\P_z(H_1)} \le 
c_{13} li  
\frac{\Phi(r)^{1/2}}{\Phi(2^{-i}r )^{1/2}}
  \frac{\psi(\Phi( r/(li^3))^{-1}) }{\psi(\Phi(2^{-i}r
)^{-1})}\nn\\
&\le  
c_{14} i  
\frac{\Phi(r)^{1/2}}{\Phi(2^{-i}r )^{1/2}}
  \frac{\Phi(2^{-i}r
)^{\gamma_1}}{\Phi( r/(li^3))^{\gamma_1} }= 
c_{14}i  \frac{\Phi(r)^{1/2}}{\Phi( r/(li^3))^{1/2} }
  \frac{\Phi(2^{-i}r
)^{\gamma_1-1/2}}{\Phi( r/(li^3))^{\gamma_1-1/2} }  \nn \\
&\le  c_{15} i   (li^3)^{\delta_2}
(2^{-i}li^3)^{(2\gamma_1-1)\delta_1}=
c_{16} i^{1+3\delta_2+ 3(2\gamma_1-1)\delta_1}2^{-i(2\gamma_1-1)\delta_1}.
\end{align} 
By this and (\ref{q12q}), 
for $z\in J_i$,
\begin{align}  
&\frac{\P_z( H_2)}{\P_z(H_1)} \leq  \sup_{1\leq k\leq i-1}  d_k   +\frac{\P_z( Y^{D}_{\tau_i}\in
D_0 (  2 r , 2r) \setminus \cup_{k=1}^{i-1}J_k)}{\P_z(H_1)} \nn\\
\leq& \sup_{1\leq k\leq i-1}d_k+
c_{16} \frac{i^{1+3\delta_2+ 3(2\gamma_1-1)\delta_1}}{2^{i(2\gamma_1-1)\delta_1}} \nonumber \\
\le & \sup_{1\leq k\leq i-1}d_k+c_{16}
\frac{
i^{10}}{2^{i(2\gamma_1-1)\delta_1}}.
\nonumber
\end{align}
This implies that
\begin{align} 
d_i& \leq  \sup_{1\leq k\leq i_0-1} d_k
+c_{16}\sum_{k=1}^i\frac{
k^{10}}{2^{k(2\gamma_1-1)\delta_1}}
\leq
\sup_{1\leq k\leq i_0-1 \atop r \le \kappa^{-1} R/2}d_k (r)
+c_{16}\sum_{k=1}^{\infty}
\frac{k^{10}}{2^{k(2\gamma_1-1)\delta_1}}
=:c_{17} <\infty.\nonumber
\end{align}
Thus the claim above is valid, since $D_0 ( 2^{-3} r , 2^{-4} r ) 
\subset \cup_{k=1}^\infty J_k$. The proof is now complete.
\qed

We are now ready to prove the boundary Harnack principle with explicit decay rate near the boundary of $D$
when the lower weak scaling index $\gamma_1$ of $\psi$ is strictly larger than $1/2$.
\begin{thm}\label{t:main-bhp}
Assume that $\gamma_1>1/2$.
Let $D\subset\R^d$ be a bounded $C^{1,1}$ open set with $C^{1,1}$ characteristics  $(R,\Lambda)$.
There exists a constant $C=C(d, \Lambda, R, \phi, \psi)>0$ 
such that for any $r \in (0, R]$, $Q\in \partial D$, and any non-negative function $f$ in $D$ which is harmonic  in $D \cap B(Q, r)$ with respect to $Y^D$ and vanishes continuously on $ \partial D \cap B(Q, r)$, we have
\begin{equation}\label{e:bhp_m}
\frac{f(x)}{\Phi(\delta_D(x))^{1/2}}\,\le C\,\frac{f(y)}{\Phi(\delta_D(y))^{1/2}} \qquad
\hbox{for all } x, y\in  D \cap B(Q, r/2).
\end{equation}
\end{thm}

\noindent
\pf
 In this proof, the constants $\eta$ and $c_i$
are always independent of $r$.

Note that, 
  since $D$ is a $C^{1,1}$ 
  open set and $r<R$, by Theorem \ref{HP2}, it
suffices to prove \eqref{e:bhp_m} for $x,y \in D \cap B(Q,2^{-7} \kappa_0 r)$. Throughout the remainder of the proof 
we assume that  $x \in D \cap B(Q,2^{-7} \kappa_0 r)$.

Let $Q_x$ be
the point $Q_x \in \partial D$ so that $|x-Q_x|=\delta_{D}(x)$ and
let $x_0:=Q_x+\frac{r}{8}(x-Q_x)/|x-Q_x|$. We choose a
$C^{1,1}$ function $\varphi: \bR^{d-1}\to \bR$ satisfying $\varphi
(\wt 0)= 0$, $\nabla\varphi (\wt 0)=(0, \dots, 0)$, $\| \nabla \varphi
 \|_\infty \leq \Lambda$, $| \nabla \varphi (\wt y)-\nabla \varphi (\wt z)|
\leq \Lambda |\wt y-\wt z|$, and an orthonormal coordinate system $CS$ with
its origin at $Q_x$ such that
$$
B(Q_x, R)\cap D=\{ y=(\wt y, y_d) \in B(0, R) \mbox{ in } CS: y_d >
\varphi (\wt y) \}.
$$
In the coordinate system $CS$ we have $\wt x = \wt 0$ and $x_0=(\wt
0, r/8)$. For any $b_1, b_2>0$, we define
$$
\widehat{D}(b_1, b_2):=
\left\{ y=(\wt y, y_d) \mbox{ in } CS: 0<y_d-\varphi(\wt
y)< 2^{-2}\kappa_0 r b_1, \ |\wt y| <  2^{-2}\kappa_0 r b_2\right\}.
$$

By \eqref{e:Uzr}, 
we have that
$\widehat{D}(2, 2)\subset D\cap B(Q_x, r/2)\subset D\cap B(Q, r)$.
Thus, since
$f$ is
harmonic in $D\cap B(Q,r)$ and vanishes continuously in
$\partial D\cap B(Q, r)$, by Lemma \ref{l:regularity}, $f$ is regular
harmonic in 
$\widehat{D}(2, 2)$ and vanishes continuously in
$\partial D\cap \widehat{D}(2, 2)$.

Recall that $V(1):=V_{Q_x}(2^{-2}\kappa_0 r)$ is a $C^{1,1}$
domain   with $C^{1,1}$ characteristics $(rR/L, \Lambda L/r)$
such that $\widehat{D}( 1/2, 1/2) \subset V(1)  \subset \widehat{D}( 1, 1) $,
where $L=L(R, \Lambda, d)>0$.

By Theorem \ref{uhp}  (or Theorem \ref{HP2})
and Lemma \ref{L:ggggg}, we have
\begin{align}\label{e:BHP2}
&f(x) =  \E_x\left[f\big(Y^D(\tau_{ V(1)})\big)\right] \ge
\E_x\left[f\big(Y^D(\tau_{ V(1)})\big); Y^D_{ \tau_{ V(1)}} \in  \widehat{D}(2,1)\setminus \widehat{D}(3/2,1)\right] \nn \\
&\ge c_{16} f(x_0) \P_x\Big( Y^D(\tau_{ V(1)}) \in \widehat{D}(2,1)\setminus \widehat{D}(3/2,1)\Big) 
\ge  c_{17} \frac{\Phi(\delta_D(x))^{1/2}}{\Phi(r)^{1/2}} f(x_0)\, .
\end{align}

Take $w=(\wt{0}, 2^{-6}\kappa_0 r)$. Then there exists 
$\epsilon\in (0,1/8)$ such that 
$$
B(w, \epsilon 2^{-5}\kappa_0 r)\subset \widehat{D}(1/2,1/2)\subset V(1).
$$ 
 Hence 
\begin{align}\label{e:f-lower-bound}
& f(w) \ge \E_w\left[f \left(Y^D(\tau_{ V(1)})\right); \,Y^D(\tau_{ V(1)})\notin  \widehat{D}(2,2)\right] \nn \\
& = \E_w \int_0^{\tau_{ V(1)}} \int_{D\setminus \widehat{D}(2,2)}     J^{Y^D}(Y^D_t, y)f(y)\, dy\, dt \nn \\
& \ge \E_w \int_0^{\tau_{ B(w, \epsilon 2^{-5}\kappa_0 r)}} \int_{D\setminus \widehat{D}(2,2)}     J^{Y^D}(Y^D_t, y)f(y) \, dy\, dt \nn \\
&\ge  c_{18} \E_w \tau_{B(w, \epsilon 2^{-5}\kappa_0 r)} 
\int_{D\setminus \widehat{D}(2,2)} J^{Y^D}(w,y)f(y)\, dy \nn \\
&\ge c_{19} \frac{1}{\psi(\Phi(r)^{-1})}  \int_{D\setminus \widehat{D}(2,2)}  
J^{Y^D}(w,y)f(y)\, dy\,,
\end{align}
 where in the fourth line we used Proposition \ref{p:estimate-of-J-away}, and 
in the last line we used Lemma \ref{l:tauB} and \eqref{e:weak-scaling-psi-phi}.

Further, note that for any $z\in V(1)$ and $y\in D\setminus \widehat{D}(2,2)$ we have that $\delta_D(z)\le c_{23} r \le c_{24}\delta_D(w)$ and $|z-y|\asymp^{c_{25}} |w-y|$. By using these two observations and Proposition  \ref{p:est-JY} we see that
\begin{align}
&J^{Y^D}(z,y)\le c_{26}\left( 
\int_0^1
 \left(\frac{\Phi(\delta_D(z))^{1/2}}{s^{1/2}}\wedge 1\right)\left(\frac{\Phi(\delta_D(y))^{1/2}}{s^{1/2}}\wedge 1\right)\right. \nn\\
&\quad \left. \times \left(\Phi^{-1}(s)^{-d}\wedge \frac{s}{|z-y|^d \Phi(|z-y|)}\right)\nu(s)ds+\Phi(\delta_D(z))^{1/2}\Phi(\delta_D(y))^{1/2}\right) \nn\\
&\le c_{27}\left( 
\int_0^1 \left(\frac{\Phi(\delta_D(w))^{1/2}}{s^{1/2}}\wedge 1\right)\left(\frac{\Phi(\delta_D(y))^{1/2}}{s^{1/2}}\wedge 1\right) \right.\nn\\
&\quad \left. \times 
\left(\Phi^{-1}(s)^{-d}\wedge \frac{s}{|w-y|^d \Phi(|w-y|)}\right)\nu(s)ds
+\Phi(\delta_D(w))^{1/2}\Phi(\delta_D(y))^{1/2}\right)\nn \\
&\le c_{28}J^{Y^D}(w,y)\, . \label{e:ffff1}
\end{align}
Hence, combining  Lemma \ref{l:ltwdnw} with \eqref{e:f-lower-bound}--\eqref{e:ffff1} we now have
\begin{align}\label{e:D2_6-not}
&\E_{x}\left[f \left(Y^D(\tau_{ V(1)})\right); \, Y^D(\tau_{ V(1)}) \notin  \widehat{D}(2,2)\right] \nn \\
&=\E_x \int_0^{\tau_{ V(1)}} \int_{D\setminus \widehat{D}(2,2)} J^{Y^D}(Y_t^D,y)f(y)\, dy\, dt \nn \\
&\le c_{28} \E_x \tau_{V(1)} \int_{D\setminus \widehat{D}(2,2)} J^{Y^D}(w,y)f(y)\, dy \nn \\
&\le c_{29}  \frac{\Phi(\delta_D(x))^{1/2}}{\Phi(r)^{1/2}\psi(\Phi(r
)^{-1})} \int_{D\setminus \widehat{D}(2,2)}  J^{Y^D}(w,y) f(y)\, dy \nn \\
&\le c_{30} \frac{\Phi(\delta_D(x))^{1/2}}{\Phi(r)^{1/2}} f(w)\, .
\end{align} 

On the other hand, by 
Theorems \ref{HP2}, 
\ref{t:carleson} and Lemma \ref{L:2}, we have
\begin{align}\label{e:D2_6}
&\E_x\left[f\left(Y^D(\tau_{ V(1)})\right);\,  Y^D(\tau_{V(1)}) \in \widehat{D}(2,2)\right] \nn \\
& \le\, 
c_{31} \, f(x_0) \P_x\left(Y^D(\tau_{V(1)}) \in  \widehat{D}(2,2)\right)\le\, c_{32} \, f(x_0)\frac{\Phi(\delta_D(x))^{1/2}}{\Phi(r)^{1/2}}\, .
\end{align}
Combining \eqref{e:D2_6-not} and \eqref{e:D2_6} and using 
Theorems \ref{HP2} and \ref{t:carleson} again, 
we get
 \begin{align} \label{e:BHP1}
&f(x) =  \E_x\left[f (Y^D(\tau_{ V(1)})); \,Y^D(\tau_{ V(1)}) \in \widehat{D}(2,2)\right]+\,\E_x\left[ f(Y^D(\tau_{ V(1)})); \, Y^D(\tau_{ V(1)}) \notin  \widehat{D}(2,2)\right] \nn  \\
&\le  c_{33}\left(\frac{\Phi(\delta_D(x))^{1/2}}{\Phi(r)^{1/2}} f(x_0) + \frac{\Phi(\delta_D(x))^{1/2}}{\Phi(r)^{1/2}} f(w) \right) \le c_{34} \frac{\Phi(\delta_D(x))^{1/2}}{\Phi(r)^{1/2}} f(x_0)\, .
\end{align}
Together with \eqref{e:BHP2} we get that
\begin{equation}\label{e:BHP-yes}
f(x)\asymp^{c_{35}} \frac{\Phi(\delta_D(x))^{1/2}}{\Phi(r)^{1/2}} f(x_0)\, .
\end{equation}
For any $y\in D \cap B(Q,2^{-7} \kappa_0 r)$, we have the same estimate with $f(y_0)$ instead of $f(x_0)$ where $y_0=Q_y+\frac{r}{8}(y-Q_y)/|y-Q_y|$ and $Q_y\in \partial D$ with $|y-Q_y|=\delta_{D}(y)$. 
 Since $D$ is $C^{1,1}$, 
 using Theorem \ref{HP2}, 
$f(y_0)\asymp^{c_{36}} f(x_0)$. It follows therefore from \eqref{e:BHP-yes} that for every $x,y \in D \cap B(Q,2^{-7} \kappa_0 r)$,
$$
\frac{f(x)}{f(y)}\,\le \,
c_{37}\,\frac{\Phi(\delta_D(x))^{1/2}}{\Phi(\delta_D(y))^{1/2}}\, ,
$$
which  proves the theorem.
\qed

%%%%%%%%%%%%%%%%%%%%%%%%%%%%%%%%%%%%%%%%%%%%%%%%%%%%%%%%%%%%%%%%%%%%%%%%%%%%%%%%%%
%%%%%%%%%%%%                                                Jumping kernel estimates                          %%%%%%%%%%%%%%%%%%%%%%%%%%%%%%%%%%
%%%%%%%%%%%%%%%%%%%%%%%%%%%%%%%%%%%%%%%%%%%%%%%%%%%%%%%%%%%%%%%%%%%%%%%%%%%%%%%%%%
\section{Jumping  kernel estimates}\label{s:j}

In this section we continue assuming that 
$d\ge 2$ and that $D$ is a bounded
$C^{1,1}$ open set in  $\R^d$ 
with $C^{1, 1}$ characteristics $(R, \Lambda)$.

The goal of the section is to derive sharp two-sided estimates for the jumping kernel $J^{Y^D}$. Somewhat surprisingly, 
obtaining explicit bounds not involving integral terms
does not seem to be possible without additional assumptions on the Laplace exponent $\psi$. In case when $\psi(\lambda)=\lambda^{\gamma}$, this can be explained by 
the different qualitative boundary behaviors
of $J^{Y^D}$ in cases 
$\gamma\in (0,1/2)$, $\gamma=1/2$ and $\gamma\in (1/2,1)$, 
cf. Example \ref{e:stable}.

Let 
$$
\theta(t):=\Phi(t)\psi(\Phi(t)^{-1}) \quad \text{and} \quad   \eta(t):=\Phi(t)^{1/2}\psi(\Phi(t)^{-1})\, ,\quad t>0\, .
$$
 It follows  that,  as a composition of two increasing functions, cf.~\eqref{e:t-psi-t-1}, 
$\theta(t)$  is an increasing function.

It is also straightforward to see that
\begin{equation}\label{e:Phi(theta-t)}
(1\wedge \lambda^2)\theta(t)\le \theta(\lambda t)\le (1\vee \lambda^2) \theta(t)\, , \qquad \lambda, t>0\, .
\end{equation}

We will say that a function $f:(0,\infty)\to [0,\infty)$ is \emph{almost increasing} near 0 if for every $T>0$ there exists a constant $C=C(T)>0$ such that $f(s)\le C f(t)$ 
for all $0<s<t\le T$. An almost decreasing function is defined analogously.

Recall that $r(t,x,y)$, $g$ and $j$ are defined in \eqref{e:rtxy},  \eqref{e:defg} and \eqref{e:defj} respectively. 
\begin{lemma}\label{l:jumping-main-integral}
\begin{itemize}
\item[\rm (i)]
Suppose that $r\mapsto r^{1/2} \psi (r^{-1})$ is almost decreasing 
near $0$ and that for each $T>0$ there is 
a constant  $C_4=C_4(T, \psi)>0$  such that
\bee\label{e:6.9}
\int_r^T  s^{-1/2} \psi (s^{-1}) \, ds \leq  C_4  \, r^{1/2} \psi (r^{-1}) \qquad \hbox{for }r\in (0, T].
\eee
Then there exists 
$C\ge 1$ such that for all $x,y\in D$,
$$
\int_0^{\Phi(|x-y|)}r(t,x,y)\nu(t)\, dt \asymp^C
\frac{ \psi(\Phi(|x-y|)^{-1})}{|x-y|^d} \left(    \frac{ \theta (\delta_D(x) \wedge \delta_D(y))}{ \theta(|x-y|)}  \wedge 1\right)\, .
$$

\item[\rm (ii)]
Suppose that $r\mapsto r^{1/2} \psi (r^{-1})$ is almost increasing 
near $0$ and that  for every $T>0$, there is 
a constant  $C_5=C_5(T, \psi)>0$ such that
\bee\label{e:6.10}
\int_0^r s^{-1/2} \psi (s^{-1})  ds 
\le C_5  r^{1/2} \psi (r^{-1}) \qquad \hbox{for every } r\in (0, T].
\eee
Then there exists $C\ge 1$ such that for all $x,y\in D$,
\begin{align*}
&\int_0^{\Phi(|x-y|)}r(t,x,y)\nu(t)\, dt\\
&\asymp^C \frac{ \psi(\Phi(|x-y|)^{-1}
)}{|x-y|^d} 
 \left(    \frac{ \Phi ( \delta_D(x) \wedge \delta_D(y))^{1/2}  }{ \Phi(|x-y|)^{1/2}}\wedge 1\right)
 \left(    \frac{ \eta ( \delta_D(x) \vee \delta_D(y)) }{ \eta(|x-y|) }\wedge 1\right)\, .
\end{align*}
\end{itemize}
\end{lemma}
\begin{remark}\label{r:moved}
{\rm
It is easy to see that if $\psi$ satisfies \eqref{e:weak-scaling-psi} with $\gamma_1>1/2$, then the assumptions in 
Lemma \ref{l:jumping-main-integral}(i) hold true. Similarly, if $\psi$ satisfies \eqref{e:weak-scaling-psi} with $\gamma_2<1/2$, then the assumptions in Lemma \ref{l:jumping-main-integral}(ii) are true. 
}
\end{remark}

\noindent {\bf Proof of Lemma \ref{l:jumping-main-integral}} 
Put $T:= 2  \Phi({\rm diam} (D))$
By using \eqref{e:v-nu-asymp},  \eqref{e:rtxy} and \eqref{e:Phi-less-x-y} we see that
\begin{eqnarray}\label{e:jumping-main-integral-upper}
\lefteqn{\int_0^{\Phi(|x-y|)}r(t,x,y)\nu(t)\, dt } \nn \\
&\asymp&\frac{1}{|x-y|^d \Phi(|x-y|)} \int_0^{\Phi(|x-y|)}
 \left(\frac{\Phi(\delta_D(x))^{1/2}}{t^{1/2}}\wedge 1\right) \left(\frac{\Phi(\delta_D(y))^{1/2}}{t^{1/2}}\wedge 1\right) 
\psi(t^{-1})\, dt \nn \\
&=:&\frac{1}{|x-y|^d \Phi(|x-y|)} \, I\, .
\end{eqnarray}
Let 
$$
I_1:=\int_0^{ \Phi(\delta_D(x))  }
 \psi(t^{-1})\, dt, \quad 
 I_2:=\Phi (\delta_D(x))^{1/2} 
 \int_{\Phi(\delta_D(x))}^{\Phi(\delta_D(y))} t^{-1/2} \psi(t^{-1})\, 
dt,
$$
$$
I_3:=\Phi (\delta_D(x))^{1/2} \Phi (\delta_D(y))^{1/2} 
\int_{\Phi(\delta_D(y))}^{\Phi(|x-y|)} t^{-1} \psi(t^{-1})\, dt
$$
and
$$
I_4:=\Phi (\delta_D(x))^{1/2} \int_{\Phi(\delta_D(x))}^{\Phi(|x-y|)}
 t^{-1/2} \psi(t^{-1})\, dt. 
$$

\medskip
\noindent
\textbf{Upper bound:}
Without loss of generality we assume $\delta_D(x) \le \delta_D(y)$ and consider three cases:

\noindent
(1) $\delta_D(x) \le \delta_D(y) \le |x-y|$:  Then
$I=I_1+I_2+I_3.$

By \eqref{e:weak-scaling-psi}, we have
 $
  I_1\asymp^{c_1} \psi(\Phi (\delta_D(x))^{-1}) \Phi (\delta_D(x))
 $
 and
 $$I_3
\le c_2 \Phi (\delta_D(x))^{1/2} \Phi (\delta_D(y))^{1/2}
 \psi(\Phi (\delta_D(y))^{-1}). 
 $$
In case (i), by \eqref{e:6.9},
$$
I_2 \le \Phi (\delta_D(x))^{1/2} \int_{\Phi(\delta_D(x))}^{2T} t^{-1/2} \psi(t^{-1})\, dt \le c_3 \psi(\Phi (\delta_D(x))^{-1}) \Phi (\delta_D(x))\, .
$$
By using that $r^{1/2} \psi (r^{-1})$ is almost decreasing 
near $0$ we see that 
\begin{align}\label{e:case-i-1}
I &\le  c_4 
\left[\psi(\Phi (\delta_D(x))^{-1}) \Phi (\delta_D(x))+ \Phi (\delta_D(x))^{1/2} \Phi (\delta_D(y))^{1/2}
\psi(\Phi (\delta_D(y))^{-1})
\right] \nn \\
 & \le c_5 \psi(\Phi (\delta_D(x))^{-1}) \Phi (\delta_D(x)) =c_5 \psi(\Phi (\delta_D(x) \wedge \delta_D(y))^{-1}) \Phi (\delta_D(x) \wedge \delta_D(y))\, .
 \end{align}
In case (ii), by \eqref{e:6.10}, we have
$$I_2 \le \Phi (\delta_D(x))^{1/2} \int_{0}^{\Phi(\delta_D(y))} t^{-1/2} \psi(t^{-1})\, dt \le c_6 \Phi (\delta_D(x))^{1/2} \Phi (\delta_D(y))^{1/2}
 \psi(\Phi (\delta_D(y))^{-1}).$$
By using that $r^{1/2} \psi (r^{-1})$ is almost increasing 
near $0$ we see that
\begin{align}\label{e:case-ii-1}
 I
 &\le
  c_7 
  \left[\psi(\Phi (\delta_D(x))^{-1}) \Phi (\delta_D(x))+ \Phi (\delta_D(x))^{1/2} \Phi (\delta_D(y))^{1/2}\psi(\Phi (\delta_D(y))^{-1})
  \right] \nn\\
  & \le c_8\Phi (\delta_D(x) \wedge \delta_D(y))^{1/2} \Phi (\delta_D(x) \vee \delta_D(y))^{1/2}
 \psi(\Phi (\delta_D(x) \vee \delta_D(y))^{-1})\, .
 \end{align}

 \noindent
(2) $\delta_D(x)  \le |x-y| \le \delta_D(y)$: Then
$I=I_1+I_4 \, .$

In case (i), by \eqref{e:6.9}, we have
$$
I_4 \le \Phi (\delta_D(x))^{1/2} \int_{\Phi(\delta_D(x))}^{2T} t^{-1/2} \psi(t^{-1})\, dt \le c_9 \psi(\Phi (\delta_D(x))^{-1}) \Phi (\delta_D(x))\, .
$$
Thus  we see that 
\begin{align}\label{e:case-i-2}
&I  \le c_{10} \psi(\Phi (\delta_D(x))^{-1}) \Phi (\delta_D(x)) =c_{10} \psi(\Phi (\delta_D(x) \wedge \delta_D(y))^{-1}) \Phi (\delta_D(x) \wedge \delta_D(y))\, .
\end{align}

In case (ii), by \eqref{e:6.10}, we have
$$
I_4 \le \Phi (\delta_D(x))^{1/2} \int_{0}^{\Phi(|x-y|)} t^{-1/2} \psi(t^{-1})\, dt \le c_{11} \Phi (\delta_D(x))^{1/2}\Phi(|x-y|)^{1/2}
 \psi(\Phi(|x-y|)^{-1})\, .
 $$
By using that $r^{1/2} \psi (r^{-1})$ is almost increasing 
near $0$, we see that
\begin{align}\label{e:case-ii-2}
I &\le c_{12} (\psi(\Phi (\delta_D(x))^{-1}) \Phi (\delta_D(x))+ \Phi (\delta_D(x))^{1/2}  \Phi(|x-y|)^{1/2} \psi(\Phi(|x-y|)^{-1})) \nn \\
 & \le c_{13}\Phi (\delta_D(x))^{1/2}  \Phi(|x-y|)^{1/2} \psi(\Phi(|x-y|)^{-1})\, .
\end{align}

 \noindent
(3) $ |x-y| \le \delta_D(x)  \le \delta_D(y)$: Then by \eqref{e:weak-scaling-psi},
 \begin{align}\label{e:case-i-ii-3}
&I= \int_0^{\Phi(|x-y|)}
 \psi(t^{-1})\, dt  \asymp^{c_{14}} \Phi(|x-y|)
 \psi(\Phi(|x-y|)^{-1}) .
 \end{align}

\medskip
\noindent 
\textbf{Lower bound:}
Again, we assume $\delta_D(x) \le \delta_D(y)$.
Let $M :=(2/a_1(T^{-2}))^{1/(2 \delta_1)} \vee 2 $ 
(where $a_1(T^{-2})$ is the constant in the extended version of \eqref{e:weak-scaling-phi}) so that
\begin{align}
\label{e:DefM}
\Phi(r ) \ge 2 \Phi(r/M), \quad \hbox{for all } r \le \text{diam}(D).
\end{align}
 and consider three cases separately:

\noindent
(1) $\delta_D(x) \le \delta_D(y) \le |x-y|/M$:  Then in case (i), by \eqref{e:weak-scaling-psi}, we have
\begin{align}\label{e:case-i-1-lower}
I&= \int_0^{\Phi(|x-y|)}\left(1\wedge \frac{\Phi (\delta_D(x))}t \right)^{1/2} \left(1\wedge \frac{\Phi (\delta_D(y))}t \right)^{1/2} \psi(t^{-1})\, dt \nn  \\
& \ge  \int_0^{ \Phi(\delta_D(x))  } \psi(t^{-1})\, dt \asymp^{c_{15}} \psi(\Phi (\delta_D(x))^{-1}) \Phi (\delta_D(x))\, , 
\end{align}
 while in case (ii), using \eqref{e:weak-scaling-psi} and
\eqref{e:DefM} we have
 \begin{align}\label{e:case-ii-1-lower}
 I &\ge \Phi (\delta_D(x))^{1/2} \Phi (\delta_D(y))^{1/2}  \int_{\Phi(\delta_D(y))}^{\Phi(|x-y|)} t^{-1} \psi(t^{-1})\, dt \nn \\
 & \ge \Phi (\delta_D(x))^{1/2} \Phi (\delta_D(y))^{1/2}  \int_{\Phi(\delta_D(y))}^{ \Phi(M\delta_D(y))} t^{-1} \psi(t^{-1})\, dt\nn \\
 & \ge \Phi (\delta_D(x))^{1/2} \Phi (\delta_D(y))^{1/2}  \int_{\Phi(\delta_D(y))}^{ 2\Phi(\delta_D(y))} t^{-1} \psi(t^{-1})\, dt\nn \\
 & \ge c_{16} \Phi (\delta_D(x))^{1/2} \Phi (\delta_D(y))^{1/2} \psi(\Phi (\delta_D(y))^{-1})\, . 
 \end{align}

 \noindent
(2) $\delta_D(x)  \le 
|x-y| /M \le \delta_D(y)$: Then in  case (i), 
using  \eqref{e:weak-scaling-psi}, we have
\begin{align}\label{e:case-i-2-lower}
I \asymp^{c_{17}}& \int_0^{\Phi(|x-y|)}\left(1\wedge \frac{\Phi (\delta_D(x))}t \right)^{1/2} 
 \psi(t^{-1})\, dt  
 \ge I_1 \asymp^{c_1} \psi(\Phi (\delta_D(x))^{-1}) \Phi (\delta_D(x))\, ,
 \end{align}
while in case (ii), using  \eqref{e:weak-scaling-phi}, \eqref{e:weak-scaling-psi} and the fact that  $r^{1/2} \psi (r^{-1})$ is almost increasing 
near $0$, we have
\begin{align}\label{e:case-ii-2-lower}
& I \ge c_{18} I_4 = c_{18} \Phi (\delta_D(x))^{1/2} 
\int_{\Phi(\delta_D(x))}^{\Phi(|x-y|)} 
t^{-1/2}
 \psi(t^{-1})\, 
dt\nn\\
&\ge c_{18} \Phi (\delta_D(x))^{1/2} \int_{\Phi(|x-y|/M)}^{\Phi(|x-y|)} t^{-1} (t^{1/2} \psi(t^{-1}))\, 
dt\nn\\
&\ge c_{19 } \Phi (\delta_D(x))^{1/2} 
\Phi(|x-y|/M)^{1/2} \psi(\Phi(|x-y|/M)^{-1})
\int_{\Phi(|x-y|/M)}^{\Phi(|x-y|)} t^{-1} \, 
dt\nn\\
&\ge c_{20 } \log\left(\frac{\Phi(|x-y|)}{\Phi(|x-y|/M)}\right) \Phi (\delta_D(x))^{1/2} 
\Phi(|x-y|)^{1/2}
 \psi(\Phi(|x-y|)^{-1}) 
\nn\\
& \ge c_{21} (\log2) \Phi (\delta_D(x))^{1/2} \Phi(|x-y|)^{1/2}
 \psi(\Phi(|x-y|)^{-1})\, .
\end{align}
Here the last line follows from 
\eqref{e:DefM}.

\noindent
(3) $ |x-y|/M \le \delta_D(x)  \le \delta_D(y)$: Here by 
\eqref{e:weak-scaling-phi}  and \eqref{e:weak-scaling-psi}, we have
 \begin{align}\label{e:case-i-ii-3-lower}
&I \asymp^{c_{22}} \int_0^{\Phi(|x-y|)}
 \psi(t^{-1})\, dt  \asymp^{c_{23}} \Phi(|x-y|)
 \psi(\Phi(|x-y|)^{-1} ) .
 \end{align}
 
We summarize the above calculations as follows:

\noindent
\textbf{Case (i)}:  By combining \eqref{e:case-i-1}, \eqref{e:case-i-2},  \eqref{e:case-i-ii-3}, \eqref{e:case-i-1-lower}, \eqref{e:case-i-2-lower}  and  \eqref{e:case-i-ii-3-lower}, and using the fact that $r \psi (r^{-1})$ is increasing, cf.~\eqref{e:t-psi-t-1}, we get that
$$
I \asymp^{c_{24}} \Big[\Phi (\delta_D(x) \wedge \delta_D(y))\psi(\Phi (\delta_D(x) \wedge \delta_D(y))^{-1}) \Big] \wedge  \Big [\Phi(|x-y|) \psi(\Phi(|x-y|)^{-1})\Big] \, .
$$
Therefore
\begin{align*}
 \frac{1}{|x-y|^d\Phi(|x-y|)} I \asymp^{ c_{25}} \frac{ \psi(\Phi(|x-y|)^{-1}))}{|x-y|^d} \left(    \frac{ \Phi (\delta_D(x) \wedge \delta_D(y)) \psi(\Phi (\delta_D(x) \wedge \delta_D(y))^{-1}) }{ \Phi(|x-y|) \psi(\Phi(|x-y|)^{-1})) }\wedge 1\right)\, ,
\end{align*}
which together with \eqref{e:jumping-main-integral-upper} proves (i).

\noindent
\textbf{Case (ii):} By combining \eqref{e:case-ii-1}, \eqref{e:case-ii-2}, \eqref{e:case-i-ii-3}, 
\eqref{e:case-ii-1-lower}, \eqref{e:case-ii-2-lower}  and  \eqref{e:case-i-ii-3-lower} and by using \eqref{e:6.10}, we see that
\begin{align*}
&I \asymp^{ c_{26} }\left[\Phi (\delta_D(x)\wedge \delta_D(y) )^{1/2}   \wedge \Phi(|x-y|)^{1/2} \right] \\
&\ \ \ \times  \Big[\left[\Phi ( \delta_D(x) \vee \delta_D(y) )^{1/2} \psi(\Phi ( \delta_D(x)  \vee \delta_D(y) )^{-1}  )\right] \wedge   \left[\Phi(|x-y|)^{1/2} \psi(\Phi(|x-y|)^{-1}) \right] \Big]\, .
\end{align*}
Thus
 \begin{align*}
 &\frac{1}{|x-y|^d\Phi(|x-y|)} I \\
 &  \asymp^{ c_{27}} \frac{ \psi(\Phi(|x-y|)^{-1}))}{|x-y|^d} 
 \left(    \frac{ \Phi ( \delta_D(x) \wedge \delta_D(y))^{1/2}  }{ \Phi(|x-y|)^{1/2}}\wedge 1\right)
 \left(    \frac{ \eta ( \delta_D(x) \vee \delta_D(y))}{ \eta(|x-y|)} \wedge 1\right).  
 \end{align*}
 Again, together with \eqref{e:jumping-main-integral-upper}, this gives (ii).
 \qed
\begin{lemma}\label{l:comparison-boundary-terms}
\begin{itemize}
\item[\rm (i)]
Suppose that $r\mapsto r^{1/2} \psi (r^{-1})$ is almost decreasing 
near $0$. Then
there exists $C>0$ depending on the diameter of $D$ 
such that  for all $x,y\in D$,
 \begin{align}\label{e:cbt1}
\left(\frac{\Phi(\delta_D(x))^{1/2}}{\Phi(|x-y|)^{1/2}} \wedge 1\right) \left(\frac{\Phi(\delta_D(y))^{1/2}}{\Phi(|x-y|)^{1/2}} \wedge 1\right) 
\le C
\left(\frac{ \theta (\delta_D(x) \wedge \delta_D(y))}{ \theta(|x-y|)}  \wedge 1\right)\, .
 \end{align}

\item[\rm (ii)]
Suppose that $r\mapsto r^{1/2} \psi (r^{-1})$ is almost increasing 
near $0$. Then
there exists $C>0$ depending on the diameter of $D$ such that  for all $x,y\in D$,
 \begin{align}\label{e:cbt2}
\lefteqn{\left(\frac{\Phi(\delta_D(x))^{1/2}}{\Phi(|x-y|)^{1/2}} \wedge 1\right)\left(\frac{\Phi(\delta_D(y))^{1/2}}{\Phi(|x-y|)^{1/2}} \wedge 1\right) } \nn \\
&\le C
 \left(    \frac{ \Phi ( \delta_D(x) \wedge \delta_D(y))^{1/2}  }{ \Phi(|x-y|)^{1/2}}\wedge 1\right)
 \left(    \frac{ \eta ( \delta_D(x) \vee \delta_D(y)) }{ \eta(|x-y|) }\wedge 1\right)\, .
 \end{align}
\end{itemize}
\end{lemma}
\pf (i) 
Denote the left-hand side of \eqref{e:cbt1} by $I$ and the right-hand side of \eqref{e:cbt1}  by $II$.
 Without loss of generality we assume that $\delta_D(x)\le \delta_D(y)$. We consider three cases:

\noindent
(1) $\delta_D(x)\le \delta_D(y)\le |x-y|$: Then
\begin{eqnarray*}
I&=& \frac{\Phi(\delta_D(x))^{1/2}\Phi(\delta_D(y))^{1/2}}{\Phi(|x-y|)}=\frac{\Phi(\delta_D(x))\psi(\Phi(\delta_D(x))^{-1})}{\Phi(|x-y|)\psi(\Phi(|x-y|)^{-1})} \frac{\Phi(\delta_D(y))^{1/2}\psi(\Phi(|x-y|)^{-1})}{\Phi(\delta_D(x))^{1/2}\psi(\Phi(\delta_D(x))^{-1})}\\
&\le &\frac{\theta(\delta_D(x))}{\theta(|x-y|)} \, \frac{\Phi(\delta_D(y))^{1/2}\psi(\Phi(\delta_D(y))^{-1})}{\Phi(\delta_D(x))^{1/2}\psi(\Phi(\delta_D(x))^{-1})} \le c_1\frac{\theta(\delta_D(x))}{\theta(|x-y|)} =c_1 II\, , 
\end{eqnarray*}
where the last inequality follows from the assumption that $r\mapsto ^{1/2} \psi (r^{-1})$ is almost decreasing 
near $0$.

\noindent
(2) $\delta_D(x)\le |x-y|\le \delta_D(y)$: Then
\begin{eqnarray*}
I&=& \frac{\Phi(\delta_D(x))^{1/2}}{\Phi(|x-y|)^{1/2}}=\frac{\Phi(\delta_D(x))\psi(\Phi(\delta_D(x))^{-1})}{\Phi(|x-y|)\psi(\Phi(|x-y|)^{-1})} \frac{\Phi(|x-y|)^{1/2}\psi(\Phi(|x-y|)^{-1})}{\Phi(\delta_D(x))^{1/2}\psi(\Phi(\delta_D(x))^{-1})}\\
&\le &c_2 \frac{\theta(\delta_D(x))}{\theta(|x-y|)} =c_2 II\, , 
\end{eqnarray*}
again because $r\mapsto r^{1/2} \psi (r^{-1})$ is almost decreasing 
near $0$.

\noindent
(3) $|x-y|\le \delta_D(x)\le \delta_D(y)$: Then both $I$ and $II$ are equal to 1.

\medskip
\noindent
(ii) Again, 
denote  the left-hand side of \eqref{e:cbt2} by $I$ and the right-hand side of \eqref{e:cbt2}  by $II$, and assume that $\delta_D(x)\le \delta_D(y)$.

\noindent
(1) $\delta_D(x)\le \delta_D(y)\le |x-y|$: Then
\begin{eqnarray*}
I&=& \frac{\Phi(\delta_D(x))^{1/2}\Phi(\delta_D(y))^{1/2}}{\Phi(|x-y|)}\le \frac{\Phi(\delta_D(x))^{1/2}}{\Phi(|x-y|)^{1/2}}\frac{\Phi(\delta_D(y))^{1/2}}{\Phi(|x-y|)^{1/2}}  \frac{\psi(\Phi(\delta_D(y))^{-1})}{\psi(\Phi(|x-y|)^{-1})} \\
&=& \frac{\Phi(\delta_D(x))^{1/2}}{\Phi(|x-y|)^{1/2}} \frac{\eta(\delta_D(y))}{\eta(|x-y|)} = II\, .
\end{eqnarray*}

\noindent
(2) $\delta_D(x)\le |x-y| \le \delta_D(y)$: By the assumption that $r\mapsto r^{1/2}\psi(r^{-1})$ is 
almost increasing near 0, we see that $\eta$ is also almost increasing near 0. 
Therefore, 
$$
I= \frac{\Phi(\delta_D(x))^{1/2}}{\Phi(|x-y|)^{1/2}} \le c_3  \frac{\Phi(\delta_D(x))^{1/2}}{\Phi(|x-y|)^{1/2}}\frac{\eta(\delta_D(y))}{\eta(|x-y|)} =c_3 II\, .
$$

\noindent
(3) $ |x-y|\le \delta_D(x)\le \delta_D(y)$: Now $I=1$, while 
$II\ge c$, because $\eta$ is 
almost increasing near 0. 
\qed

 Recall that
$
j(r)={\psi(\Phi(r)^{-1})}r^{-d}$ for  $r>0$.
\begin{thm}\label{t:jumping-function-estimate}
\begin{itemize}
\item[\rm (i)]
Suppose that $r\mapsto r^{1/2} \psi (r^{-1})$ is almost decreasing 
near $0$ and  
for each $T>0$ there is a constant  $C_4=C_4(T, \psi)>0$  such that
\eqref{e:6.9} holds true. 
Then
there exists $C \ge 1$ such that  for all $x,y\in D$,
$$
J^{Y^D}(x, y)\asymp^C
\left(\frac{ \theta (\delta_D(x) \wedge \delta_D(y))}{ \theta(|x-y|)}  \wedge 1\right) j(|x-y|)\, .
$$

\item[\rm (ii)]
Suppose that $r\mapsto r^{1/2} \psi (r^{-1})$ is almost increasing 
near $0$ and  
for every $T>0$, there is a constant  $C_5=C_5(T, \psi)>0$  so that
\eqref{e:6.10} holds true. 
Then there exists $C \ge 1$ such that  for all $x,y\in D$,
$$
J^{Y^D}(x, y)\asymp^C
 \left(    \frac{ \Phi ( \delta_D(x) \wedge \delta_D(y))^{1/2}  }{ \Phi(|x-y|)^{1/2}}\wedge 1\right)
 \left(    \frac{ \eta ( \delta_D(x) \vee \delta_D(y))}{ \eta(|x-y|) }\wedge 1\right) j(|x-y|)\, .
$$
\end{itemize}
\end{thm}
\pf 
Let $T:=2\Phi(\mathrm{diam}(D))$. 
By \eqref{e:JY}, \eqref{e:rtxy} and \eqref{e:p-r-comparable}, we have
\begin{align*}
&J^{Y^D}(x,y)=\int_0^{\infty}p^D(t,x,y)\nu(t)\, dt \\
&\le  c_1 \left(\int_0^{\Phi(|x-y|)} r(t,x,y)\nu(t)\, dt
 +\int_{\Phi(|x-y|)}^Tr(t,x,y)\nu(t)\, dt +\int_T^{\infty} p^D(t,x,y)v(t)\, dt\right) \\
&=:c_1(I_1+I_2+I_3)\, .
\end{align*}

\noindent
(i) By Lemma \ref{l:jumping-main-integral}(i), we get
$$
I_1\le c_2 \left(\frac{ \theta (\delta_D(x) \wedge \delta_D(y))}{ \theta(|x-y|)}  \wedge 1\right) j(|x-y|)\, .
$$
By using Lemma \ref{l:second-integral} and the fact that $\nu$ is decreasing in the first line, and 
Lemma \ref{l:comparison-boundary-terms}(i), 
the definition of $j$ and \eqref{e:v-nu-asymp}  in the second line,
\begin{eqnarray*}
I_2&\le &c_3\left(\frac{\Phi(\delta_D(x))^{1/2}}{\Phi(|x-y|)^{1/2}}\wedge 1\right) \left(\frac{\Phi(\delta_D(y))^{1/2}}{\Phi(|x-y|)^{1/2}}\wedge 1\right)
\frac{\Phi(|x-y|)\nu(\Phi(|x-y|)}{|x-y|^d}\\
&\le &c_4 \left(\frac{ \theta (\delta_D(x) \wedge \delta_D(y))}{ \theta(|x-y|)}  \wedge 1\right) j(|x-y|)\, .
\end{eqnarray*}
By using Lemma \ref{l:third-integral} in the first line, and again Lemma \ref{l:comparison-boundary-terms}(i), the definition of $j$ and \eqref{e:v-nu-asymp}  in the second line, we have
\begin{eqnarray*}
I_3 &\le & c_5\left(\frac{\Phi(\delta_D(x))^{1/2}}{\Phi(|x-y|)^{1/2}}\wedge 1\right) \left(\frac{\Phi(\delta_D(y))^{1/2}}{\Phi(|x-y|)^{1/2}}\wedge 1\right)
j(|x-y)\\
&\le &c_6 \left(\frac{ \theta (\delta_D(x) \wedge \delta_D(y))}{ \theta(|x-y|)}  \wedge 1\right) j(|x-y|)\, .
\end{eqnarray*}
 The three displays above prove  the upper bound.
For the lower bound, it suffices to use the lower bound for $I_1$ coming from Lemma  \ref{l:jumping-main-integral}(i).

\noindent
(ii) This is proved in the same way as (i) by using part (ii) of
 Lemmas  \ref{l:jumping-main-integral} and \ref{l:comparison-boundary-terms}.
\qed
\begin{example}\label{e:stable}{\rm
Assume that $\psi(\lambda)=\lambda^{\gamma}$ where $\gamma\in (0,1)$. 
As already mentioned in
Remark \ref{r:moved},
when $\gamma>1/2$ the assumptions in (i) hold true. Since $\theta(t)=\Phi(t)^{1-\gamma}$ we get that
$$
J^{Y^D}(x,y)\asymp^c \left(\frac{   \Phi(\delta_D(x) \wedge \delta_D(y))}{ \Phi(|x-y|)}  \wedge 1\right)^{1-\gamma} \frac{\Phi(|x-y|)^{-\gamma}}{|x-y|^d}\, .
$$
When $\gamma <1/2$, the assumptions in (ii) hold true, $\eta(t)=\Phi(t)^{1/2-\gamma}$ and
\begin{align}\label{e:jgamma}
J^{Y^D}(x, y)\asymp^c
 \left(    \frac{ \Phi ( \delta_D(x) \wedge \delta_D(y))  }{ \Phi(|x-y|)}\wedge 1\right)^{1/2}
  \left(    \frac{ \Phi ( \delta_D(x) \vee \delta_D(y)) }{ \Phi(|x-y|) }\wedge 1\right)^{1/2-\gamma} 
 \frac{\Phi(|x-y|)^{-\gamma}}{|x-y|^d}\, .
\end{align}
The case $\gamma=1/2$ is not covered by Theorem \ref{t:jumping-function-estimate}, but 
by following  
the proofs of Lemmas \ref{l:jumping-main-integral} and  
\ref{l:comparison-boundary-terms} step by step, it is straightforward to deduce that
\begin{align}
&J^{Y^D}(x, y)\nn\\
\asymp^c& \left(    \frac{ \Phi ( \delta_D(x) \wedge \delta_D(y))  }{ \Phi(|x-y|)}\wedge 1\right)^{1/2} 
\log\left(1+\frac{\Phi(\delta_D(x)\vee \delta_D(y))\wedge \Phi(|x-y|)}{\Phi(\delta_D(x)\wedge \delta_D(y))\wedge\Phi (|x-y|)}\right)
\frac{\Phi(|x-y|)^{-1/2}}{|x-y|^d}\, .\label{e:Jgamma12}
\end{align}
In particular, with $y\in D$ fixed, as $\delta_D(x)\to 0$, we have
$$
J^{Y^D}(x, y)\asymp^c \begin{cases}
\Phi(\delta_D(x))^{1/2}, & 0<\gamma <1/2,\\
\Phi(\delta_D(x))^{1/2}\,  \log(1/ \Phi(\delta_D(x))), & \gamma=1/2,\\
\Phi(\delta_D(x))^{1/2}\, \Phi(\delta_D(x))^{1/2-\gamma}, & 1/2<\gamma <1.
\end{cases}
$$
}
\end{example}

\section{Failure of BHP in the case  of $\gamma_2 \le 1/2$}\label{s:counterexample}
In this section 
we assume that $d\ge 2$ and  that $D$ is a bounded
$C^{1,1}$ open set in  $\R^d$ 
with $C^{1, 1}$ characteristics $(R, \Lambda)$ with $R <1$.
The goal of this section is  
to give an example showing that even the non-scale invariant
boundary Harnack principle does not hold when $\gamma_2 \le 1/2$.
For simplicity we  consider the case $\psi(t)=t^\gamma$ and
$\gamma_2=\gamma \le 1/2$. 
The example works for any bounded $C^{1, 1}$ open set.

 Suppose that the non-scale invariant boundary Harnack principle holds near the boundary of $D$. That is, 
{\it 
there is a constant $\wh{R}\in (0,1)$ 
such that for any $r \in (0, \wh{R}\, ]$, there exists 
a constant $c=c(r)\ge 1$ such that for every $Q\in \partial D$ and any non-negative functions $f, g$ in $D$ which are harmonic  in $D \cap B(Q, r)$ with respect to $Y^D$ and 
vanish continuously on $ \partial D \cap B(Q, r)$, we have}
\begin{equation}\label{e:bhp_mfg}
\frac{f(x)}{f(y)}\,\le c\,\frac{g(x)}{g(y)} \qquad
\hbox{for all } x, y\in  D \cap B(Q, r/2).
\end{equation}
Note that we can take $g(\cdot)=G^{Y^D}(\cdot, w)$ with $w\notin D \cap B(Q, r)$.
Thus
by Theorem \ref{t:green-function-estimate},
we have that for any $r \in (0, \wh{R}\, ]$ there exists a constant $C_6=C_6(r)>0$ such that for every $Q\in \partial D$ and any non-negative function $f$ in $D$ which is harmonic  in $D \cap B(Q, r)$ with respect to $Y^D$ and vanishes continuously on $ \partial D \cap B(Q, r)$, 
\begin{equation}\label{e:bhp_mfx}
\frac{f(x)}{f(y)}\,\le \,
C_6 \,\frac{\Phi(\delta_D(x))^{1/2}}{\Phi(\delta_D(y))^{1/2}}, \quad  \hbox{for all } x, y\in  D \cap B(Q, r/2).
\end{equation}

Fix $Q \in \partial D$.
We choose a
$C^{1,1}$ function $\varphi: \bR^{d-1}\to \bR$ satisfying $\varphi
(\wt 0)= 0$, $\nabla\varphi (\wt 0)=(0, \dots, 0)$, $\| \nabla \varphi
 \|_\infty \leq \Lambda$, $| \nabla \varphi (\wt y)-\nabla \varphi (\wt z)|
\leq \Lambda |\wt y-\wt z|$, and an orthonormal coordinate system $CS$ with
its origin at $Q$ such that
$$
B(Q, R)\cap D=\{ y=(\wt y, y_d) \in B(0, R) \mbox{ in } CS: y_d >
\varphi (\wt y) \}.
$$

Recall that $\kappa_0=(1+(1+\Lambda)^2)^{-1/2}$.
Since $D$ satisfies the interior ball condition, there exist
$r_0 \le \wh{R}\wedge (2^{-4}\kappa_0 R)$ and $x^{(1)} \in B(Q, R)\cap D$ with $\delta_D(x^{(1)})=r_0$ such that 
$\delta_D(x^{(s)})=|x^{(s)}-Q|$ for all $s \le 1$
where $x^{(s)}=Q+s(x^{(1)}-Q)$.
 
 In the coordinate system $CS$ we have $\wt x^{(s)} = \wt 0$ and $x^{(1)}=(\wt
0, r_0)$. For any $b_1, b_2>0$, we define
$$
D^*(b_1, b_2):=\left\{ y=(\wt y, y_d) \mbox{ in } CS: 0<y_d-\varphi(\wt
y)< 2^{-2}\kappa_0 r_0 b_1, \ |\wt y| <  2^{-2}\kappa_0 r_0 b_2\right\}.
$$
By \eqref{e:Uzr}, 
we have that
$D^*(2, 2)\subset D\cap B(Q, r_0/2)$.
Recall that $V_{Q}(2^{-2}\kappa_0 r_0)$ is a $C^{1,1}$
domain  with $C^{1,1}$ characteristics $(r_0 R/L, \Lambda L/r_0)$
such that $D^*( 1/2, 1/2) \subset V_{Q}(2^{-2}\kappa_0 r_0)  \subset  D^*( 1, 1) $,
where $L=L(R, \Lambda, d)>0$.
Let $V=V_{Q}(2^{-2}\kappa_0 r_0)$ and $U=D^*(2,2)$.

Recall that 
\begin{align}\label{e:ggamma}
g(r)=\frac{1}{r^d \psi(\Phi(r)^{-1})}=\frac{\Phi(r)^\gamma}{r^d}.
\end{align}
\begin{lemma}\label{l:Elow} 
If $\gamma \le 1/2$, then there exists $C>0$ such that 
$$
\E_x \left[\int_0^{\tau_V}\Phi(\delta_D(Y^D_t))^{\frac12-\gamma}dt \right] \ge 
C \Phi(\delta_D(x))^{1/2} \log(r_0/\delta_D(x)) 
$$
for all  $x=x^{(s)}=(\wt 0, s)$ in   $CS$ with $s \in (0, 2^{-7} \kappa_0 r_0)$.
\end{lemma}
\pf
By 
Proposition \ref{p:gfcnlb-scaling},
\begin{align}\label{e:BL}&\E_x[\int_0^{\tau_V}\Phi(\delta_D(Y^D_t))^{\frac12-\gamma}dt]=\int_{V } G^{Y^D}_{V}(x,z) \Phi(\delta_D(z))^{\frac12-\gamma} dz
 \nn\\
&\ge \int_{V } G^{Y^V} (x,z)  \Phi(\delta_D(z))^{\frac12-\gamma}   dz
\nn\\
&\ge c_1\int_{V }  \Phi(\delta_D(z))^{\frac12-\gamma} \left(\frac{\Phi(\delta_{V}(x))^{1/2}}{\Phi(|x-z|)^{1/2}}\wedge 1\right) \left(\frac{\Phi(\delta_{V}(z))^{1/2}}{\Phi(|x-z|)^{1/2}}\wedge 1\right) g(|x-z|) dz.
\end{align} 

Recall from \eqref{e:Ceta} that 
there exists $\varepsilon \in (0, (2(1+ \Lambda))^{-2})$ such that 
the cone 
$$
\sC_\varepsilon(x):=
\{y=(\wt y,y_d)\in B(x,2^{-6} \kappa_0 r_0 ) \mbox{ in }  CS : y_d>x_d ,|\wt y|<\varepsilon (y_d -x_d)\}
$$ is contained in $D^*(2^{-2}, 2^{-2})$.
Moreover,  \eqref{e:D2_421} implies
$$
\left(\frac{\Phi(\delta_{V}(z))}{\Phi(|x-z|)}\wedge 1 \right) 
\asymp 1, \quad z\in \sC_{\varepsilon/2}(x).
$$
If $z \in \sC_{\varepsilon/2}(x) \setminus  B(x,   \delta_D(x)/2)$,
then $|x-z| \ge \delta_D(x)/2$, 
so by \eqref{e:D2_421} and \eqref{e:Phi(lambda-t)},
\begin{align*}
\left(\frac{\Phi(\delta_{D}(x))}{\Phi(|x-z|)}\wedge 1 \right) 
\Phi(\delta_{D}(z))^{1-2\gamma}  \ge 
c_2\left(\frac{\Phi(\delta_{D}(x))}{\Phi(|x-z|)}\wedge 1\right)\Phi(|x-z|)^{1-2\gamma}\ge c_3 \frac{
\Phi(\delta_{D}(x))}{\Phi(|x-z|)^{2\gamma}} \, .
\end{align*}
Therefore, using \eqref{e:ggamma},
\begin{align*}&
\int_{V }  \Phi(\delta_D(z))^{\frac12-\gamma} \left(\frac{\Phi(\delta_{V}(x))^{1/2}}{\Phi(|x-z|)^{1/2}}\wedge 1\right) \left(\frac{\Phi(\delta_{V}(z))^{1/2}}{\Phi(|x-z|)^{1/2}}\wedge 1\right) g(|x-z|) dz\\
&\ge  c_{4} \Phi(\delta_D(x))^{1/2}  
\int_{\sC_{\varepsilon/2}(x) \setminus  B(x,   \delta_D(x)/2)}
\frac{g(|x-z|) }{\Phi(|x-z|)^{\gamma}}dz
  \nn \\
&  \ge  c_{5}\Phi(\delta_D(x))^{1/2} \int^{2^{-6} \kappa_0 r_0}_{\delta_D(x)/2}
\frac1{s}ds\ge c_{6}
 \Phi(\delta_D(x))^{1/2} \log(r_0/\delta_D(x)).
\end{align*} 
\qed

Choose a point $z_0 \in \partial D \setminus \overline{D \cap B(Q, 2r_0)}$ 
with $|z_0-Q| \le 1$
(such $z_0$ exists since $2r_0 \le 2^{-3}\kappa_0R$).
  For  $n \in \N$ large enough so that $B(z_0,1/n)$ does not 
intersect $B(Q, 2r_0)$, we define
$$
f_n(y):=\Phi(\delta_D(y))^{-1/2} {\bf 1}_{D \cap B(z_0, 1/n)}(y) \times
\begin{cases}
|D \cap B(z_0, 1/n)|^{-1}, & \text{for } \gamma <1/2 \\
K_n^{-1}, & \text{for } \gamma =1/2,
\end{cases}
$$
where
$$
K_n:=\int_{D \cap B(z_0, 1/n)} \log(\frac{1}{\Phi(\delta_D(y))})dy.
$$ 
Define
$$
g_n(x):=\E_x[f_n(Y^D_{\tau_V})]= \int_{D\setminus U}  \int_{V} G^{Y^D}_{V}(x,z) J^{Y^D}(z,y)f_n(y)dz dy, \quad x \in V.
$$
\begin{lemma}\label{l:gnlow} 
If $\gamma \le 1/2$,  then there exists $C_7>0$ such that 
\begin{align}
\label{e:gnlow}
\liminf_{n \to \infty}g_n(x) 
\ge C_7 \Phi(\delta_D(x))^{1/2} \log(r_0/\delta_D(x))
\end{align}
for all  $x=x^{(s)}=(\wt 0, s)$   in   $CS$ with $s \in (0, 2^{-7} \kappa_0 r_0)$.
\end{lemma}
\pf 
{\it (i) Case $\gamma <1/2$:}  
Since 
$$
1\asymp^{c} |z-y|  \ge  c_1 \delta_D(z) \quad \text{and} \quad |y-z| \ge c_2  \delta_D(y), \quad 
\text{ for } (y, z)  \in (D \cap B(z_0, 1/n)) \times   V,
$$
using \eqref{e:jgamma} we have 
$$
J^{Y^D}(z, y)\asymp^c  \Phi ( \delta_D(y))^{1/2} \Phi ( \delta_D(z))^{1/2}\Phi ( \delta_D(z) \vee \delta_D(y))^{-\gamma}     \quad 
\text{ for } (y, z)  \in (D \cap B(z_0, 1/n)) \times   V.
$$
Thus, 
\begin{align}\label{e:notD213}
&g_n(x)= \int_{D\setminus U}  \int_{V} G^{Y^D}_{V}(x,z) J^{Y^D}(z,y)f_n(y)dz dy
\nn \\
&\asymp^c   |D \cap B(z_0, 1/n)|^{-1}   \int_{V} \int_{D \cap B(z_0, 1/n)} \Phi ( \delta_D(z) \vee \delta_D(y))^{-\gamma}     dy \Phi ( \delta_D(z))^{1/2}G^{Y^D}_{V}(x,z) dz.
\end{align} 
Since
$$
\lim_{n \to \infty}
  |D \cap B(z_0, 1/n)|^{-1}
 \int_{D \cap B(z_0, 1/n)} \Phi ( \delta_D(z) \vee \delta_D(y))^{-\gamma}     dy =
 \Phi ( \delta_D(z))^{-\gamma}, \quad z \in V,     
 $$
 $$
  |D \cap B(z_0, 1/n)|^{-1}
 \int_{D \cap B(z_0, 1/n)} \Phi ( \delta_D(z) \vee \delta_D(y))^{-\gamma}     dy \le    \Phi ( \delta_D(z))^{-\gamma},     \quad z \in V
 $$ 
 and
 $$
 \int_{V} \Phi ( \delta_D(z))^{\frac12-\gamma}   G^{Y^D}_{V}(x,z) dz  < \infty, $$
 by the dominated convergence theorem,  for all  $x=x^{(s)}=(\wt 0, s)$   in   
 $CS$ with $s \in (0, 2^{-7} \kappa_0 r_0)$,
 \begin{align}\label{e:notD231}
&\lim_{n \to \infty}  |D \cap B(z_0, 1/n)|^{-1}   \int_{V} \int_{D \cap B(z_0, 1/n)} \Phi ( \delta_D(z) \vee \delta_D(y))^{-\gamma}     dy \Phi ( \delta_D(z))^{1/2}G^{Y^D}_{V}(x,z) dz \nn \\
&=  \int_{V} \Phi ( \delta_D(z))^{\frac12-\gamma}   G^{Y^D}_{V}(x,z) 
dz.
\end{align} 
Combining
\eqref{e:notD231} with Lemma \ref{l:Elow} we conclude that \eqref{e:gnlow} holds true for $\gamma<1/2$.

\noindent
{\it (ii) Case $\gamma =1/2$:}  
Using \eqref{e:Jgamma12} and following the same argument in (i), we have
\begin{align}\label{e:notD213h}
&g_n(x)= \int_{D\setminus U}  \int_{V} G^{Y^D}_{V}(x,z) J^{Y^D}(z,y)f_n(y)dz dy
\nn \\
\asymp^c &  K_n^{-1}   \int_{V} \int_{D \cap B(z_0, 1/n)} 
\left(\frac{\Phi ( \delta_D(z) \wedge \delta_D(y))}
{\Phi ( \delta_D(y))}
\right)^{1/2}    \log\left(1+\frac{\Phi(\delta_D(y)\vee \delta_D(z))}{\Phi(\delta_D(y)\wedge \delta_D(z))}\right)
 dy G^{Y^D}_{V}(x,z) dz.
\end{align} 
We have 
$
 \int_{V}  G^{Y^D}_{V}(x,z) dz  < \infty$ and 
\begin{align*}
&  K_n^{-1}
 \int_{D \cap B(z_0, 1/n)}\left(\frac{\Phi ( \delta_D(z) \wedge \delta_D(y))}
{\Phi ( \delta_D(y))}
\right)^{1/2}    \log\left(1+\frac{\Phi(\delta_D(y)\vee \delta_D(z))}{\Phi(\delta_D(y)\wedge \delta_D(z))}\right)
 dy  \nn\\
\le & K_n^{-1}
 \int_{D \cap B(z_0, 1/n)\cap \{\delta_D(z) \le \delta_D(y)\}}\left(\frac{\Phi ( \delta_D(z))}
{\Phi ( \delta_D(y))}
\right)^{1/2}    \log\left(1+\frac{\Phi(\delta_D(y))}{\Phi( \delta_D(z))}\right)
 dy \nn\\
&+ K_n^{-1}
 \int_{D \cap B(z_0, 1/n)\cap \{\delta_D(z) > \delta_D(y)\}} \log\left(1+\frac{\Phi( \delta_D(z))}{\Phi(\delta_D(y))}\right)
 dy \nn\\
 \le& c K_n^{-1}
 \int_{D \cap B(z_0, 1/n)}   \log\left(\frac{1}{\Phi(\delta_D(y))}\right)
 dy   \le c,     \quad z \in V.
\end{align*} 
Moreover, since 
 $$
\limsup_{n \to \infty}
  K_n^{-1}
 \int_{D \cap B(z_0, 1/n)}   \log\left(1+\frac{\Phi(\delta_D(y))}{\Phi(\delta_D(z))}\right)dy\le
  \limsup_{n \to \infty}
  K_n^{-1}|D \cap B(z_0, 1/n) |  =0, \quad z \in V,
 $$
for each $z \in V$,
\begin{align*}
&\lim_{n \to \infty}
  K_n^{-1}
 \int_{D \cap B(z_0, 1/n)}\left(\frac{\Phi ( \delta_D(z) \wedge \delta_D(y))}
{\Phi ( \delta_D(y))}
\right)^{1/2}    \log\left(1+\frac{\Phi(\delta_D(y)\vee \delta_D(z))}{\Phi(\delta_D(y)\wedge \delta_D(z))}\right)
 dy\\   
 = &\lim_{n \to \infty}
  K_n^{-1}
 \int_{D \cap B(z_0, 1/n)}    \log\left(1+\frac{\Phi(\delta_D(z))}{\Phi(\delta_D(y))}\right)
 dy\\
  = &\lim_{n \to \infty}
  K_n^{-1}
 \int_{D \cap B(z_0, 1/n)}   \log\left(1+\frac{\Phi(\delta_D(y))}{\Phi(\delta_D(z))}\right)
 dy\\
&+ \lim_{n \to \infty}
  K_n^{-1}
 \int_{D \cap B(z_0, 1/n)}   \left( \log\frac{1}{\Phi(\delta_D(y))} -  \log\frac{1}{\Phi(\delta_D(z))}    \right) dy
  \\=&1- \log\frac{1}{\Phi(\delta_D(z))} \lim_{n \to \infty}  K_n^{-1}|D \cap B(z_0, 1/n)|
 =1.
\end{align*}
Thus by the dominated convergence theorem,  for all  $x=x^{(s)}=(\wt 0, s)$   
in  $CS$ with $s \in (0, 2^{-7} \kappa_0 r_0)$,
 \begin{align}\label{e:notD231h}
&\lim_{n \to \infty} K_n^{-1}   \int_{V} \int_{D \cap B(z_0, 1/n)} 
\left(\frac{\Phi ( \delta_D(z) \wedge \delta_D(y))}
{\Phi ( \delta_D(y))}
\right)^{1/2}    \log\left(1+\frac{\Phi(\delta_D(y)\vee \delta_D(z))}{\Phi(\delta_D(y)\wedge \delta_D(z))}\right)
 dy G^{Y^D}_{V}(x,z) dz \nn \\
&=  \int_{V}   G^{Y^D}_{V}(x,z) dz.
\end{align} 
Combining
\eqref{e:notD231h} with Lemma \ref{l:Elow} we conclude that \eqref{e:gnlow} holds true for $\gamma=1/2$.\

\qed

\bigskip

By Lemmas \ref{l:gamma-small} and \ref{l:gamma-1/2}, we have that for large $n$ anf $y \in D \cap B(Q, 2^{-7} \kappa_0 r_0)$
 \begin{align}
 \label{e:gnupper}
 g_n(y) &\le c_1  \int_{V} \Phi ( \delta_D(z))^{\frac12-\gamma}   G^{Y^D}_{V}(y,z) dz 
\nn\\
& \le c_2 \Phi ( r_0)^{\frac12-\gamma}  \int_{D} G^{Y^D}(y,z) dz \nn\\
& \le c_3  \Phi ( r_0)^{\frac12-\gamma} \Phi( \delta_D(y))^{\gamma}
 \begin{cases}
1,
 & \text{for } \gamma <1/2 \\
\log(1/\delta_D(y)), & \text{for } 
\gamma =1/2.
\end{cases}
 \end{align}
Thus  we see that $g_n$'s are non-negative functions 
in $D$ which are harmonic  in $D \cap B(Q,2^{-7} \kappa_0 r_0)$ with respect to $Y^D$ and vanish continuously on $ \partial D \cap B(Q,2^{-7} \kappa_0 r_0)$.
Therefore, by \eqref{e:bhp_mfx},
$$
\frac{g_n(y)}{g_n(w)} \le C_6\frac{\Phi(\delta_D(y))^{1/2} }{\Phi(\delta_D(w))^{1/2} } \quad \text{for all } y\in D \cap B(Q,2^{-8} \kappa_0 r_0)
 $$
  where $w=(\wt{0}, 2^{-9}\kappa_0 r_0)$ and $C_6=C_6(2^{-7} \kappa_0 r_0)$.
 Thus by \eqref{e:gnupper}, for all $y\in D \cap B(Q,2^{-8} \kappa_0 r_0)$,
 $$
 \limsup_{n \to \infty}g_n(y) \le 
 C_6 \limsup_{n \to \infty}g_n(w) \frac{\Phi(\delta_D(y))^{1/2} }{\Phi(\delta_D(w))^{1/2} } \le c_4\log(c_5/r_0) \Phi(\delta_D(y))^{1/2}.$$
 This and \eqref{e:gnlow} imply that  for all  $x=x^{(s)}=(\wt 0, s)$  in   $CS$ with $s \in (0, 
  2^{-8} \kappa_0 r_0)$, 
  $$
 \log(r_0/\delta_D(x)) \le  (c_4/C_7) \log(c_5/r_0), $$
which  gives a contradiction.

\bigskip
\noindent
{\bf Acknowledgements:}
We thank the referee for carefully reading the
manuscript and providing some useful suggestions.

\bigskip
\noindent

\begin{singlespace}

%%%%%%%%%%%%%%%%%%%%%%%%%%%%%%%%%
%%%%%%%%%%%%%%%%%%%%%%%%%%%%%%%%%%
%%%%%%%%            References
%%%%%%%%%%%%%%%%%%%%%%%%%%%%%%%%%
%%%%%%%%%%%%%%%%%%%%%%%%%%%%%%%%%%

\end{singlespace}

\end{doublespace}
\vskip 0.1truein

\parindent=0em

{\bf Panki Kim}

Department of Mathematical Sciences and Research Institute of Mathematics,
Seoul National University, Building 27, 1 Gwanak-ro, Gwanak-gu Seoul 151-747, Republic of Korea

E-mail: \texttt{pkim@snu.ac.kr}

\bigskip

{\bf Renming Song}

Department of Mathematics, University of Illinois, Urbana, IL 61801,
USA, and\\
School of Mathematical Sciences, Nankai University, Tianjin 300071, PR China

E-mail: \texttt{rsong@illinois.edu}

\bigskip

{\bf Zoran Vondra\v{c}ek}

Department of Mathematics, Faculty of Science, University of Zagreb, Zagreb, Croatia

E-mail: \texttt{vondra@math.hr}
\end{document}